\documentclass[a4paper,12pt]{article}
\usepackage{amssymb}
\usepackage{amsthm}
\usepackage{tikz}
\usepackage{tikz,pgfplots}
\usetikzlibrary{decorations.markings}
\usepackage{amsmath,amssymb}
\usepackage{todonotes}
\usepackage{url}
\usepackage[T1]{fontenc}
\usepackage{cite}
\usepackage{color}
\usepackage{graphicx}
\usepackage{placeins,caption}
\usepackage[figure,tworuled,linesnumbered,vlined]{algorithm2e} 

\DeclareMathOperator{\Circ}{Circ}
\DeclareMathOperator{\mob}{Mob_{gp}}
\DeclareMathOperator{\gp}{gp}

\DeclareMathOperator{\gpe}{gp_e}
\DeclareMathOperator{\vp}{vp}
\DeclareMathOperator{\g}{g}
\DeclareMathOperator{\s}{s}
\DeclareMathOperator{\cut}{c}
\DeclareMathOperator{\I}{I}
\DeclareMathOperator{\ip}{ip}
\DeclareMathOperator{\eq}{eq}
\DeclareMathOperator{\gex}{gex}
\DeclareMathOperator{\mex}{mex}
\DeclareMathOperator{\rad}{rad}
\DeclareMathOperator{\diam}{diam}
\DeclareMathOperator {\mono} {mp}
\DeclareMathOperator {\Pl} {P\ell}
\DeclareMathOperator {\sgp} {sgp}
\DeclareMathOperator{\cp}{\,\square\,}
\newcommand{\GG}{G\overline{G}}
\newcommand{\gps}{{\rm gp}{_{\rm S}}}
\newcommand{\gpsl}{\underline{{\rm gp}}{_{\rm S}}}
\DeclareMathOperator{\gpt}{gp_t}
\DeclareMathOperator{\gpd}{gp_d}
\DeclareMathOperator{\gpo}{gp_o}
\DeclareMathOperator{\sr}{_{SR}}
\DeclareMathOperator{\dpn}{dp}
 \newcommand{\sink }{{\rm Sink}}
\newcommand{\nsink }{{\rm Sink'}}
\newcommand{\source }{{\rm Source}}
\newcommand{\nsource }{{\rm Source'}}
\DeclareMathOperator{\Ka}{Ka}
\DeclareMathOperator{\Pe}{P}
\newcommand{\ugp}{{\rm gp^{\rightarrow }}}
\newcommand{\lgp}{{\rm gp_{\rightarrow }}}

\usepackage[top=3cm,bottom=3cm,left=2.8cm,right=2.8cm]{geometry}
	\newtheorem{theorem}{Theorem}[section]
	\newtheorem{lemma}[theorem]{Lemma}
	\newtheorem{corollary}[theorem]{Corollary}
	\newtheorem{proposition}[theorem]{Proposition}
	
	\newtheorem{conjecture}[theorem]{Conjecture}

	\theoremstyle{definition}

	\newtheorem{definition}[theorem]{Definition}

	\newcommand{\move}{\rightsquigarrow}

	\newcommand{\address}[1]{#1}

\pgfplotsset{compat=1.18} 

\begin{document}
	
		\title{The General Position Problem: A Survey}
		\author{ Ullas Chandran S.V. $^{a}$ \\ \texttt{\footnotesize svuc.math@gmail.com}        
        \and Sandi Klav\v{z}ar $^{b,c,d}$ \\ \texttt{\footnotesize sandi.klavzar@fmf.uni-lj.si}
		\and
		James Tuite $^{e,f}$ \\ \texttt{\footnotesize james.t.tuite@open.ac.uk}
  }

\maketitle	

\address{
\noindent
$^a$  Department of Mathematics, Mahatma Gandhi College,   Thiruvananthapuram-695004, \newline Kerala, India \\
    $^b$ Faculty of Mathematics and Physics, University of Ljubljana, Slovenia\\
	$^c$ Institute of Mathematics, Physics and Mechanics, Ljubljana, Slovenia\\
	$^d$ Faculty of Natural Sciences and Mathematics, University of Maribor, Slovenia\\
	$^e$ School of Mathematics and Statistics, Open University, Milton Keynes, UK\\
    $f$ Department of Informatics and Statistics, Klaip\.{e}da University, Lithuania\\
    
}

\begin{abstract}
Inspired by a chessboard puzzle of Dudeney, the general position problem in graph theory asks for a largest set $S$ of vertices in a graph such that no three elements of $S$ lie on a common shortest path. The number of vertices in such a largest set is the \emph{general position number} of the graph. This paper provides a survey of this rapidly growing problem, which now has an extensive literature. We cover exact results for various graph classes and the behaviour of the general position number under graph products and operations. We also discuss interesting variations of the general position problem, including those corresponding to different graph convexities, as well as dynamic, fractional, colouring and game versions of the problem. 
\end{abstract}


\noindent
AMS Subj.\ Class.\ (2020): 05C12, 05C69, 68Q25
\newpage
\tableofcontents

\section{Introduction}\label{sec:intro}

Many significant mathematical problems involve several objects lying on a line, or else how to avoid such configurations. For example, a collection of vectors is in general position if no three lie on a common line, so that such problems relate to the fundamental notions of dimension and independence. The Hales-Jewett Theorem, Van der Waerden's Theorem and other results from Ramsey theory on the integers and combinatorial set theory~\cite{RamseyTheoryontheIntegers} can be interpreted as no-$k$-in-line problems, and the study of points with no $k$ on a line has a long history in discrete geometry.

No-three-in-line problems are also familiar from games that span multiple cultures from the start of recorded history to the modern day. The ancient Egyptians and Romans played games in which the winner had to get three in a line~\cite{Zaslavsky}, leading up to our modern day games Tic-Tac-Toe (known to the third author in the UK as `Noughts and Crosses') and `Nine Men's Morris' (also known as `Mill'). For Tic-Tac-Toe and Nine Men's Morris, optimal play leads to a draw~\cite{Gasser}, i.e.\ a configuration with no three counters in line to which no further counters can be added. The game `Connect Four' requires a player to get four counters in a line to win, whereas in the Japanese game Gomoku the winner must achieve five in a row (the former game is a first player win~\cite{connectfour} and the latter is a first player win on a $15 \times 15$ board~\cite{Go}). We could mention many other examples. Some mathematical discussion of these games can be found in~\cite{positionalgames}. 

The mathematical investigation of such problems in combinatorics can be dated to a puzzle of the famous recreational mathematician Henry Dudeney (1857-1930). As a keen chess player, many of his puzzles are set on a chessboard, and one such problem from 1900 (see~\cite{dudeney-1917}) can be regarded as the origin of the general position problem. The puzzle asks for the largest number of pawns that can be placed on an $8 \times 8$ chessboard without three pawns lying on any straight line in the plane (not just rows, columns or diagonals of the board). The answer is 16. The more general No-Three-In-Line Problem for an $n \times n$ chessboard has been called by Brass, Moser and Pach ``one of the oldest and most extensively studied geometric questions concerning lattice points''~\cite{Brass}. It is listed as Problem 72 in Ben Green's list of 100 unsolved problems~\cite{Green}.

\begin{center}
			\begin{figure}[ht!]
				\centering
				\includegraphics[scale=0.35]{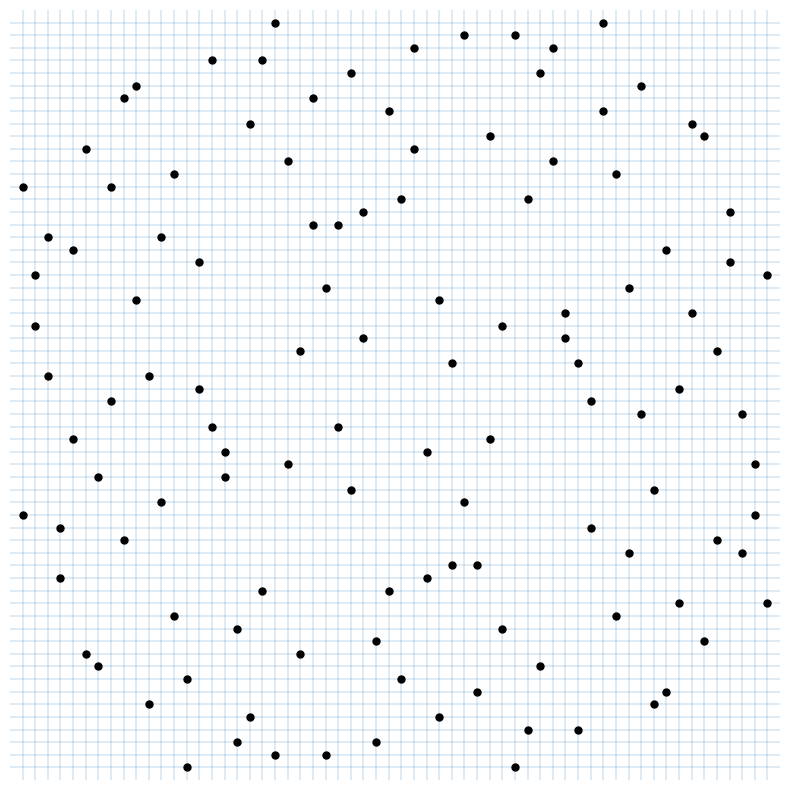}
                \quad
                \includegraphics[scale=0.35]{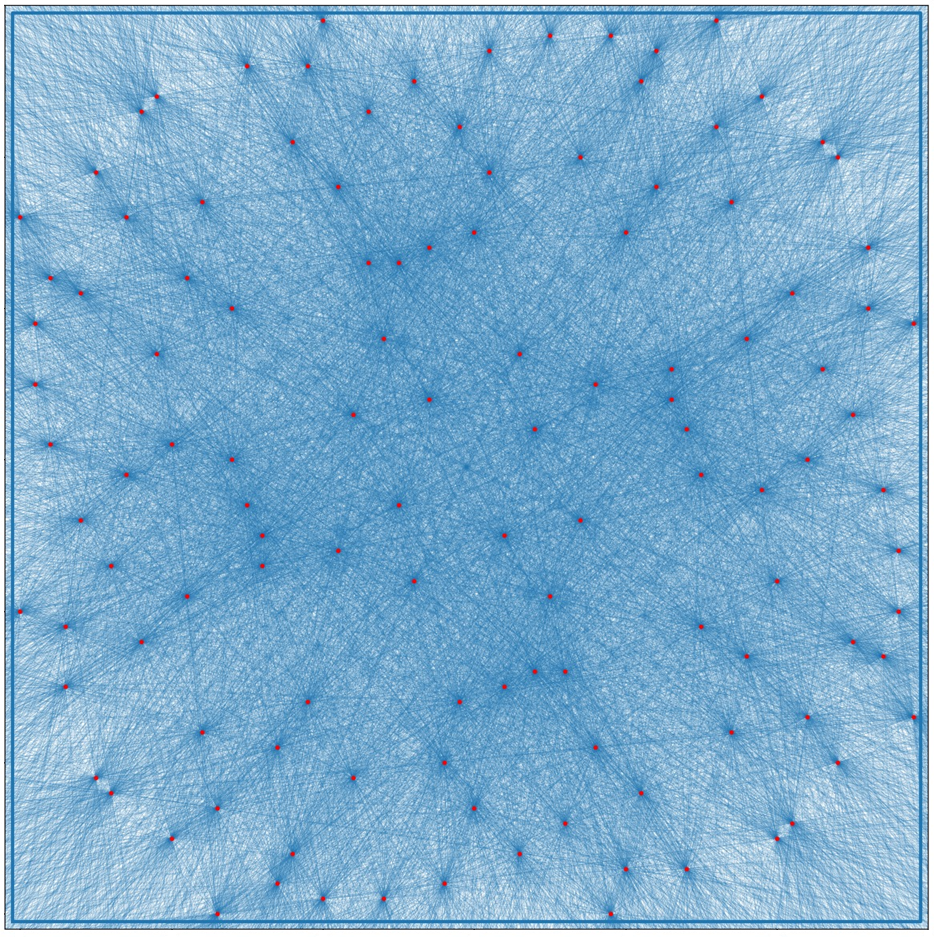}
				\caption{A no-three-in-line configuration of 120 pawns on a $60 \times 60$ chessboard from~\cite{Prellberg}, reproduced with kind permission from T.~Prellberg (with lines connecting the vertices on the right to show the complexity of the problem)}
				\label{fig:nothreeinline52}
			\end{figure}
		\end{center}

For the $n \times n$ board, since any row or column can contain at most two pawns, $2n$ pawns is a trivial upper bound. For many values of $n$ this upper bound can be achieved, including all odd $n$ up to 69 and all even $n$ up to 74, making $n = 71$ the smallest currently open case; see the papers~\cite{Flammenkamp1,Flammenkamp2,Guy-1968,Prellberg2026} and the databases~\cite{Flammenkampwebsite,Prellberg} for more recent developments, the constructions for the largest values of $n$ being due to Prellberg and Heule. A solution for $n = 60$ with 120 pawns due to Prellberg is shown in Figure~\ref{fig:nothreeinline52}. 62563 different configurations meeting the $2n$ upper bound can be found at~\cite{Wolfram}. There are typically many maximum configurations (the number of solutions with $2n$ pawns is given in OEIS entry A000769), and the aforementioned computational results are often concerned with finding maximum configurations with certain symmetry properties. Applications of the problem include optimal graph drawings with edges as straight line segments~\cite{Brassetal,wooddrawing}.

For a lower bound on the $n \times n$ board, Erd\H{o}s gave a construction containing $n-o(n)$ pawns in general position, which was published by Roth in~\cite{Roth} (interestingly Erd\H{o}s investigated this as a way to attack the Heilbronn triangle problem). This estimate was later improved to $\frac{3n}{2}-o(n)$ pawns by Hall et al.~\cite{Hall-1975}. A probabilistic argument by Guy and Kelly~\cite{Guy,Guy-1968} (with a subsequent slight correction, see~\cite{Voutier} for details) suggests that the true answer is $\frac{\pi n}{\sqrt{3}}+O(n)$, although see~\cite{Kaplan} for the summary of a talk by Kaplan that gives a counterargument. It is noted in~\cite{Prellberg2026} that the heuristic probabilistic argument only applies for $n \geq 493$, and so the fact that the $2n$ upper bound can be met for small $n$ is not strong evidence against the conjecture of Guy and Kelly. An early approach to the No-Three-In-Line Problem using neural networks can be found in~\cite{Tsuchiya95}, whilst a recent approach using machine learning and AI is presented and contrasted with the results of integer linear programming in~\cite{Ramanathan2025}.

A more general problem is to find the largest number $f_k(n)$ of pawns that can be placed on an $n \times n$ grid without $(k+1)$-in-line, with the trivial upper bound $f_k(n) \leq kn$. The paper~\cite{Kovacs-2025 algebraic} gave a randomised algebraic construction that places $\left ( 1-\frac{2}{k-1}\right ) kn$ or $\left ( 1-\frac{3}{k-1}\right ) kn$ pawns on the grid with no $(k+1)$-in-line for odd and even $k$ respectively, with better bounds for $k \leq 21$. A recent series of papers has made significant progress in proving that equality holds in $f_k(n) \leq kn$ using random graph theory; equality is proven for $k \geq \frac{5\sqrt{35}}{2}\sqrt{ n\log n }$ in~\cite{Kovacs-2025}, for all $k \geq 10^{37}$ in~\cite{GrebennikovKwan} and all $k \geq 3$ and sufficiently large $n$ in~\cite{Ghosal2026}.

The No-Three-In-Line Problem has been investigated in other geometrical settings. The article~\cite{Nagy-2023} examines the density of no-three-in-line sets on the infinite grid (known as the \emph{Extensible No-Three-In-Line Problem}); a construction is presented that contains $\Theta \left (\frac{n}{\log ^{1+\epsilon }n} \right )$ vertices from $[1,n]^2$ for every $n$ and it is conjectured that a linear bound is possible. Their lower bound is improved to $\Theta  \left (\frac{n}{\sqrt{ \log n}} \right )$ in~\cite{Ghosal} using a random construction. A relaxation of the extensible problem is examined in~\cite{Gabriel}, in which the number of vertices of the square grid $[1,n]^2$ that are not allowed to lie on a common line is a function of $n$. A variation in which the pawns must all lie on squares of the same colour (or, one could say, no three checkers in a line) is presented in~\cite{Prellbergcheckerboard}, which derives strong upper bounds via linear programming using the horizontal, vertical and diagonal lines, and presents strong evidence that the answer is $\alpha n+o(n)$, where $\alpha \approx 1.5768$.

The paper~\cite{Por} treats the No-Three-In-Line Problem on a three-dimensional grid and a higher dimensional version can be found in~\cite{Suk}, whereas the papers~\cite{Cooper,Ku,Misiak} treat the problem on a torus. The two articles~\cite{DongXu,Suk-spheres} investigate the problem of finding large subsets of the $d$-dimensional grid with no $d+2$ vertices lying on a common hyperplane or $d$-dimensional sphere, with an extension to quadratic surfaces in~\cite{Szabo}. Of particular interest is the number $f_{circ}(n)$ of pawns that can be placed on the $n \times n$ grid without four pawns lying on a common line or circle. Estimates for $f_{circ}(n)$ are given in~\cite{DongXu,Ghosal2026,GhosalGoenkaKeevash,Thiele,Thiele2}, with the state of the art being $2n-o(n) \leq f_{circ}(n) \leq 2.5n-1.5$, where the lower bound is from~\cite{Ghosal2026} and the upper bound is from~\cite{Thiele2}.

A generalised version in which configurations other than straight lines are forbidden is discussed in~\cite{Janosik}. The number of ways to arrange $k$ points in an $n \times n \times n$ triangular grid with no three in line is given in OEIS sequence A194136. Other geometrical investigations into no-three-in-line sets include~\cite{Balogh,Erdos}, as well as a random version in~\cite{BaloghLuo,Chen,Chen2}. The general position subset selection problem in discrete geometry also has an extensive literature, including~\cite{Froese,Payne}, of which we will not attempt to give an exhaustive overview.

In this paper we give a survey of the general position problem in the context of graph theory. We are therefore interested in sets of vertices of a graph $G = (V(G),E(G))$ which have a no-three-in-line property. This raises the question of what we mean by a `line' in a graph. As will be seen in Sections~\ref{sec:distancevariations} and~\ref{sec:variantsnotdistance}, this ambiguity has led to many interesting variations on the general position problem for graphs. The `classical' definition of the general position problem for graphs takes `line' to mean `shortest path.'

\begin{definition}\label{def:gp problem}
    For a graph $G$, a subset $S \subseteq V(G)$ is in \emph{general position}, or is a \emph{general position set}, if no three vertices of $S$ are contained in a common shortest path of $G$. The \emph{general position number} $\gp (G)$ of $G$, or the \emph{$\gp$-number} for short, is the number of vertices in a largest general position set of $G$. A general position set of $G$ containing $\gp (G)$ vertices is a \emph{$\gp $-set}. 
\end{definition}

For example, largest general position sets for the Clebsch graph and Frucht graph are shown in Figures~\ref{fig:Clebsch} and~\ref{fig:Frucht} respectively, with vertices of the set in red (up to symmetry, there is a unique solution for the Clebsch graph and two distinct solutions for the Frucht graph; the $\gp $-sets were verified computationally by G.~Erskine).

\begin{figure}\centering
\begin{tikzpicture}[scale=0.7,x=0.2mm,y=0.2mm,very thick,vertex/.style={circle,draw,minimum size=10,fill=lightgray}]
	\node at (-289.9,0) [vertex,fill=red] (v1) {};
	\node at (-205,205) [vertex,fill=red] (v2) {};
	\node at (-267.9,111) [vertex,fill=white] (v3) {};
	\node at (0,289.9) [vertex,fill=white] (v4) {};
	\node at (-111,-267.8) [vertex,fill=white] (v5) {};
	\node at (-111,267.8) [vertex,fill=white] (v6) {};
	\node at (267.9,111) [vertex,fill=red] (v7) {};
	\node at (205,205) [vertex,fill=red] (v8) {};
	\node at (-205,-205) [vertex,fill=white] (v9) {};
	\node at (-267.9,-110.9) [vertex,fill=white] (v10) {};
	\node at (110.9,-267.8) [vertex,fill=red] (v11) {};
	\node at (111,267.9) [vertex,fill=red] (v12) {};
	\node at (0,-289.9) [vertex,fill=red] (v13) {};
	\node at (267.9,-110.9) [vertex,fill=red] (v14) {};
	\node at (205,-205) [vertex,fill=white] (v15) {};
	\node at (289.9,0) [vertex,fill=white] (v16) {};
	\draw [bend right=18,color=black] (v1) to (v2);
	\draw [color=black] (v1) to (v3);
	\draw [bend right=45,color=black] (v1) to (v5);
	\draw [bend right=-18,color=black] (v1) to (v9);
	\draw [color=black] (v1) to (v16);
	\draw [bend right=18,color=black] (v2) to (v4);
	\draw [color=black] (v2) to (v6);
	\draw [bend right=45,color=black] (v2) to (v10);
	\draw [color=black] (v2) to (v15);
	\draw [bend right=-45,color=black] (v3) to (v4);
	\draw [color=black] (v3) to (v7);
	\draw [color=black] (v3) to (v11);
	\draw [color=black] (v3) to (v14);
	\draw [bend right=18,color=black] (v4) to (v8);
	\draw [color=black] (v4) to (v12);
	\draw [color=black] (v4) to (v13);
	\draw [color=black] (v5) to (v6);
	\draw [color=black] (v5) to (v7);
	\draw [color=black] (v5) to (v12);
	\draw [color=black] (v5) to (v13);
	\draw [bend right=-45,color=black] (v6) to (v8);
	\draw [color=black] (v6) to (v11);
	\draw [color=black] (v6) to (v14);
	\draw [color=black] (v7) to (v8);
	\draw [color=black] (v7) to (v10);
	\draw [bend right=-45,color=black] (v7) to (v15);
	\draw [color=black] (v8) to (v9);
	\draw [bend right=18,color=black] (v8) to (v16);
	\draw [color=black] (v9) to (v10);
	\draw [bend right=45,color=black] (v9) to (v11);
	\draw [bend right=-18,color=black] (v9) to (v13);
	\draw [color=black] (v10) to (v12);
	\draw [color=black] (v10) to (v14);
	\draw [color=black] (v11) to (v12);
	\draw [color=black] (v11) to (v15);
	\draw [bend right=-45,color=black] (v12) to (v16);
	\draw [bend right=45,color=black] (v13) to (v14);
	\draw [bend right=-18,color=black] (v13) to (v15);
	\draw [color=black] (v14) to (v16);
	\draw [bend right=-18,color=black] (v15) to (v16);
\end{tikzpicture}
\caption{A gp-set in the Clebsch graph}\label{fig:Clebsch}
\end{figure}

\begin{figure}\centering
\begin{tikzpicture}[x=0.2mm,y=0.2mm,very thick,vertex/.style={circle,draw,minimum size=5}]
	\node at (98,-28.2) [vertex,fill=red] (v1) {};
	\node at (45.8,-83.3) [vertex,fill=white] (v2) {};
	\node at (144.1,-130.6) [vertex,fill=red] (v3) {};
	\node at (178,54.7) [vertex,fill=white] (v4) {};
	\node at (-132.5,-138.1) [vertex,fill=white] (v5) {};
	\node at (-107.2,-7.6) [vertex,fill=red] (v6) {};
	\node at (-187.1,43.2) [vertex,fill=white] (v7) {};
	\node at (-95.3,168.4) [vertex,fill=red] (v8) {};
	\node at (-5.7,109.1) [vertex,fill=red] (v9) {};
	\node at (-47.5,29.1) [vertex,fill=white] (v10) {};
	\node at (89.7,168.4) [vertex,fill=white] (v11) {};
	\node at (3.9,-195.5) [vertex,fill=white] (v12) {};
	\draw [color=black] (v1) to (v2);
	\draw [color=black] (v1) to (v3);
	\draw [color=black] (v1) to (v4);
	\draw [color=black] (v2) to (v10);
	\draw [color=black] (v2) to (v12);
	\draw [color=black] (v3) to (v4);
	\draw [color=black] (v3) to (v12);
	\draw [color=black] (v4) to (v11);
	\draw [color=black] (v5) to (v6);
	\draw [color=black] (v5) to (v7);
	\draw [color=black] (v5) to (v12);
	\draw [color=black] (v6) to (v7);
	\draw [color=black] (v6) to (v10);
	\draw [color=black] (v7) to (v8);
	\draw [color=black] (v8) to (v9);
	\draw [color=black] (v8) to (v11);
	\draw [color=black] (v9) to (v10);
	\draw [color=black] (v9) to (v11);
\end{tikzpicture}
\caption{A gp-set in the Frucht graph}\label{fig:Frucht}
\end{figure}

The general position problem for a graph $G$ therefore consists in finding the largest set of vertices $S$ of $G$ such that any shortest path of $G$ contains at most two vertices of $S$. We will regard as the origin of the general position problem for graphs the article~\cite{Korner-1995} by K\"{o}rner from 1995, which concerns the general position problem in the hypercube. However, the article notes that ``this problem has been known for a long time in a different form'' and shows that it is related to several problems in combinatorial set theory. Interestingly, the authors also note that ``ours is not the only possible definition of a line in a combinatorial setting''; for example, another natural definition of `line' in a hypercube leads to the fascinating world of Ramsey problems on arithmetic progressions. 

At the time, the article~\cite{Korner-1995} did not lead to further work on the general position problem for general graphs. This problem started to gain wider attention twenty years later with the appearance of two independent articles,~\cite{ullas-2016} by Chandran and Parthasarathy in 2016 (under the name \emph{geodetic irredundant sets}), and~\cite{ManuelKlavzar-2018a} by Manuel and Klav\v{z}ar in 2018 under the present name.

The arXiv version~\cite{manuel-2017arxiv} of~\cite{ManuelKlavzar-2018a} illustrates the general position problem in a useful and suggestive form. Suppose that a collection of robots is stationed on the nodes of a network. The robots communicate with each other by sending signals along the shortest paths between them. To ensure that these signals are not blocked, we require that no shortest path between a pair of robots be occupied by a third robot, and, subject to this condition, we wish to maximise the number of robots. This problem has important applications in robotic navigation; for example, there is an extensive literature in computer science on how a swarm of robots can reach such a configuration (see for example~\cite{Adhikary2019,Bhagat2020,Bhagat2019,Bhagat2017,CiDiDiNav2023,DiLuna-2017,DiLuna2014}).

The general position problem and its relatives are also connected to graph convexities, as discussed in~\cite{Araujo-2025}. The study of convexities in graph theory has a longer history than the general position problem, and has a correspondingly larger literature; we refer the reader to the book~\cite{Pelayo-2013} for more details. Given a collection $\Pi $ of paths of a graph $G$, the \emph{interval} $\I_{\Pi }[u,v]$ is $\{ u,v\} $ together with the set of vertices that lie on some $u,v$-path belonging to $\Pi $. The \emph{closure} $\I_{\Pi }[S]$ of a subset $S \subseteq V(G)$ is $\bigcup _{u,v \in S}\I_{\Pi }[u,v]$ and $S$ is \emph{convex} if $\I_{\Pi }[S] = S$. The \emph{convex hull} of a subset $S$ is the smallest convex set containing $S$. The main question of graph convexity is then to find the smallest subsets with closure or convex hull equal to $V(G)$. The family of convex subsets constitutes a \emph{convexity} with associated \emph{interval function} $\I$; these notions can also be defined axiomatically in a more abstract setting. Associated with each convexity is a position problem, namely to find large subsets $S \subseteq V(G)$ such that for any $u,v \in S$ the interval $\I _{\Pi }[u,v]$ does not contain any vertices of $S$ apart from $u$ and $v$. This connection has served as the inspiration for several developments in position problems, including the monophonic position problem defined in Subsection~\ref{subsec:mp} and the all-path position number in Subsection~\ref{subsec:all paths}.  

The plan of this survey is as follows. In Subsection~\ref{subsec:terminology} we briefly define some terminology that will be used in our discussion. Section~\ref{sec:the-number} deals with the `classical' general position problem, covering bounds for the general position number in Subsection~\ref{subsec:bounds}, complexity results in Subsection~\ref{subsec:complexity}, exact results for various graph classes in Subsection~\ref{subsec:exact} and behaviour of the general position number under graph products in Subsection~\ref{subsec:products} and other graph operations in Subsection~\ref{subsec:otheroperations}. Section~\ref{sec:distancevariations} describes some of the variations on the general position number that are based on distance, whilst Section~\ref{sec:variantsnotdistance} deals with variations that are not based on distance (principally the \emph{monophonic position problem}, in which `shortest path' is replaced by `induced path'). Section~\ref{sec:otherdirections} surveys some of the recent new directions for research into position problems, including a dynamic version of the general position problem, smallest maximal general position sets, general position colourings, fractional general position and general position games. Finally, Section~\ref{sec:openproblems} gives a list of important open problems.

There is one important variation on the general position problem for which we do not give a survey here, namely the mutual visibility problem. This problem was introduced to graph theory in~\cite{DiStefano-2022} by Di Stefano and was inspired by the applications to robotic navigation discussed above. Whereas in the general position problem we require all shortest paths between a pair of robots to be free of a third robot, in the mutual visibility problem we demand only that at least one shortest path is left free for communication. This field has such a wealth of literature and is growing so quickly that it merits its own survey paper.

\subsection{Terminology}\label{subsec:terminology}

In this subsection we define terminology that is peculiar to the general position problem. For standard graph theory notation we refer the reader to~\cite{Graphsdigraphs}. 

For a positive integer $k$ we will use the notation $[k] = \{1,\ldots, k\}$.
The {\em order} of a graph $G$ is the number of its vertices, denoted by $n(G)$, and the {\em size} of $G$ is the number of its edges, denoted by $m(G)$. If $u\in V(G)$, then $N_G(u)$ denotes the open neighbourhood of $u$ and $N_G[u]$ its closed neighbourhood, that is, $N_G(u)$ is the set of neighbours of $u$ and $N_G[u] = N_G(u) \cup \{u\}$. Vertices $u$ and $v$ of a graph $G$ are {\em true twins} if $N_G[u] = N_G[v]$. A vertex is a \emph{support vertex} if it is adjacent to a leaf. A vertex $u$ of $G$ is \emph{simplicial} if its neighbourhood induces a clique. The set of all simplicial vertices of a graph will be denoted by $S(G)$ and their number by $\s(G)$. The number of leaves of a tree $T$ will be written as $\ell(T)$ (so $\ell(T) = \s(T)$). The {\em join} $G \vee H$
of graphs $G$ and $H$ is obtained from the disjoint union of $G$ and $H$ by adding all possible edges between $G$ and $H$. For $X\subseteq V(G)$, $G[X]$ is the subgraph of $G$ induced by $X$. By $\alpha(G)$ and $\omega(G)$ we represent the independence number and the clique number of $G$ respectively. 

Unless stated otherwise, the graph distance considered here is the shortest-path distance; for $u,v\in V(G)$, their distance in $G$ is denoted by $d_G(u,v)$. A subgraph $H$ of $G$ is {\em isometric} if $d_H(u,v) = d_G(u,v)$ holds for any $u,v\in V(H)$. $H$ is {\em convex} if for all $u,v\in V(H)$ any shortest $u,v$-path in $G$ lies in $H$. The {\em interval} between vertices $u$ and $v$ of a graph $G$ is defined as $I_G(u,v) = \{w:\ d_G(u,v) = d_G(u,w) + d_G(w,v)\}$. A subset $S \subset V(G)$ is a \emph{geodetic set} if for any vertex $x \in V(G)-S$ there are $u,v \in S$ such that $x$ lies on a shortest $u,v$-path. The number of vertices in a smallest geodetic set is the \emph{geodetic number} $\g (G)$.

The {\em complete multipartite graph} $K_{n_1,\dots,n_t}$, $t \geq 2$, $n_1 \geq \dots \geq n_t$, is the graph of order $\sum _{i=1}^tn_i$ with vertex set $\bigcup _{i=1}^tV_i$, where $|V_i| = n_i$, $i \in [t]$ and $V_i \cap V_j = \emptyset $ for $i \neq j$, where two vertices $u \in V_i$ and $v \in V_j$ are adjacent if and only if $i \neq j$.

If $G$ and $H$ are graphs, then their {\em Cartesian product}, {\em direct product}, {\em strong product}, and {\em lexicographic product} are respectively denoted by $G\cp H$, $G\times H$, $G\boxtimes H$, and $G\circ H$. Each of these products has the vertex set $V(G)\times V(H)$. Vertices $(g,h)$ and $(g',h')$ are adjacent in $G\cp H$ if either $g=g'$ and $hh'\in E(H)$, or $gg'\in E(G)$ and $h=h'$. Vertices $(g,h)$ and $(g',h')$ are adjacent in $G\times H$ if $gg'\in E(G)$ and $hh'\in E(H)$. For the strong product we have $E(G\boxtimes H) = E(G\cp H)  \cup E(G\times H)$. Finally, vertices $(g,h)$ and $(g',h')$ are adjacent in $G\circ H$ if either $gg'\in E(G)$, or $g=g'$ and $hh'\in E(H)$. If $h\in V(H)$ and $*\in \{\cp, \times, \boxtimes, \circ\}$, then the subgraph of $G * H$ induced by the vertex set $\{(g,h):\ g\in V(G)\}$ is called a {\em $G$-layer} of $G*H$ and denoted by $G^h$. For $g\in V(G)$ the $H$-layer $^gH$ is defined analogously.    

Given two graphs $G$ and $H$ with $V(G) = \{v_1, \ldots ,v_{n}\}$, the {\em corona product} graph $G\odot H$ is formed by taking one copy of $G$ and $n$ disjoint copies of $H$, call them $H^1,\ldots, H^{n}$, and for each $i \in [n]$ adding all the possible edges between $v_i\in V(G)$ and every vertex of $H^i$. 

If $G = (V(G), E(G))$ is a graph and $X \subseteq V(G)$, then we say that vertices $u, v \in V (G)$ are {\em $X$-positionable} if for any shortest $u,v$-path $P$ we have $V(P)\cap X \subseteq \{u,v\}$. Note that if $uv\in E(G)$, then $u$ and $v$ are $X$-positionable for any $X$. Using this terminology we can say that the set $X$ is a general position set if every pair $u, v \in X$ is $X$-positionable. 

\section{General Position Number}\label{sec:the-number}

General position sets have a structure that can be characterised. To do this, we need some definitions. Let $G$ be a connected graph, $X\subseteq V(G)$, and ${\cal P} = \{X_1, \ldots, X_p\}$ a partition of $X$. Then ${\cal P}$ is \emph{distance-constant} if for any $i,j\in [p]$, $i\ne j$, the distance $d_G(u,v)$, $u\in X_i$, $v\in X_j$, is independent of the selection of $u$ and $v$ in $X_i$ and $X_j$. If ${\cal P}$ is a distance-constant partition, and $i,j\in [p]$, $i\ne j$, then the distance $d_G(X_i, X_j)$ can be defined as the distance between a vertex of $X_i$ and a vertex of $X_j$. A distance-constant partition ${\cal P}$ is {\em in-transitive} if $d_G(X_i, X_k) \ne d_G(X_i, X_j) + d_G(X_j,X_k)$ holds for arbitrary pairwise distinct $i,j,k\in [p]$. 

\begin{theorem} {\rm \cite[Theorem~3.1]{AnaChaChaKlaTho}}
\label{thm:gpsets-characterised}
Let $G$ be a connected graph. Then $S\subseteq V(G)$ is a general position set if and only if the components of $G[S]$ are cliques, the vertices of which form an in-transitive, distance-constant partition of $S$.
\end{theorem}

Using Theorem~\ref{thm:gpsets-characterised}, one can deduce that in a complete multipartite graph, an arbitrary general position set is either a subset of one partite set or contains at most one vertex from each of the partite sets. This fact is a widely known folklore result, but has been stated explicitly in sources including~\cite{Thomas-2024b,rather-2026+}.

\begin{corollary} 
\label{cor:complete-multi}
If $G = K_{n_1,\dots, n_t}$, $t\ge 2$, and $M = \max_{i\in [t]}n_i$, then $\gp(G) = \max \{M, t\}$. 
\end{corollary}

\subsection{Bounds}\label{subsec:bounds}

An obvious upper bound for the general position number of a graph is the total number of its vertices. Since any pair of vertices is always in general position in a graph $G$, we have $2\leq \gp(G)\leq n(G)$. Both of these bounds are sharp. A trivial lower bound is $\gp (G) \geq \omega (G)$, as any clique is in general position.

Graphs with $\gp(G) \in \{2, n(G), n(G) - 1\}$ are characterised in the paper~\cite{ullas-2016}. In particular, $\gp(G) = n(G)$ if and only if $G = K_n$, and $\gp(G) = 2$ if and only if $G \in \{C_4, P_n\}$ \cite[Theorem 2.10]{ullas-2016}. By~\cite[Theorem 2.4]{ullas-2016}, for any connected graph $G$ we have $\gp(G) \leq n(G) - \mathrm{diam}(G) + 1$. Consequently,~\cite[Theorem 3.1]{ullas-2016} proves that $\gp(G) = n(G) - 1$ if and only if
$G = K_1 + \bigcup_j m_j K_j$ with $\sum m_j \geq 2$ or $G = K_n - \{e_1, \ldots, e_k\}$ with $k \in [n - 3]$, where the $e_i$'s are edges in $K_n$ that are all incident to a common vertex.

The paper~\cite{Thomas-2020} characterises graphs $G$ satisfying $\gp(G) = n(G) - 2$ into eight families of graphs formed by modifying cycles, paths, and complete graphs. Four families ($\mathcal{F}_1$ to $\mathcal{F}_4$) consist of graphs with diameter three, while the other four ($\mathcal{F}_5$ to $\mathcal{F}_8$) have diameter two. 

\begin{theorem}{\rm\cite[Theorem 3.1]{Thomas-2020} }
If $G$ is a connected graph of order at least four, then ${\rm gp}(G) = n(G)-2$  if and only if $G\in \cup_{i=1}^8{\cal F}_i.$
\end{theorem}

Returning our attention to upper bounds on the general position number, we say that a set of subgraphs $\{H_1, \ldots, H_k\}$ of a graph $G$ is an \emph{isometric cover} of $G$ if each $H_i$ ($i \in [k]$) is isometric in $G$, and $\bigcup_{i=1}^k V(H_i) = V(G)$. Each isometric cover of $G$ provides an upper bound on $\gp(G)$.

\begin{theorem}{\rm \cite[Theorem 3.1]{ManuelKlavzar-2018a} ({\bf Isometric Cover Lemma})} If $\left\{H_1, \ldots, H_k\right\}$ is an isometric cover of $G$, then
$$
\gp(G) \leq \sum_{i=1}^k \gp\left(H_i\right).
$$
\end{theorem}
The \emph{isometric-path number} of a graph $G$, which we will write as $\rho (G)$, is the smallest number of isometric paths (shortest paths) needed to cover all the vertices of $G$ (note that $\ip (G)$ is often used for this invariant, but we wish to avoid a notation clash with the independent position number). Similarly, the \emph{isometric cycle number} of $G$, written as $\rho _C(G)$, is the smallest number of isometric cycles required to cover all the vertices of $G$. If $G$ cannot admit a cover by isometric cycles, we set $\rho _C(G) = \infty$. Then, the Isometric Cover Lemma yields that for any graph $G$, $\gp(G) \leq 2 \rho(G)$ and $\gp (G) \leq 3 \rho _C(G)$ (see~\cite[Corollary 3.2]{ManuelKlavzar-2018a}).

If $v$ is a vertex of a graph $G$, let \( \rho(v, G) \) denote the minimum number of isometric paths, all starting at \( v \), that cover \( V(G) \). With this concept in hand, we have the following implicit bound in~\cite{ManuelKlavzar-2018a}.

\begin{theorem}{\rm \cite[Theorem 3.3]{ManuelKlavzar-2018a}} If $S$ is a general position set of a graph $G$ and $v \in S$, then
$$
|S| \leq \rho(v, G)+1
$$
In particular, if $v$ belongs to some $\gp $-set, then $\gp(G) \leq \rho(v, G)+1$.
\end{theorem}

Note that if $S$ is a general position set of $G$, and $W$ and $W'$ are two cliques from the induced subgraph $G[S]$, then by Theorem~\ref{thm:gpsets-characterised}, $d_G(W,W') \ge 2$, that is, $W$ and $W'$ are {\em independent cliques}. We set $\alpha ^{\omega }(G)$ to be the maximum number of vertices in a union of pairwise independent cliques of $G$. Here we allow the union of cliques to consist of a single clique, so $\alpha ^{\omega }(G) \geq \omega (G)$. 

\begin{corollary} {\rm \cite[Theorem~4.1]{Ghorbani-2021}}
\label{thm:rho-diameter-2}
If $G$ is a connected graph, then $\gp(G)\le \alpha ^{\omega }(G)$. Moreover, if $\diam(G) \in [2]$, then $\gp(G) = \alpha ^{\omega }(G)$. 
\end{corollary}

Since a join $G \vee H$ has diameter at most two, by Corollary~\ref{thm:rho-diameter-2} the $\gp $-number of $G \vee H$ is given by $\max \{ \alpha ^{\omega } (G),\alpha ^{\omega } (H),\omega (G)+\omega (H)\} $. If $G$ is bipartite and $S$ is a general position set of $G$ with $|S|\ge 3$, then Theorem~\ref{thm:gpsets-characterised} implies that $S$ must be independent. For small diameter graphs more can be concluded as the next result asserts. 

\begin{theorem} {\rm \cite[Theorem~5.1]{AnaChaChaKlaTho}}
If $G$ is a connected, bipartite graph on at least three vertices, then $\gp(G) \leq \alpha(G)$.  Moreover, if $\diam(G) \in \{2,3\}$, then $ \gp(G)=\alpha (G)$.
\end{theorem}

\subsubsection{Strong resolving graphs}
\label{sec:strong-resolving}

A vertex $x$ of a connected graph $G$ is \emph{maximally distant} from a vertex $y$ if every $z\in N_G(x)$ satisfies $d_G(y,z)\le d_G(y,x)$. If $x$ is maximally distant from $y$, and $y$ is maximally distant from $x$, then $x$ and $y$ are \emph{mutually maximally distant} (MMD for short). The {\em strong resolving graph} $G\sr$ of $G$ has $V(G\sr) = V(G)$ and two vertices are adjacent in $G\sr$ if they are MMD in $G$, see~\cite{Oellermann2007}. The notion of the strong resolving graph was introduced as a useful tool for studying the strong metric dimension of graphs. In particular, the problem of determining the strong metric dimension of a graph $G$ can be reduced to the problem of finding the vertex cover number of its strong resolving graph $G{\sr}$.

There is a fundamental connection between general position sets and strong resolving graphs. 

\begin{theorem} {\rm \cite[Theorem~3.1]{Klavzar-2019}}
\label{thm:gp-versus-strong-resolving-graphs}
If $G$ is a connected graph, then $\gp(G)\ge \omega(G\sr)$. Moreover, equality holds if and only if $G$ contains a $\gp $-set that induces a complete subgraph of $G\sr$.
\end{theorem}

For diameter two graphs with no true twins, Theorem~\ref{thm:gp-versus-strong-resolving-graphs} can be strengthened. 

\begin{proposition} {\rm \cite[Proposition~2.4]{Klavzar-2019}}
If $G$ has no true twins and $\diam(G) = 2$, then $\gp(G) = \omega(G\sr)$ if and only if $\gp(G) = \alpha(G)$.
\end{proposition}

On the other hand, as demonstrated in~\cite[Proposition~3.5]{Klavzar-2019}, for any integers $r\ge t\ge 2$, there exists a graph $G$ such that $\gp(G) = r$ and $\omega(G\sr) = t$.

A structural characterisation of graphs achieving equality in Theorem~\ref{thm:gp-versus-strong-resolving-graphs} is not known and  seems to be elusive because of the great variety of different structures that can appear. Among other classes of graphs, it contains block graphs (so in particular complete graphs and trees), complete multipartite graphs, special corona products, and direct products of complete graphs. 

\subsection{Complexity}\label{subsec:complexity}

The basic complexity question concerning the decision problem for the general position number is formally defined as follows:
\begin{definition}\label{def:probgp}
        {\sc General position set} \\
        {\sc Instance}: A graph $G$, a positive integer $k\leq |V(G)|$. \\
        {\sc Question}: Is there a general position set $S$ for $G$ such that $|S|\geq k$?
\end{definition}
The general position subset selection problem from discrete geometry has been proven to be NP-hard (see~\cite{Froese, Payne}). Manuel and Klav\v{z}ar~\cite{ManuelKlavzar-2018a} showed that the {\sc General Position Set} problem is NP-complete for arbitrary graphs. The proof of this is based on a reduction from the NP-complete {\sc Maximum Independent Set} problem.

\begin{theorem}{\rm \cite[Theorem 5.1]{ManuelKlavzar-2018a}}\label{thm:fundamental_complexity_result}
{\sc General position set} is NP-Complete. 
\end{theorem}

Since determining the general position number of a graph is difficult in general, it is natural to consider different algorithmic approaches. Kor\v{z}e and Vesel~\cite{KorzeVesel-2023} described a reduction from the problem of finding a general position set in a graph to the satisfiability problem and applied this approach to the case of general position sets of hypercubes. In~\cite{Hamed-2026} three different heuristic algorithms for the general position problem are proposed, namely integer linear programming, a genetic algorithm, and simulated annealing. These algorithms were tested on relatively large graphs from different areas of graph theory, including chemical graphs and Cayley graphs; all three algorithms found the correct answer in all cases, but it appears that integer linear programming is the slowest of the three.

\subsection{Exact results}\label{subsec:exact}

As the general position problem is NP-complete, it is natural to consider the restriction of the problem to some simple graph classes. From the seminal papers~\cite{ullas-2016,ManuelKlavzar-2018a} we recall these results. If $n\ge 2$, then $\gp(K_n) = n$. If $n\ge 5$, then $\gp(C_n) = 3$, while $\gp(C_4) = 2$. If $n\ge 2$, then $\gp(P_n) = 2$. Further, if $G$ is a block graph, then $\gp(G) = \s(G)$~\cite[Theorem 3.6]{ManuelKlavzar-2018a}; in particular, if $T$ is a tree, then $\gp(T) = \ell(T)$~\cite[Theorem 2.5]{ullas-2016}. 

\subsubsection{Kneser graphs}
The \emph{Kneser graph} $K(n,k)$ is the graph whose vertex set consists of all $k$-subsets of the set $[n]$, with an edge between any two subsets if and only if they are disjoint.
Ghorbani et al.~\cite{Ghorbani-2021} computed $\gp(K(n, 2))$ and $\gp(K(n, 3))$ for all $n$. In particular,~\cite[Theorem 2.2]{Ghorbani-2021} gives $\gp(K(n, 2)) = n - 1$ when $n \geq 7$, and~\cite[Theorem 2.4]{Ghorbani-2021} gives $\gp(K(n, 3)) = \binom{n-1}{2}$ when $n \geq 9$. In~\cite[Theorem 2.3]{Ghorbani-2021}, it is proved that for any fixed $k$, if $n$ is sufficiently large, the equality $\gp(K(n, k)) = \binom{n-1}{k-1}$ holds. 
\begin{theorem}{\rm\cite[Theorem 2.3]{Ghorbani-2021}}\label{kneser_G_thm2.3} Let $n$ and $k$ be positive integers with $n \geq 3k - 1$. If, for all $t$ such that $2 \leq t \leq k$, the inequality $k^t \binom{n-t}{k-t} + t \leq \binom{n-1}{k-1}$ holds, then $\gp(K(n, k)) = \binom{n-1}{k-1}$.
    
\end{theorem}
For fixed $k$ and $t = 2$, the above inequality is satisfied when $n \geq k^3 - k^2 + 2k - 1$, while Patk\'os~\cite[Theorem 1.2]{Patkos-2019} improved this result by proving that the same conclusion holds for $n \geq 2.5k - 0.5$. 
\begin{theorem}{\rm \cite[Theorem 1.2]{Patkos-2019}}
     If $n,k \geq 4$ are integers with $n \geq 2k+1$, then $\gp(K(n,k)) \leq \binom{n-1}{k-1}$.
Moreover, if $n \geq 2.5k-0.5$, then $\gp(K(n,k)) = \binom{n-1}{k-1}$, while if $2k+1 \leq n < 2.5k-0.5$, then $\gp(K(n,k)) < \binom{n-1}{k-1}$.
\end{theorem}
As noted in~\cite{Patkos-2019}, it remains an open problem to determine $\gp(K(n, k))$ for $2k + 1 \leq n < 2.5k - 0.5$.


\subsubsection{Maximal outerplanar graphs}
\emph{Outerplane graphs} are planar graphs with a plane embedding such that every vertex lies on the boundary of the exterior region. An outerplane graph is \emph{maximal outerplane} if adding any edge results in a non-outerplane graph. A \emph{straight linear $2$-tree} on $n$ vertices is a graph $G_n$ with vertex set $V(G_n) = [n]$, where vertices $i$ and $j$ are adjacent in $G_n$ if and only if $0 < |i - j| \leq 2$.

The paper~\cite{TianXuChao-2023} focuses on the general position number of maximal outerplane graphs and presents several key results. It is shown in~\cite[Lemma 2.8]{TianXuChao-2023} that the general position number of a maximal outerplane graph is bounded in terms of the maximum degree as:
$\gp(G) \geq \lfloor \frac{2(\Delta(G) + 1)}{3}\rfloor$. In the next result, maximal outerplane graphs with $\gp(G)=3$ are characterised.

\begin{theorem}{\rm \cite[Theorem 3.1]{TianXuChao-2023}}\label{max_out_plane_gp3}  
If $G$ is a maximal outerplane graph with $n(G)\ge 7$, then $\gp(G) = 3$ if and only if $G$ is a straight linear $2$-tree.
\end{theorem}
In~\cite[Theorem 3.4]{TianXuChao-2023}, it was proved that for any maximal outerplane graph $G$ of order $n \geq 6$, $\gp(G) \leq \left\lfloor \frac{2n}{3} \right\rfloor$. The paper~\cite{TianXuChao-2023} then focuses on all maximal outerplane graphs that attain this bound.

\subsubsection{Additional exact results}

The general position numbers of glued binary trees $GT(r)$, $r \ge 2$, and of Bene\v{s} networks $BN(r)$, $r\ge 1$, were determined in~\cite{ManuelKlavzar-2018a} and~\cite{ ManuelKlavzar-2018b} respectively. We do not give the definitions of these networks here, but recall that if $r \ge 2$, then $\gp(GT(r)) = 2^r$~\cite[Proposition 3.8]{ManuelKlavzar-2018a}, and that if $r\ge 1$, then $\gp(BN(r)) = 2^{r+1}$~\cite[Theorem 4.1]{ManuelKlavzar-2018b}. Similar research has been done in~\cite{Prabha-2020, Prabha-2021, Prabha-2023}. In~\cite{Prabha-2020}, the general position number is determined for hyper trees and shuffle hyper trees. In~\cite{Prabha-2021}, the general position number is obtained for hexagonal networks, honeycomb networks, silicate networks, and oxide networks, while~\cite{Prabha-2023} reports the general position number of butterfly networks. 

Let $H$ be a graph and let $K_n - E(H)$ be the graph obtained from $K_n$ by deleting the edges of a copy of $H$. Applying Corollary~\ref{thm:rho-diameter-2}, the next exact results were presented in~\cite{AnaChaChaKlaTho}.  
\begin{itemize}
\item $\gp(K_n - E(K_k)) = \max \{ k, n- k + 1 \}$, where $2\leq k < n$.
\item $\gp(K_n - E(K_{1,k})) = \max \{ k+1, n-1 \}$, where $ 2\leq k < n$. 
\item $\gp(K_n - E(P_k)) = \max \{3, n- k + \lceil\frac{k}{2}\rceil \}$, where $3\leq k< n$. 
\item $\gp(K_n - E(K_{r,s})) = \max \{ r+s, n-r \}$, where $ 2\leq r\leq s$ and $r +s < n$. 
\item $\gp(K_n - E(W_k)) = \max \{3, n- k + \lfloor\frac{k-1}{2}\rfloor \}$, where $ 5\leq k< n$. 
\item $\gp(K_n - E(C_k)) =
   \left\{
     \begin{array}{ll}
       \max \{ 3, n- k + \lfloor\frac{k}{2}\rfloor\}, & 5 \le k < n; \\
       \max \{ 4, n-2 \}, & k = 4.
     \end{array}
   \right.$
\end{itemize}

Let $G$ be a bipartite graph and $A, B$ its bipartition. Furthermore, set 
$$M_G = \{u\in A:\ {\rm deg}_G(u) = |B|\}\, \cup\, \{u\in B:\ {\rm deg}_G(u) = |A|\}\,,$$
and let $\eta(G)$ be the maximum order of an induced complete bipartite subgraph of $G$. With these concepts the general position number of the complement of a bipartite graph can be determined.  

\begin{theorem}{\rm \cite[Theorem 5.2]{AnaChaChaKlaTho}}
\label{thm:bip-complement}
If $G$ is a bipartite graph with bipartition $A,B$, then $\gp(\overline{G})$ is equal to 
$$   \left\{
     \begin{array}{ll}
       n(G); & \hspace*{-1.1cm}\diam(\overline{G}) \in \{1, \infty\}\,, \\ 
       \max \{\alpha(G), \eta(G)\}; & \hspace*{-0.1cm} \diam (\overline{G}) = 2\,, \\
\max \{\alpha(G), \eta(G\setminus (M_G\cap A)), \eta(G \setminus (M_G\cap B)), |M_G|\}; & \hspace*{-0.1cm} \diam(\overline{G}) = 3\,.
\end{array}
\right.$$
\end{theorem}
Theorem~\ref{thm:bip-complement} implies that if $T$ is a tree, then 
$\gp(\overline{T}) = \max \{\alpha(T), \Delta(T) +1\}$. Another interesting consequence of Theorem~\ref{thm:bip-complement} is: 

\begin{corollary}{\rm \cite[Corollary 5.4]{AnaChaChaKlaTho}}
\label{cor:grids}
If $n,m\ge 2$, then 
$$\gp(\overline{P_n\cp P_m}) =
   \left\{
     \begin{array}{ll}
       4; &  n = m = 2,\\
       \left\lceil \frac{n}{2}\right\rceil \left\lceil \frac{m}{2}\right\rceil + \left\lfloor \frac{n}{2}\right\rfloor \left\lfloor \frac{m}{2}\right\rfloor; & \text{otherwise}\,.
     \end{array}
   \right.$$
\end{corollary}

The paper~\cite{Ghorbani-2021} also solves the general position problem for line graphs of complete graphs $L(K_n)$.

\begin{theorem}{\rm \cite[Theorem 4.4]{Ghorbani-2021}} 
If $n\ge 3$, then 
\[\gp(L(K_n)) = \left\{ {\begin{array}{ll}
n; & 3\mid n\,, \\
n-1; & 3\nmid n\,.
\end{array}} \right.\]
\end{theorem}

The {\em wheel graph} $W_n$, $n\ge 3$, is the graph obtained from $C_n$ by adding a new vertex and connecting it to all vertices of $C_n$.

\begin{theorem}{\rm \cite[Theorem 4.1]{yao-2022}} 
If $n\ge 3$, then 
\[\gp(W_n) = \left\{ {\begin{array}{ll}
4; & n = 3\,, \\
3; & n\in \{4,5\}\,, \\
\left\lfloor\frac{2}{3}n\right\rfloor; & n\ge 6\,.
\end{array}} \right.\]
\end{theorem}

The main focus of the paper~\cite{yao-2022} is on the general position number of cactus graphs, for which lower and upper bounds are proved, as well as partial exact results. The next theorem is one of the notable results of the paper. 

\begin{theorem}{\rm \cite[Theorem 2.5]{yao-2022}}\label{thm:cactus} 
If $G$ is a cactus graph with $k$ cycles and $t$ pendant edges, then $\gp(G) \le \max \{3, 2k+t\}$.
\end{theorem}

Moreover,~\cite[Theorem 2.5]{yao-2022} also characterises the class of cactus graphs that attain equality in Theorem~\ref{thm:cactus} as the cactus graphs in which all cycles are `good'; see the paper for the definition of a `good cycle'. 

If $p \geq 3$ and $n \geq 1$, then the {\em Sierpi\'nski graph} $S_p^n$ has vertex set $V (S_p^n) = \{0,1,\ldots, p-1\}^n$, whilst vertices $i_1\cdots i_n$ and $j_1 \cdots j_n$ are adjacent if there exists an index $h  \in [n]$, such that (i) if $t < h$, then $i_t = j_t$, (ii) $i_h \neq j_h$, and (iii) if $t > h$, then $i_t = j_h$ and $j_t = i_h$, see~\cite{hinz-2017}. In~\cite[Corollary~3.2]{Roy-2026} it is proved that if $p\geq 3$ then, 
$$\gp(S_p^2) = \begin{cases}
\frac{(p+1)^2}{4}; & p \text{ odd}\,,\\[5pt]
\frac{p(p+2)}{4}; & p \text{ even}\,,
\end{cases} 
\quad \text{and} \quad \# \gp(S_p^2) = \begin{cases}
\binom{p}{\frac{p+1}{2}}; & p \text{ odd}\,,\\[5pt]
\binom{p+1}{\frac{p+2}{2}}; & p \text{ even}\,,
\end{cases}$$
where  $\#\gp(G)$ indicates the number of $\gp$-sets of a graph $G$. Moreover, \cite[Corollary~4.5]{Roy-2026} asserts that if $n\geq 2$, then $\gp(S_3^n) = 3^{n-2} + 3$.

The Sierpi\'nski triangle graphs, also called Sierpi\'nski gasket graphs, are an interesting family of graphs that converge towards the Sierpi\'nski triangle fractal. They can be defined iteratively as follows: the graph $ST_3^0$ is the triangle $K_3$, and when for $n \geq 1$ the graph $ST_3^n$ has been constructed, $ST_3^{n+1}$ can be formed from three copies of $ST_3^n$ arranged into a triangle by identifying three corner vertices (see for example \url{mathworld.wolfram.com/SierpinskiGasketGraph.html}). Hence the graph $ST_3^n$ has order $\frac{3(3^{n-1}-1)}{2}$. The general position numbers of these graphs are studied (along with the mutual visibility numbers) in~\cite{KorzeVesel-2024}, which shows that $ST_3^n$ has a unique $\gp $-set, which has cardinality three for $n = 1$ (trivially) and cardinality $3^{n-1}+3$ for $n \geq 2$.

\subsubsection{General position polynomials}

A useful graph theoretical technique is to study polynomials, the coefficients of which count combinatorial objects of a certain type; the analytical properties of the polynomial can then be used to shed light on these structures. Important examples include \emph{chromatic polynomials}~\cite{fialho-2026}, \emph{matching polynomials}~\cite{chen-2026}, \emph{domination polynomials}~\cite{mertens-2024} and \emph{independence polynomials}~\cite{kadrawi-2025}. Inspired by this, the \emph{general position polynomial} $\psi (G)$ of a graph $G$ is defined in~\cite{irsic-2024} to be the polynomial $\Sigma _{i\ge 0} a_ix^i$, where $a_i$ is the number of distinct general position sets of $G$ with cardinality $i$. For any graph with order at least two, this polynomial will begin $1+nx+{n \choose 2}x^2+\cdots $. By setting  $\psi(X) = (1+x)^{|X|}$ for $X\subseteq V(G)$, this polynomial can be calculated by applying the inclusion/exclusion principle to the maximal general position sets of $G$ as follows.  

\begin{proposition}
\label{prop:inclusion-exclusion}
Let $G$ be a graph and let $X_1,\ldots,X_n$ be the maximal general position sets of $G$. Then 
$$\psi(G) = \sum_{k = 1}^n  (-1)^{k-1} \sum_{\{i_1,\ldots,i_k\}\subseteq [n]} \psi(X_{i_1} \cap \cdots \cap X_{i_k})\,.$$
\end{proposition}

For example, the article~\cite{irsic-2024} uses Proposition~\ref{prop:inclusion-exclusion} to show that $\psi(P) = 1 + 10x + 45x^2 + 90x^3 + 80x^4 + 30x^5+5x^6$, where $P$ is the Petersen graph. In~\cite{irsic-2024} this polynomial is determined for some simple families of graphs, such as complete graphs, paths, cycles and grid graphs. We include here the result for grid graphs.

\begin{theorem}{\rm \cite[Theorem 3.5]{irsic-2024}}
    \label{thm:grids-large}
    If $r, s \geq 3$, then 
    \begin{align*}
        \psi(P_r\cp P_s) & = 
        \frac{rs(r - 1)(r - 2)(s - 1)(s - 2)(r(s - 3) - s + 7)}{144} \; x^4 \\ & + \frac{1}{18}(r-1)r(s-1)s(r(2s-1)-s-4) \; x^3 \\ & + {rs \choose 2} x^2 + rs x + 1.
    \end{align*}
\end{theorem}

An interesting question is whether the general position polynomial determines a graph. This turns out not to be the case, with $C_4$ and $P_4$ as simple counterexamples. Less trivially, \cite{irsic-2024} constructs a family of trees with identical general position polynomials. The article also determines the general position polynomial of the disjoint union and join of two graphs in terms of the polynomials of the factors.

A common topic in graph polynomials is \emph{unimodality}. A polynomial $a_0+a_1x+a_2x^2+\ldots a_rx^r$ is unimodal if the sequence $a_0,a_1,\ldots ,a_r$ of coefficients increases to some maximum value and then decreases again. The polynomial is \emph{log-concave} if $a_k^2 \geq a_{k-1}a_{k+1}$ for $1 \leq k \leq r-1$, and log-concavity implies unimodality. The general position polynomial is in general not unimodal, and in fact this property fails to hold even for all trees (a broom graph on 24 vertices is used to witness this in~\cite{irsic-2024} and an infinite family of such brooms is constructed in~\cite{rather-2026+}). However, three families of graphs are shown to have unimodal polynomials, namely Kneser graphs $K(n,2)$, combs of paths and complete bipartite graphs $K_{r,r}$ with a perfect matching deleted.

The result for combs is extended in~\cite{rather-2026+} by proving that their general position polynomials are log-concave, although the same cannot be said for $K(n,2)$. The paper~\cite{rather-2026+} also reports that the general position polynomial of balanced complete multipartite graphs is log-concave (and hence unimodal) if the parts are of size at most four, but these properties do not hold generally for larger part size; for example $\psi (K_{8,8})$ is not unimodal. The paper~\cite{rather-2026+} also corrects a mistaken formula from~\cite{irsic-2024} for the general position polynomial of a family of graphs used in~\cite{TuiThoCha-2025} to give an estimate for the largest number of edges in a graph with given order and monophonic position number (see Subsection~\ref{subsec:mp}).

A question raised in~\cite{irsic-2024} is whether all corona graphs have unimodal general position polynomials. This was confirmed for some families in~\cite{rather-2026+}, although it was shown that they are not generally log-concave. A complete multipartite counterexample was found in~\cite{Wang2026} and later a complete bipartite counterexample in~\cite{Wang2026-2}.

\subsection{Graph products}\label{subsec:products}

In this subsection we examine the general position numbers of five types of graph product: Cartesian products in Subsubsection~\ref{subsubsec:Cartesian}, strong and direct products in Subsubsection~\ref{subsubsec:strongdirect}, generalised lexicographic products in Subsubsection~\ref{subsubsec:genlexicographic} and Sierpi\'nski products in Subsubsection~\ref{subsubsec:sierpinski}.

\subsubsection{Cartesian product}\label{subsubsec:Cartesian}

General position sets have been extensively researched on Cartesian products. As a first result, it was proved in~\cite[Corollary 3.2]{ManuelKlavzar-2018b} that $\gp(P_\infty \cp P_\infty) = 4$, where $P_\infty$ is the two-way infinite path. In the same paper it was further proved that $10\le \gp(P_\infty^{\cp, 3}) \le 16$, where $G^{\cp, n}$ denotes the $n$-fold Cartesian product of $G$. These results were completed by the next  result, where the Erd\H{o}s-Szekeres Theorem on monotone sequences was applied to prove the upper bound.  

\begin{theorem} {\rm \cite[Theorem~1]{KlavzarRus-2021}}
\label{thm:Cartesain-powers-of infinite paths}
If $n\in {\mathbb N}$, then $\gp(P_\infty^{\cp, n}) = 2^{2^{n-1}}$.
\end{theorem}

For general Cartesian products, we have the following lower bound. 

\begin{theorem} {\rm \cite[Theorem~3.1]{Ghorbani-2021}}
\label{thm:general-lower-Cartesian}
If $G$ and $H$ are connected graphs, then 
$$\gp(G \cp H) \geq \gp(G) + \gp(H) - 2\,.$$
\end{theorem}

Tian and Xu followed with a related result, where a graph is a {\em generalised complete graph} if it is obtained by the join of an isolated vertex with a disjoint union of one or more complete graphs.

\begin{theorem} {\rm \cite[Theorems~3.7 and 3.8]{Tian-2021a}}
\label{thm:general-lower-Cartesian-order}
If $G$ and $H$ are connected graphs with $n(G)\ge 3$ and $n(H)\ge 3$, then 
$$\gp(G \cp H) \le n(G) + n(H) - 2\,.$$
Moreover, the equality holds if and only if both $G$ and $H$ are generalised complete graphs.
\end{theorem}

The equality in the lower bound of Theorem~\ref{thm:general-lower-Cartesian} for the case when one factor is $K_2$ has been characterised when the second factor is bipartite. 

\begin{theorem} {\rm \cite[Theorem~3.5]{Tian-2021a}}
\label{thm:gp-prism-over-bipartite}
If $G$ is a connected bipartite graph with $n(G)\ge 3$, then $\gp(G \cp K_2) = \gp(G)$ if and only if G is a complete bipartite graph.
\end{theorem}

The equality case of the next result also demonstrates that the bound of Theorem~\ref{thm:general-lower-Cartesian} is sharp.

\begin{theorem} {\rm \cite[Theorem~3.2]{Ghorbani-2021}}
\label{thm:Hamming-gp}
If $k\ge 2$ and $n_1, \ldots, n_k\ge 2$, then 
$$\gp(K_{n_1}\cp \cdots \cp K_{n_k}) \ge n_1 + \cdots + n_k - k\,.$$  
Moreover, $\gp(K_{n_1}\cp K_{n_2}) = n_1 + n_2 - 2$.  
\end{theorem}

We have two exact results for the Cartesian products when trees are involved. 

\begin{theorem} 
\label{thm:main}
If $T$ and $T'$ are trees of order at least three, then the following hold. 
\begin{enumerate}
\item[(i)] {\rm \cite[Theorem~2.1]{Tian-2021b}} $\gp(T\cp T') = \gp(T) + \gp(T')$.
\item[(ii)] {\rm \cite[Theorem~3.2]{Tian-2021a}}
$\gp(T \cp K_m) = \gp(T) + \gp(K_m) - 1$.
\end{enumerate}
\end{theorem}

From Theorem~\ref{thm:Cartesain-powers-of infinite paths} we can deduce that if $r\ge 3$ and $s\ge 3$, then $\gp(P_r\cp P_s) = 4$. This result has been expanded as follows, where  $\#\gp(G)$ indicates the number of $\gp$-sets of a graph $G$. 

\begin{theorem} {\rm \cite[Theorem~2.1]{KlaPatRusYero-2021}}
\label{thm:gp-grids-enumerate}
If $2\le r\le s$, then
$$
\#\gp(P_r\cp P_s) = \left\{
\begin{array}{ll}
\vspace*{2mm}
6; & r = s = 2\,, \\
\vspace*{2mm}
\displaystyle{\frac{s(s-1)(s-2)}{3}}; & r = 2, s\ge 3\,, \\
\displaystyle{\frac{rs(r - 1)(r - 2)(s - 1)(s - 2)(r(s - 3) - s + 7)}{144}}; & r, s\ge 3\,.
\end{array}
\right.
$$
\end{theorem}

For cylinders we have: 

\begin{theorem} {\rm \cite[Theorem~3.2]{KlaPatRusYero-2021}}
\label{thm:gp-cylinders}
If $r\ge 2$ and $s\ge 3$, then
$$
\gp(P_r\cp C_s) = \left\{
\begin{array}{ll}
3; & r = 2, s = 3\,, \\
5; & r \ge 5,\ \text{and}\ s=7\ \text{or}\ s\ge 9\,, \\
4; & \text{otherwise}\,.
\end{array}
\right.
$$
\end{theorem}

In~\cite{KlaPatRusYero-2021} it was also proved that if $r\ge 3$ and $s\ge 3$, then $\gp(C_r\cp C_s)\leq 7$, 
and that if $r\ge s \ge 3$, $s\ne 4$, and $r\ge 6$, then $\gp(C_r\cp C_s) \in \{6, 7\}$. This research was subsequently concluded by Kor\v ze and Vesel in the next theorem.

\begin{theorem} {\rm \cite[Theorem~3.4]{KorzeVesel-2023}}
\label{thm:gp-cylindersKV}
If $r$ and $s$ are integers with $r\ge s\ge 4$, then
$$
\gp(C_r\cp C_s) = \left\{
\begin{array}{ll}
5; & s = 4\,, \\
6; & s \in \{5,6\} \ {\rm or}\ r,s\in \{8,10,12\}\,, \\
7; & \text{otherwise}\,.
\end{array}
\right.
$$
\end{theorem}

A related exact result is: 

\begin{theorem} {\rm \cite[Theorem~3.4]{Tian-2021a}}
\label{thm:gp-cycle-complete}
If $r\ge 4$, then 
$$
\gp(C_r\cp K_m) = \left\{
\begin{array}{ll}
2m; & m\le \lfloor s/2 \rfloor\,, \\
m + \left\lfloor \frac{r-1}{2}\right\rfloor; & m > \lfloor s/2 \rfloor\,.
\end{array}
\right.
$$
\end{theorem}

 A formula for the general position number of the Cartesian product of a complete multipartite graph with a path is also given in~\cite[Theorem~3.6]{Tian-2021a}. 

The $n$-dimensional hypercube $Q_n$, $n \ge 1$, is defined to be $K_2^{\cp, n}$, the Cartesian product of $n$ copies of $K_2$. To determine $\gp(Q_n)$ is a notoriously difficult problem. The only exact values known so far are $\gp(Q_1) = \gp(Q_2) = 2$, $\gp(Q_3) = 4$, $\gp(Q_4) = 5$, $\gp(Q_5) = 6$, $\gp(Q_6) = 8$, and $\gp(Q_7) = 9$, see~\cite[Table~1]{KorzeVesel-2023}. Using a probabilistic construction, K\"orner~\cite{Korner-1995} obtained general position sets in $Q_n$ of size $\frac{1}{2}\frac{2^n}{\sqrt{3^n}}$. He also noted that the problem of finding $\gp(Q_n)$ is equivalent to finding the largest size of a $(2,1)$-separating system in coding theory. Defining 
$$\alpha=\limsup_{n\rightarrow \infty}\frac{\log_2\gp(Q_n)}{n}\,,$$
the above construction yields $\alpha\ge 1-\frac{1}{2}\log_23$. K\"orner also proved that $\alpha\le 1/2$, while Randriambololona~\cite{randriambololona-2013} improved the lower bound to $\alpha\ge \frac{3}{50}\log_211$ with an explicit construction. Extending $\alpha$ to general Cartesian powers, the following quantity was introduced in~\cite{KlaPatRusYero-2021}: 
$$\gp_{\cp}(G):=\limsup_{n\rightarrow \infty}\frac{\log_{n(G)} \gp(G^{\cp, n})}{n}\,,$$
and the next result proved, where $p(G)$ is the probability that if one picks a triple $(x,y,z) \in V(G)^3$ uniformly at random, then $d_G(y,z) = d_G(y,x) + d_G(x,z)$ holds.

\begin{theorem} {\rm \cite[Theorem~5.1]{KlaPatRusYero-2021}}
\label{thm:random}
If $G$ is a graph, then
$$\gp_{\cp} (G)\ge \log_{n(G)}p(G)^{-1/2}\ge 1-\log_{n(G)}(n(G)^2-n(G)+1)\,.$$
\end{theorem}

\subsubsection{Strong product and direct product}\label{subsubsec:strongdirect}

For the strong product, we have the following general bounds. 

\begin{theorem} {\rm \cite[Theorem~4.2, Corollary~4.1]{Klavzar-2019}}
\label{thm:gp-strong-bounds}
If $G$ and $H$ are connected graphs, then 
$$\gp(G)\gp(H) \le \gp(G\boxtimes H) \le \min\{n(G)\gp(H), n(H)\gp(G)\}\,.$$
\end{theorem}

As noted in~\cite{Klavzar-2019}, if $n,m\ge 2$, then $\gp(P_n \boxtimes P_m) = 4$. This demonstrates sharpness of the lower bound of Theorem~\ref{thm:gp-strong-bounds}. Tightness of both bounds follows from the next result connecting the general position number of strong products with strong resolving graphs. 

\begin{proposition} {\rm \cite[Proposition~4.3]{Klavzar-2019}}
\label{prop:strong-one-factor-complete}
If $G$ is a connected graph and $n\ge 1$, then $\gp(G\boxtimes K_n) = n\cdot \gp(G)$. Moreover, if $\gp(G) = \omega(G\sr)$, then $\gp(G\boxtimes K_n) = \omega((G\boxtimes K_n)\sr)$.
\end{proposition}

The paper~\cite{Klavzar-2019} reports several additional exact results for the general position number of strong products, for example 
$$\gp(K_{r_1,t_1}\boxtimes K_{r_2,t_2})=r_1r_2=\omega((K_{r_1,t_1}\boxtimes K_{r_2,t_2})\sr) = \alpha(K_{r_1,t_1}\boxtimes K_{r_2,t_2})$$
for any $r_1\ge t_1\ge 1$ and any $r_2\ge t_2\ge 1$. The paper ends with the problem whether it is true that $\gp(G\boxtimes H) = \gp(G)\gp(H)$ holds for arbitrary connected graphs $G$ and $H$. This was shown to be false in~\cite{Vesel}; for example, it is proven that $\gp (C_{11k} \boxtimes C_{11k}) = 11$ for $k \geq 1$. However, the paper does show that $\gp(G\boxtimes H) = \gp(G)\gp(H)$ holds when one of the factors is a path or one of the cycles $C_4, C_5, C_6$, as well as the strong product $C_s \boxtimes C_t$ when $s,t \geq 3$, $s,t \not = 4$ and $t > 15(\left \lfloor \frac{s}{2} \right \rfloor -1)$. An upper bound for $\gp (C_s \boxtimes H)$ is also derived. It would be interesting to clarify when the general position number does behave multiplicatively over the strong product.

For the direct product of graphs, the only result found in the literature concerns the direct product of complete graphs. Notice that this result is another case where equality is achieved in Theorem~\ref{thm:gp-versus-strong-resolving-graphs}. 

\begin{proposition} {\rm \cite[Proposition~3.3]{Klavzar-2019}}
\label{prop:gp-direct-complete}
If $n_1\ge n_2\ge 3$, then $$\gp(K_{n_1}\times K_{n_2}) = \omega((K_{n_1}\times K_{n_2})\sr) = n_1 = \alpha((K_{n_1}\times K_{n_2})\sr)\, ,$$ unless $n_1 = n_2 = 3$, in which case $\gp(K_{3}\times K_{3}) = 4$.
\end{proposition}

We add that the formula of Proposition~\ref{prop:gp-direct-complete} was stated in~\cite{Klavzar-2019} for all $n_1 \geq n_2 \geq 3$, but it was pointed out by Ethan Shallcross that the case $K_{3}\times K_{3}$ is exceptional. 

\subsubsection{Generalised lexicographic product}\label{subsubsec:genlexicographic}

Let $G$ be a graph with $V(G) = \{g_1,\ldots, g_n\}$ and let $H_i$, $i\in [n]$, be pairwise disjoint graphs. Then the {\em generalised lexicographic product} $G[H_1,\ldots, H_n]$ is obtained from $G$ by replacing each vertex $g_i\in V(G)$ with the graph $H_i$, and each edge $g_ig_{j}\in E(G)$ with all possible edges between $H_i$ and $H_{j}$. If all the graphs $H_i$, $i\in [n]$, are isomorphic to a fixed graph $H$, then the generalised lexicographic product $G[H_1,\ldots, H_n] = G[H,\ldots, H]$ coincides with the standard lexicographic product $G\circ H$.

With the use of Theorems~\ref{thm:gpsets-characterised} and~\ref{thm:gp-versus-strong-resolving-graphs}, the next result was proved. 

\begin{theorem} {\rm \cite[Theorem~5.1]{Klavzar-2019}}
\label{thm:blow-up}
Let $G$ be a graph with $V(G) = \{g_1,\ldots, g_n\}$ and let $k_i$, $i\in [n]$, be positive integers. If $S$ is a gp-set of $G$ that induces a complete subgraph of $G\sr$, and  $\min\{k_i:\ g_i\in S\} \ge \max\{k_i:\ g_i\notin S\}$, then
$$\gp(G[K_{k_1},\ldots, K_{k_n}]) = \sum_{i:g_i\in S}k_i = \omega((G[K_{k_1},\ldots, K_{k_n}])\sr)\,.$$
\end{theorem}

\subsubsection{Sierpi\'nski product}\label{subsubsec:sierpinski}

Motivated by the classical Sierpi\'{n}ski graphs~\cite{hinz-2017}, the Sierpi\'{n}ski product of graphs was introduced in~\cite{kovic-2023} as follows. Let $G$ and $H$ be graphs and let $f \colon V(G)\rightarrow V(H)$ be a function. The {\em Sierpi\'{n}ski product} of $G$ and $H$ (with respect to $f$) is the graph  $G \otimes _f H$ with vertices $V(G \otimes _f H) = V(G)\times V(H)$, and with edges
\begin{itemize}
\item $(g,h)(g,h')$, where $g\in V(G)$ and $hh' \in E(H)$, and
\item $(g,f(g'))(g',f(g))$, where $gg' \in E(G)$.
\end{itemize}
In~\cite{TianKlavzar-2024++}, the \textit{Sierpi\'{n}ski general position number} $\gps(G,H)$, and the {\em lower Sierpi\'nski general position number} $\gpsl(G,H)$, were introduced as the cardinality of a largest, resp.\ smallest, $\gp $-set in $G \otimes _f H$ over all possible functions $f \colon V(G)\rightarrow V(H)$, that is, 
\[ \gps(G, H) = \max_{f\in H^G}\{\gp(G\otimes _f H)\} \qquad
{\rm and} \qquad
\gpsl(G, H) = \min_{f\in H^G}\{\gp(G\otimes _f H)\}\,.
\]
To deal with them, a few new concepts have been introduced. If $u\in V(G)$ and $S\subseteq V(G)$, then $S$ is a {\em $u$-colinear set} if $S$ is a general position set such that $u\notin S$ and $y\notin I_G[x,u]$ for any $x,y\in S$. In the terminology of Section~\ref{subsec:vertex position}, we can say that $S$ is a $u$-colinear set if $S$ is both a general position set and a $u$-position set. Also set 
\begin{align*}
\xi_G(u) & = \max \{ |S|:\ S\ {\rm is\ a}\ u\mbox{-}{\rm colinear\ set}\}\,, \\
\xi^{-}(G) & = \min \{\xi_G(u):\ u\in V(G)\}\,,\\
\xi(G) & = \max \{\xi_G(u):\ u\in V(G)\}\,.
\end{align*}
The hierarchy of these invariants on an arbitrary graph is given in the next result. 

\begin{theorem}{\rm \cite[Theorem~3.2]{TianKlavzar-2024++}}
\label{thm:2xi>gp}
If $G$ is a connected graph of order at least two and $u\in V(G)$, then
$$\xi^{-}(G) \le \xi_G(u) \le \xi(G)\le \gp(G) \leq 2\xi^{-}(G)\,.$$
\end{theorem}

The next theorem gives a general bound for Sierpi\'{n}ski products. 

\begin{theorem}{\rm \cite[Theorem~4.2]{TianKlavzar-2024++}}
\label{thm:Sierpinski-product-bounds}
If $G$ and $H$ are two connected graphs of order at least two, then
\[
\gp(H) \leq \gpsl(G,H)\leq \gps(G,H)\leq n(G)\gp(H).
\]
Moreover, $\gps(G,H) = n(G)\gp(H)$ if and only if $\gp(H) = \xi(H)$.
\end{theorem}
If $G=K_2$ and $H$ is a connected graph with $n(H)\ge 2$, then $\gpsl(K_2,H) = 2\xi^{-}(H)$ and $\gps(K_2,H) = 2\xi(H)$, see~\cite[Theorem~4.4]{TianKlavzar-2024++}. Finally, for Sierpi\'{n}ski products of complete graphs this result applies. 

\begin{theorem}{\rm \cite[Theorems~5.1 and 5.3]{TianKlavzar-2024++}}
If $m,n\geq 2$, then 
\begin{enumerate}
\item[(i)] $\gps(K_m,K_n) = m(n-1)$, and 
\item[(ii)] if in addition $n\ge 2m-2$, then $\gpsl(K_m,K_n) = m(n-m+1)$. 
\end{enumerate}
\end{theorem}

\subsection{Other operations}\label{subsec:otheroperations}

\subsubsection{Vertex- and edge-deleted subgraphs}
\label{sec:vertex-edge-delete}

For vertex-deleted subgraphs we have this result. 

\begin{theorem}{\rm \cite[Theorem~3.1]{Dokyeesun-2025}}
\label{thm:gp(g-x)-at-most}
If $x$ is a vertex of a graph $G$, then $\gp(G-x)\leq 2\gp(G)$. Moreover, the bound is sharp.
\end{theorem}

The sharpness examples for Theorem~\ref{thm:gp(g-x)-at-most} provided in~\cite{Dokyeesun-2025} are such that $G-x$ is not connected. Constructions are also given which demonstrate that $\gp(G-x)$ can be much larger than $\gp(G)$ even when $G-x$ is connected.

On the other hand, $\gp(G-x)$ cannot be bounded from below by a function of $\gp(G)$. For instance, consider the fan graphs $F_n$, $n\ge 3$, that is, the join between $P_n$ and $K_1$. Then $\gp(F_n) = \lceil \frac{2n}{3}\rceil$ and $\gp(F_n - x) = 2$, where $x$ is the vertex of $F_n$ of degree $n$. On the positive side, if $x$ lies in some $\gp $-set, then $\gp(G) - 1 \le \gp(G-x)$, see~\cite[Proposition~3.3]{Dokyeesun-2025}. For edge-deleted subgraphs we have: 

\begin{theorem}
\label{thm:edge-removal}
If $e$ is an edge of a graph $G$, then
$$ \frac{\gp(G)}{2}\le \gp(G-e)\leq\ 2\gp(G)\,.$$
Moreover, both bounds are sharp.
\end{theorem}

Stronger results than those above for vertex- and edge-deleted subgraphs are produced for graphs of diameter two. 

\subsubsection{Complementary prisms}

If $G$ is a graph and $\overline{G}$ its complement, then the {\em complementary prism} $G\overline{G}$ of $G$ is the graph formed from the disjoint union of $G$ and $\overline{G}$ by adding the edges of a perfect matching between the corresponding vertices of $G$ and $\overline{G}$, see~\cite{haynes-2007}. The paper~\cite{Neethu-2021} focuses on the general position number of complementary prisms, of which the following are some of the most significant results. 

\begin{theorem} {\rm \cite[Theorem~3.1]{Neethu-2021}}
\label{thm:gp-complementary-bounds} 
If $G$ is a connected graph, then $\gp(G\overline{G})\leq n(G)+1$, and if $G$ is disconnected, then $\gp(G\overline{G})\leq n(G)$.
\end{theorem}

The next result characterises the graphs that achieve the first upper bound of Theorem~\ref{thm:gp-complementary-bounds}. To state it, we need a couple of new definitions from~\cite{Neethu-2021}.  A set $S$ of vertices in a graph $G$ is a {\em $3$-general  position set} if no three vertices from $S$ lie on a common shortest path of length at most three (this was the original idea that led to the $d$-positions sets which will be considered in Section~\ref{subsec:d-position-sets}). If $u \in V(G)\cap V(G\overline{G})$, then its unique neighbour in 
$V(G\overline{G})\cap V(\overline{G})$ is denoted by $\overline{u}$ and is called the {\em partner} of $u$ in $G$. 
For a set $X \subseteq V(G)$, the set of partners of the vertices of $X$ is denoted by $\overline{X}$. Recall that a vertex $v$ of a connected graph $G$ is a \emph{central vertex} of $G$ if the eccentricity of $v$ is equal to the radius of $G$. By $N_{G}^{2}(v)$ we denote the set of vertices at distance two from $v$ in $G$.

\begin{theorem}  {\rm \cite[Theorem~3.4]{Neethu-2021}}
\label{thm:gp-complementary-sharpness}
Let $G$ be a graph with $n(G)\ge 4$ and such that both $G$ and $\overline{G}$ are connected. Then $\gp(G\overline{G})=n(G)+1$ if and only if $\rad(G)= 2$ and there exists a central vertex $v$ such that
\begin{enumerate}
\item[(i)] $\overline{N_{G}(v)}$ is a $3$-general position set in $\overline{G}$ and $N_{G}^{2}(v)$ is a $3$-general position set in $G$, and 
\item[(ii)] for each $x\in N_G(v)$ there exists $y\in N_{G}^{2}(v)$ such that $xy\notin E(G)$. 
\end{enumerate}
\end{theorem}

Setting 
$$\overline{\gp}_{3}(G) = \max \{ \gp_{3}(\overline{G}[V(\overline{G})\setminus \overline{S}]):\ S \text{ is a } \gp_{3}{\mbox -}{\rm set\ of}\ G\}\,,$$
we have the following lower bound on $\gp(\GG)$. 

\begin{theorem} {\rm \cite[Theorem~4.1]{Neethu-2021}}
\label{thm:gp-complementary-lower}
If $G$ is a graph, then 
$$\gp(\GG) \ge \max \{\gp_{3}(G) + \overline{\gp}_{3}(G),  \gp_{3}(\overline{G}) + \overline{\gp}_{3}(\overline{G})\}\,.$$
Moreover, the bound is sharp. 
\end{theorem}

It follows from Theorem~\ref{thm:gp-complementary-lower} that $\gp(G\overline{G})\geq \max\{\gp_{3}(G),\gp_{3}(\overline{G})\}$. In~\cite[Theorem~4.2]{Neethu-2021} it was proved that among connected graphs the equality holds here if and only if $G$ is a complete multipartite graph. The paper~\cite{Neethu-2021} considers several additional classes of graphs $G$ and the corresponding values of $\gp(\GG)$: bipartite graphs (in particular trees, hypercubes and Cartesian products of paths), split graphs and block graphs. For instance, if $T$ is a tree, then  
$$\gp(T\overline{T}) = 
\begin{cases}
 n(G)+1; & \diam(T) = 4\,, \\
 n(G); & otherwise\,,
 \end{cases}
$$
and if $n\geq 2$, then $\gp(Q_n\overline{Q}_n) = 2^n$.

\subsubsection{Rooted products}

By a \emph{rooted graph} we mean a connected graph having one fixed vertex called the \emph{root} of the graph. Consider now a connected graph $G$ of order $n$, and let $H$ be a rooted graph with root $v$. The \emph{rooted product graph} $G\circ_v H$ is the graph obtained from $G$ and $n$ copies of $H$, say $H_1,\dots,H_n$, by identifying the root of $H_i$ with the $i^{\rm th}$ vertex of $G$. 

\begin{theorem} {\rm \cite[Theorem~6.1]{Klavzar-2019}}
\label{thm:rooted}
Let $G$ be any connected graph of order $n\ge 2$, and let $H$ be a rooted graph with root $v$.
\begin{itemize}
  \item[(i)] $\gp(G\circ_v H)=n=\omega((G\circ_v H)\sr)$ if and only if $H$ is a path and $v$ is a leaf of $H$.
  \item[(ii)] If $H$ contains a $\gp $-set $S$ not containing $v$ and such that for each pair of vertices $u,w\in S$ neither $u\in I_H(v,w)$ nor $w\in I_H(v,u)$, then $\gp(G\circ_v H)=n\cdot \gp(H)$. Moreover, if in addition $S$ is a maximum clique in $H\sr$, then $\gp(G\circ_v H)=\omega((G\circ_v H)\sr)$.
  \item[{\rm (iii)}] Suppose $H$ is not a path rooted in one of its leaves. If every $\gp $-set $S$ of $H$ either contains the root $v$, or contains two vertices $x,y$ such that ($x\in I_H(v,y)$ or $y\in I_H(v,x)$), then $2n\le \gp(G\circ_v H)\le n(\gp(H)-1)$. In particular, if every $\gp $-set of $H$ contains the root $v$, then $\gp(G\circ_v H)=n(\gp(H)-1)$.
\end{itemize}
\end{theorem}

\subsubsection{Powers of graphs} 
If $G$ is a graph and $k \geq 1$, then the $k$-th power $G^k$ of $G$ is the graph with vertex set $V(G)$, where two vertices are adjacent in $G^k$ if and only if their distance in $G$ is at most $k$. The paper~\cite{powers} investigates the general position number of graph powers. A lower bound for $\gp(G^k)$ for an arbitrary graph $G$ was given in~\cite{powers}.

\begin{theorem}{\rm \cite[Theorem 3.1]{powers}}\label{power3.1} 
If $k \geq 2$ and $G$ is a graph with $n(G)\ge 3$, then $\gp(G^k) \geq \Delta(G) + k - 1$. 
 \end{theorem}
 
 In~\cite[Theorem 3.2]{powers} the result on paths is extended to $\gp(P_n^k) = k + 1$ for $n \geq 4$. This demonstrates the sharpness of the bound in Theorem~\ref{power3.1}.  As noted in~\cite{powers}, for \( n \leq 2k+1 \), \( C_n^k \cong K_n \), and so \( \gp(C_n^k) = n \), whilst for larger $n$ we have the following result~\cite[Lemma 3.3, Theorem 3.4]{powers}.
 
 \begin{theorem}{\rm \cite{powers}} Let \( n \geq 6 \). Then \( k+1 \leq \gp(C_n^k) \leq 2k + 1 \), with equality on the left side if and only if \( n \leq 3k + 1 \).
Moreover, for $n\geq 3k+2$, we have 
$$
\gp(C_n^k) = \left\{
\begin{array}{ll}
2k+1; & \gcd(n-1, \left\lfloor \frac{n}{2}\right\rfloor)\equiv 0\pmod{k}\,, \\
2k; & \gcd(n, \left\lfloor \frac{n}{2}\right\rfloor)\equiv 0\pmod{k}\,.
\end{array}
\right.
$$
\end{theorem}

In~\cite[Proposition 3.5]{powers} it was also proved that if $n\geq 3k+2$ and $n\equiv 0\pmod{k}$, then $\gp(C_n^k)\leq 2k$. The paper~\cite{powers} further focuses on the general position number of the square of graphs. In the next result, graphs with $\gp(G^2) = 3$ are characterised.

\begin{theorem}
  {\rm \cite[Lemma 4.1]{powers}}\label{gp3_square} Let $G$ be a graph of order $n \geq 3$. Then $\gp(G^2) = 3$ if and only if
 $G\in\{P_n,C_6,C_7\}$.  
\end{theorem}
Consequently, the general position numbers of the square of cycles are derived in~\cite[Theorem 4.4]{powers}. 
\begin{theorem} {\rm \cite[Theorem 4.4]{powers}}
If $n\geq 8$, then 
$$\gp(C_n^2) = \left\{
\begin{array}{ll}
4; & n\ \text{is even or} \ n=11\,, \\
5; & \text{otherwise}\,.
\end{array}
\right.$$
\end{theorem}

The paper~\cite{powers} then presents some additional results for the square of block graphs. For instance, $\gp(G^2) \leq \gp(G) + 1$ for any block graph \( G \), as shown in~\cite[Lemma 4.6]{powers}. In addition, the next theorem provides a characterisation of block graphs for which equality holds.

\begin{theorem}{\rm \cite[Theorem 4.8]{powers}}\label{square_block} 
Let $G$ be a block graph with $\delta(G) \geq 2$ and $S$ be the set of simplicial
 vertices of $G$. Then $\gp(G^2) = \gp(G) + 1$ if and only if there exists a cut-vertex $u$ such that $d_G(u,s)$ is odd for every vertex $s \in S$.
 \end{theorem}    

In particular, for the square of trees we have: 

\begin{theorem}{\rm\cite[Theorem 4.9]{powers}}  
If $T$ is a tree with $n(T)\ge 3$, then $\gp(T^2) \in \{\gp(T),\gp(T)+1\}$.  Moreover, $\gp(T^2) = \gp(T)$ if and only if $T$ has two non-trivial blocks $Q,Q'$ such that $d(Q,Q')$ is odd.  
\end{theorem}

\subsubsection{Mycielskian graphs, double graphs and shadow graphs}

The \emph{Mycielskian} $\mathcal{M}(G)$ of a graph $G = (V,E)$ is the graph with vertex set $V \cup V' \cup \{ u^*\} $, where $V = \{ u_1,\dots ,u_n\} $ is the vertex set of $G$, $V' = \{ u_1',\dots,u_n'\} $ is a copy of $V$, and $u^*$ is the \emph{root vertex}. The edge set of $\mathcal{M}(G)$ is defined to be $\{ uv:uv \in E\} \cup \{ uv':uv \in E\} \cup \{ v'u^*:v'\in V'\}$. This construction was introduced by Mycielski as a means of finding triangle-free graphs with arbitrarily large chromatic number, but it has been studied in its own right in many papers.

The general position number of Mycielskians was first studied in~\cite{Thomas-2024a}. It is shown that there is a $\gp $-set of $\mathcal{M}(G)$ containing the root vertex $u^*$ if and only if $G$ is a clique or $G$ is the join of $K_1$ with a disjoint union of $t \geq 2$ cliques $W_1,\dots ,W_t$, where one of the cliques has order one and all of the other cliques have order at least three. Furthermore the root $u^*$ belongs to every $\gp $-set only if $G$ is a complete graph, in which case $\gp (\mathcal{M}(G)) = n+1$ and $V \cup \{ u^*\} $ is the unique $\gp $-set.  Hence for non-complete graphs it suffices to examine $\gp $-sets not containing the root, and such $\gp $-sets are studied in terms of associated partitions of $V$ into four parts. 

A simple lower bound is $\gp (\mathcal{M}(G)) \geq \max \{ n,2\ip _4(G)\} $, where $\ip _4(G)$ is the number of vertices in a largest independent set $S$ of $G$ such that no shortest path of length at most four passes through three vertices of $S$. Graphs that meet this bound are called \emph{meagre} and graphs for which $\gp (\mathcal{M}(G)) > \max \{ n,2\ip _4(G)\}$ are called \emph{abundant}. For $n \geq 3$ the largest size of an abundant graph is ${n \choose 2}+1$ and the abundant graphs that achieve this extremal size are an $(n-1)$-clique with a leaf attached and the gem graph $K_1 \vee P_4$. Complete multipartite graphs are meagre, as are paths, all cycles apart from $C_3$ and $C_5$, all cubic graphs with two exceptions and all sufficiently large regular graphs. An upper bound on $\gp (\mathcal{M}(G))$ is $\gp(\mathcal{M}(G)) \leq n+\max\{ 0,\ip _4(G)-\delta(G) +1\} $, so that a weaker upper bound is $n+\alpha(G) -1$; the graphs meeting these upper bounds are characterised. 

The case of regular graphs is particularly interesting. If $G$ is a $d$-regular graph, then $\mathcal{M}(G) \leq n+\left \lfloor \frac{d-1}{2}+\frac{1}{d} \right \rfloor $. All $d$-regular graphs with order $\geq d^3-2d^2+2d+2$ are meagre, so that there are only finitely many abundant $d$-regular graphs. A $d$-regular abundant graph of order $n = 3d-1$ and with $\gp (\mathcal{M}(G)) = n+1$ is constructed, and this turns out to be the unique abundant $d$-regular graph for $d \in \{2,3\}$. It is an open question whether this construction is always the unique abundant $d$-regular graph and, if not, how large the difference $\gp (\mathcal{M}(G))-n$ can be for $d$-regular graphs.  

Trees and graphs with large girth are also studied in~\cite{Thomas-2024a}. If we denote the number of support vertices by $\sigma $, then $\gp(\mathcal{M}(T)) = n(T)+\ell (T)-\sigma (T)$, and if all support vertices are at distance at least three apart, then we have equality. If every vertex of a tree $T$ is either a leaf or a support vertex, then $\gp (\mathcal{M}(T)) = 2\ell (T)$ and $T$ is meagre. If $G$ has girth at least six and matching number $\nu(G) $, then $\gp(\mathcal{M}(G)) \leq 2n-2\nu (G)$. Hence any graph $G$ with girth at least six that contains a perfect matching has $\gp (\mathcal{M}(G)) = n$. 

The article~\cite{Roy-2025} is mainly concerned with mutual visibility, but does contain a section on the general position problem in double graphs. The \emph{double graph} $\mathcal{D}(G)$ of $G$ is formed from the disjoint union of two copies of $G$ (the vertex sets of which we will denote by $V$ and $V'$, as for Mycielskians) by adding an edge $uv'$ whenever $uv$ is an edge in $G$ (where $v'$ is the copy of $v$ belonging to $V'$). It is shown that the general position number of $\mathcal{D}(G)$ lies in the range $\gp (G) \leq \gp (\mathcal{D}(G)) \leq 2\gp (G)$. If the upper bound is met, then the $\gp $-sets of $\mathcal{D}(G)$ have the form $S \cup S'$, where $S$ is a $\gp $-set of $G$ that is also an independent set. The upper bound is sharp for paths and cycles, and the lower bound is sharp for $K_n$ and $K_n^-$. It would be interesting to determine the graphs for which the lower bound is attained.

The complexity of finding the $\gp $-number of Mycielskian and double graphs remains unknown. 

A graph class closely related to the construction of Mycielskian graphs is that of \emph{shadow graphs}. The \emph{shadow graph} $S(G)$ of a graph $G$ is obtained from the Mycielskian $\mathcal{M}(G)$ by removing the root vertex $u^*$. The article \cite{haritha-2026} studies the general position and mutual visibility problems on shadow graphs. A lower bound on $\gp(S(G))$ is given by $\gp(S(G)) \geq 2\ip (G)$, where $\ip (G)$ denotes the number of vertices in a largest independent general position set of $G$, see Subsection~\ref{subsec:ipset}. This bound is sharp for complete bipartite graphs, trees, and cycles $C_n$ $(n \geq 8)$.
An upper bound on $\gp(S(G))$ is given by
\[
\gp(S(G)) \leq n(G) + \min\Bigg\{ \ip(G) - \delta(G) + 1,\,
\left\lfloor \frac{\ip(G)\,(n(G)-1) - \delta(G)}{\ip(G) + \delta(G)} \right\rfloor \Bigg\}.
\]
   This upper bound is sharp for many classes of graphs, including complete graphs. Moreover, the bound can be further improved for regular triangle-free graphs to $\gp(S(G)) \le n(G)$, and this improved bound is sharp for all regular triangle-free graphs of diameter at most three.

\section{Distance Based Variations on General Position}\label{sec:distancevariations}

In the next three sections we give an overview of variations of the general position problem that have been studied in the literature. In the present section, we concentrate on variations of the general position problem that, like the ordinary general position problem, are based on distance. In Section~\ref{sec:intro} we mentioned that an interesting source of variants of the general position problem is to change the family of paths that should not contain three-in-a-line; in Subsection~\ref{subsec:d-position-sets} we examine the effect of restricting the path family to shortest paths of bounded length. In Subsection~\ref{subsec:ipset} we discuss general position sets that are required to have the additional property of being an independent set. Subsection~\ref{subsec:vertex position} looks at a local version of the general position problem, i.e.\ how many vertices are visible along shortest paths from some fixed vertex. It is trivial that any clique is in general position; we can extend this by observing that if a set $S$ has the property that for some $k \geq 1$ all vertices in $S$ are at distance $k$ from each other, then $S$ is in general position. This leads us to the idea of an \emph{equidistant set}, and this is the subject of Subsection~\ref{subsec:equidistant}. The traditional general position problem for graphs asks for vertex subsets $S$ such that any pair of vertices from $S$ is $S$-positionable; however, if we impose the $S$-positionability requirement on other pairs of vertices we obtain the \emph{variety of general position problems} discussed in Subsection~\ref{subsec:variety}. The \emph{edge general position problem}, addressed in Subsection~\ref{subsec:edgegp}, arises if we require that no $k$ edges are contained in a common shortest path. Subsection~\ref{subsec:steiner} discusses the \emph{Steiner general position problem}, which arises if we move away from the ordinary graph distance and instead use the Steiner distance for a set of vertices. Finally, Subsection~\ref{subsec:detour} reports on one paper that uses longest paths instead of shortest paths.

 An illustration of some of the different types of position sets discussed in this section can be seen in the Petersen graph in Figure~\ref{fig:Petersen}. 

	\begin{figure}[ht!]
		\centering
		\begin{tikzpicture}[x=0.2mm,y=-0.2mm,inner sep=0.2mm,scale=0.6,thick,vertex/.style={circle,draw,minimum size=10}]
			\node at (180,200) [vertex,fill=red] (v1) {};
			\node at (8.8,324.4) [vertex] (v2) {};
			\node at (74.2,525.6) [vertex,fill=red] (v3) {};
			\node at (285.8,525.6) [vertex,fill=red] (v4) {};
			\node at (351.2,324.4) [vertex] (v5) {};
			\node at (180,272) [vertex,fill=red] (v6) {};
			\node at (116.5,467.4) [vertex] (v7) {};
			\node at (282.7,346.6) [vertex,fill=red] (v8) {};
			\node at (77.3,346.6) [vertex,fill=red] (v9) {};
			\node at (243.5,467.4) [vertex] (v10) {};

			\node at (580,200) [vertex] (u1) {};
			\node at (408.8,324.4) [vertex,fill=green] (u2) {};
			\node at (474.2,525.6) [vertex] (u3) {};
			\node at (685.8,525.6) [vertex] (u4) {};
			\node at (751.2,324.4) [vertex,fill=green] (u5) {};
			\node at (580,272) [vertex,fill=green] (u6) {};
			\node at (516.5,467.4) [vertex] (u7) {};
			\node at (682.7,346.6) [vertex] (u8) {};
			\node at (477.3,346.6) [vertex] (u9) {};
			\node at (643.5,467.4) [vertex] (u10) {};

            \node at (980,200) [vertex] (w1) {};
			\node at (808.8,324.4) [vertex,fill=blue] (w2) {};
			\node at (874.2,525.6) [vertex] (w3) {};
			\node at (1085.8,525.6) [vertex] (w4) {};
			\node at (1151.2,324.4) [vertex,fill=blue] (w5) {};
			\node at (980,272) [vertex] (w6) {};
			\node at (916.5,467.4) [vertex,fill=blue] (w7) {};
			\node at (1082.7,346.6) [vertex] (w8) {};
			\node at (877.3,346.6) [vertex] (w9) {};
			\node at (1043.5,467.4) [vertex,fill=blue] (w10) {};
			\path
			(v1) edge (v2)
			(v1) edge (v5)
			(v2) edge (v3)
			(v3) edge (v4)
			(v4) edge (v5)
			
			(v6) edge (v7)
			(v6) edge (v10)
			(v7) edge (v8)
			(v8) edge (v9)
			(v9) edge (v10)
			
			(v1) edge (v6)
			(v2) edge (v9)
			(v3) edge (v7)
			(v4) edge (v10)
			(v5) edge (v8)

			(u1) edge (u2)
			(u1) edge (u5)
			(u2) edge (u3)
			(u3) edge (u4)
			(u4) edge (u5)
			
			(u6) edge (u7)
			(u6) edge (u10)
			(u7) edge (u8)
			(u8) edge (u9)
			(u9) edge (u10)
			
			(u1) edge (u6)
			(u2) edge (u9)
			(u3) edge (u7)
			(u4) edge (u10)
			(u5) edge (u8)

            (w1) edge (w2)
			(w1) edge (w5)
			(w2) edge (w3)
			(w3) edge (w4)
			(w4) edge (w5)
			
			(w6) edge (w7)
			(w6) edge (w10)
			(w7) edge (w8)
			(w8) edge (w9)
			(w9) edge (w10)
			
			(w1) edge (w6)
			(w2) edge (w9)
			(w3) edge (w7)
			(w4) edge (w10)
			(w5) edge (w8)
			
			;
		\end{tikzpicture}
		\caption{The Petersen graph with different types of position set. On the left in red: a largest general position set. In the centre in green: a largest monophonic position set [Section~\ref{subsec:mp}]. On the right in blue: a largest equidistant set [Section~\ref{subsec:equidistant}], independent position set [Section~\ref{subsec:ipset}] and mobile general position set [Section~\ref{subsec:mobilegp}], also a lower general position set [Section~\ref{subsec:lower}].}
		\label{fig:Petersen}
	\end{figure}

\subsection{$d$-position sets}\label{subsec:d-position-sets}

If we restrict attention to the family of shortest paths with bounded length, we obtain the following problem that was first investigated in~\cite{KlaRalYer-2021}. The general $d$-position problem asks for the largest number of vertices in a set $S$ such that if three vertices of $S$ lie on a common shortest path $P$, then $P$ has length greater than $d$. The cardinality of a largest such set is denoted by $\gp _d(G)$. For any graph $G$ we have the inequality chain
\[ \gp (G) = \gp _{\diam (G)}(G) \leq \gp _{\diam (G)-1}(G) \leq \dots \leq \gp _2(G) \]
and various examples are given in~\cite{KlaRalYer-2021} to show that it is possible to have equality throughout this chain, strict inequality in exactly one case or strict inequality throughout.

The article~\cite{KlaRalYer-2021} also gives a characterisation of the structure of general $d$-position sets along the lines of Theorem~\ref{thm:gpsets-characterised} and shows that for $d \geq 2$ the general $d$-position problem is NP-complete. The general $d$-position numbers of paths and cycles are determined, and it is also shown that any infinite graph has a general $d$-position set of infinite cardinality. Some connections with strong resolving graphs, the dissociation number and independence number are also explored.

Two interesting questions left open in~\cite{KlaRalYer-2021} are to find the general $d$-position numbers of all grid graphs and to determine the complexity of finding the general $d$-position number of trees for any $d$ (it is known that it can be found in polynomial time for $d = 2$ or $d = \diam (T)$).

A generalised version of this problem is treated in~\cite{Cody-2025}. The $k$-general $d$-position number of $G$, denoted $\gp _d^k(G)$, is the largest possible number of vertices in a set $S \subseteq V(G)$ such that no $k$ vertices of $S$ lie on a common shortest path of length at most $d$. Instead of the above chain of inequalities, one now has a two-dimensional lattice of inequalities in $d$ and $k$. The article gives bounds for this number in terms of the $k$-general $d$-position numbers of isometric subgraphs of $G$ that are sufficiently far apart. The article determines these numbers for paths, cycles and thin grid graphs. Finding the $k$-general $d$-position number of other Cartesian products is left as an open problem, along with finding the computational complexity of the $k$-general $d$-position problem and the problem of characterising the structure of such sets. The $k$-general $d$-position number of powers of paths are determined in~\cite{Cody-2026}.

\subsection{Independent general position sets}\label{subsec:ipset}
The article~\cite{Thomas-2021} studies the properties of vertex subsets of graphs that are simultaneously independent sets and in general position. Such subsets are called \emph{independent position sets} and the largest order of an independent position set is the \emph{independent position number} of $G$, denoted by $\ip (G)$. Naturally for any graph $\ip (G) \leq \max \{ \gp(G),\alpha(G)\} $.

For some families of graphs, including trees, we have equality $\gp (G) = \ip (G)$. However, it is easily seen that $\gp (G)- \ip (G)$ can be arbitrarily far apart (for example complete graphs). Realisation results from~\cite{Thomas-2021} demonstrate that for any $2 \leq a \leq b < n$ there exists a graph with $\ip (G) = a$, $\gp (G) = b$ and order $n$.  It is shown in~\cite{Thomas-2021} that for any graph with diameter at most three we will have $\ip (G) = \alpha (G)$, since any independent set is in general position. We also have equality with the independence number for bipartite graphs. The extreme values $\ip (G) = 1$ and $\ip (G) = n(G)$ are achieved only by cliques and null graphs respectively; amongst connected graphs the largest possible value of $\ip (G)$ is $n(G)-1$ and this is achieved only for star graphs.  

The independent general position numbers of three graph products are discussed in~\cite{Thomas-2021}, the Cartesian product, the lexicographic product and the corona product. For the Cartesian product, the lower bound for the independent position number is analogous to Theorem~\ref{thm:general-lower-Cartesian} for the ordinary general position number.

\begin{theorem} {\rm \cite[Theorem~3.1]{Thomas-2021}}
\label{thm:ip-Cartesian}
 For any connected graphs $G$ and $H$, \[ \ip (G \cp H) \geq \ip (G) +\ip (H) -2.\]
\end{theorem}
It appears to be unknown if this bound is tight, so perhaps it could be improved. For the lexicographic product, the lower bound is as follows.

\begin{theorem} {\rm \cite[Theorem~3.3]{Thomas-2021}}
\label{thm:ip-lexicographic}
 For any connected graphs $G$ and $H$, \[ \ip (G\circ H) \geq \ip (G)\ip (H).\]
\end{theorem}

The bound in Theorem~\ref{thm:ip-lexicographic} is tight, as can be seen by considering $C_4\circ K_2$. It would be of interest to characterise the case of equality, and to provide upper bounds for Cartesian and lexicographic products.

Finally, the exact value of the independent position number is known for corona products.
\begin{theorem} {\rm \cite[Theorem~3.4]{Thomas-2021}}
\label{thm:ip-corona}
 For any connected graphs $G$ and $H$, \[ \ip (G \odot H) = n(G)\alpha (H).\]
\end{theorem}
Independent position sets are also briefly considered in the context of graph colouring~\cite{ChaDiSreeThoTui2024+}, see Subsection~\ref{subsec:colouring}. 

\subsection{Vertex position sets}\label{subsec:vertex position}

\emph{Vertex position sets}, introduced in~\cite{ThaChaTuiThoSteErs-2024}, are a local version of general position sets that measure how many vertices are visible from a fixed vertex of a graph. It was initially inspired by the classical result of analytic number theory on the density of points on an integer lattice that are visible from the origin, as well as art gallery theorems.

For a fixed vertex $x$ of a graph $G$, an \emph{$x$-position set} is a subset $S_x \subseteq V(G)$ such that for any $y\in S_x$ no shortest $x,y$-path contains a vertex of $S_x \setminus \{ y\}$, and the \emph{$x$-position number} of $G$ is the number of vertices in a largest $x$-position set. According to this definition, $x$ does not belong to any $x$-position set of order $\geq 2$. Naturally, the $x$-position number of $G$ varies from vertex to vertex. The maximum value of this position number over all vertices of $G$ is the \emph{vertex position number} $\vp (G)$ of $G$, whereas the minimum value is the \emph{lower vertex position number} $\vp ^-(G)$ of $G$. Note that in this definition the word `lower' is used in a slightly different sense to the `lower position number' discussed in Subsection~\ref{subsec:lower}.

By choosing a vertex inside a $\gp $-set of $G$, it can be seen that $\vp (G) \geq \gp(G)-1$. It is also shown in~\cite{ThaChaTuiThoSteErs-2024} that the numbers are bounded in terms of the vertex degrees as follows: $\vp ^-(G) \geq \max\{ \delta,\left \lceil \frac{\Delta +1}{3}\right \rceil \} $ and $\vp (G) \geq \Delta $. In terms of the diameter and radius, we have 
\[ \vp ^-(G) \geq \frac{n(G)-1}{\diam(G)} \quad {\rm and}\quad \vp (G) \geq \frac{n(G)-1}{\rad (G)},\] and for graphs with $\rad (G) \geq 3$, we have $\vp (G) \leq n(G)-\rad(G)-1$. For any vertex $x$, the boundary $\partial (x)$, i.e.\ the set of all vertices maximally distant from $x$, is also an $x$-position set. If $G$ is bipartite, then $\vp (G) \leq \alpha (G)$. The argument in~\cite[Theorem~3.2]{TianKlavzar-2024++} also implies that $\vp^{-}(G) \geq \frac{\gp(G)}{2}$.

In terms of particular graph classes, we have $\vp (C_n) = 2$ for cycles, $\vp (T) = \ell (T)$ for trees and $\vp (G) = \s(G)$ when $G$ is a block graph. For a complete multipartite graph $K_{n_1,\dots ,n_r}$ where $n_1\geq \dots \geq n_r$, we have $\vp (K_{n_1, \dots ,n_r}) = n-n_r$. We have $\vp (K(n,k)) = \binom{n-k}{k}$ for Kneser graphs with sufficiently large $n$. The article also characterises graphs with very large or small values of the vertex position numbers.

An interesting question is how far apart $\vp ^-(G)$ and $\vp (G)$ can be in a graph. A family of graphs with ratio $\frac{\vp (G)}{\vp ^-(G)}$ approaching 6 is constructed, but whether this ratio can be larger is unknown. 

Regarding complexity, it is surprising in view of the fact that most position type problems are NP-hard, that for any vertex $x$ of a graph $G$ the $x$-position number can be determined in polynomial time. In particular, it is proven that for each $x \in V(G)$, the $x$-position number can be computed as an independent set calculated on a graph $G^*_x$, which is derived from $G$ through a transformation (see \cite{ThaChaTuiThoSteErs-2024}).

\begin{theorem} {\rm \cite[Theorem 33]{ThaChaTuiThoSteErs-2024}}
\label{thm:max-x-pos}
Given a graph $G$ and a vertex $x\in V(G)$, a maximum $x$-position set can be computed in $O(n m \log(n^2/m))$ time, where $n = n(G^*_x)$ and $m = m(G^*_x)$.
\end{theorem}

The graph $G^*_x$ has the same vertex set as $G$, so that $n(G) = n(G^*_x)$ and $m(G^*_x) \le \binom{n(G)}{2}$. Hence Theorem~\ref{thm:max-x-pos} yields the following consequence.

\begin{corollary} {\rm \cite[Corollary 34]{ThaChaTuiThoSteErs-2024}}
Given a graph $G$ with order $n$, $\vp^-(G)$ and $\vp(G)$ can be computed in $O(n^4 \log(n))$ time.
\end{corollary}

\subsection{Equidistant numbers}\label{subsec:equidistant}

One special type of general position set that merits study in its own right is the class of equidistant sets. A set $S \subseteq V(G)$ is \emph{equidistant} if there is a value $k \geq 1$ such that $d(u,v) = k$ for all $u,v \in S$. The largest number of vertices in an equidistant set of $G$ is the \emph{equidistant number} $\eq (G)$. Such subsets have already been studied in the Hamming space in the guise of equidistant codes in information theory. It is the discrete version of \emph{equilateral dimension} from geometry.

The equidistant number can be interpreted in terms of the \emph{exact distance $k$-power graph} $G^{[\#k]}$ of $G$, which has vertex set $V(G)$ and two vertices $u,v \in V(G)$ are adjacent in $G^{[\#k]}$ if and only if $d_G(u,v) = k$. A clique in $G^{[\#k]}$ corresponds to an equidistant set in $G$ with distance $k$ and we define $\eq _k(G) = \omega (G^{[\#k]})$. Thus we have $\eq (G) = \max \{ \eq _k(G) : 1 \leq k \leq \diam (G)\} $. Such cliques were studied by Foucaud et al.\ in~\cite{foucaud_cliques_2021}. Equidistant sets are also related to the notion of a \emph{$k$-independent set}, otherwise known as a \emph{$k$-packing}, which is a subset of $V(G)$ with all vertices at distance $> k$. For a given value of $k$ the $k$-independence number $\alpha _k(G)$ is the number of vertices in a largest $k$-independent set (so that $\alpha (G) = \alpha _1(G)$), and hence $\alpha _{k-1} \geq \eq _k(G)$ for $1\leq k \leq \diam (G)-1$.  

Equidistant sets were introduced in an unpublished manuscript with contributions from Erskine, Salia, \v{S}ir\'{a}\v{n}, Taranchuk, Tompkins and Tuite. Any equidistant set is in general position, and hence for any graph $\eq (G) \leq \gp (G)$. Also, any equidistant set is either a clique or an independent set (depending on whether the largest value is obtained for $k = 1$ or $k > 1$) and so $\eq (G) \leq \max \{ \omega (G),\alpha (G)\} $, with equality when $\diam (G) = 2$. The equidistant number of a tree $T$ is given by $\eq (T) = \Delta (T)$. Kneser graphs $K(n,2)$ for $n \geq 7$ and line graphs of complete graphs $L(K_n)$ when $3 \nmid n$ are examples of graphs for which the equidistant number equals the general position number. One interesting result due to \v{S}ir\'{a}\v{n} is that the equidistant number of cubic graphs can be arbitrarily large.

The major works to study this problem are~\cite{Abiad2024} and~\cite{Reijnders}. These papers show that the equidistant number cannot be approximated to within a constant factor in polynomial time unless P = NP. They also study the growth with $n$ of the function $AE(n,k) = \max \{ \alpha _{k-1}(G)-\eq _k(G): G \text{ has order } n\}$. It turns out that 
\[ \lim _{n \rightarrow \infty }\frac{AE(n,k)}{n} = \lim _{n \rightarrow \infty }\max \left\{ \frac{\alpha _{k-1}(G)}{n}: G \text{ has order } n\right\}.\]
This limit is known to be one when $k = 1$, but remains open for larger $k$.

The authors then give a series of bounds on the equidistant number from spectral theory, including the two following results. 

\begin{theorem} {\rm \cite[Proposition 13, Inertial-type bound]{Abiad2024}}
Let $G$ be a graph with adjacency eigenvalues $\lambda _1 \geq \dots \geq \lambda_n$ and adjacency matrix $A$. Let $p \in \mathbb{R}_t[x]$ with corresponding parameters $W(p) := \max _{u \in V(G)}\{ (p(A))_{uu}\}$ and $w(p) := \min _{u \in V(G)}\{ (p(A))_{uu}\}$. Then the $(t+1)$-equidistant number of $G$ satisfies the bound \[ \eq _{t+1}(G) \leq \min \{ |\{ i:p(\lambda _i) \geq w(p)\} |,|\{ i:p(\lambda _i) \leq W(p)\}| \}. \]
\end{theorem}

\begin{theorem} {\rm \cite[Proposition 14, Ratio-type bound]{Abiad2024}}
Let $t \geq 1$ and $G$ be a regular graph with adjacency eigenvalues $\lambda _1 \geq \dots \geq \lambda _n$ and adjacency matrix $A$. Let $p \in \mathbb{R}_t[x]$ with corresponding parameters $W(p) := \max _{u \in V(G)}\{ (p(A))_{uu}\}$ and $\lambda (p) := \min _{i \in [2,n]}\{ p(\lambda _i)\} $ and assume $p(\lambda _1) > \lambda (p)$. Then \[ \eq _{t+1}(G) \leq n(G)\frac{W(p)-\lambda(p)}{p(\lambda _1)-\lambda (p)}.\]
\end{theorem}
The paper closes with a computational study to compare the bounds. It is an interesting question whether spectral theory may have further applications in the general position problem.

One question for future research is to find the largest size of a graph with given order and equidistant number. For small values of $n$ and the equidistant number, the extremal graphs bear an intriguing resemblance to Ramsey graphs. 

\subsection{Variety of general position sets}~\label{subsec:variety}

If $X \subseteq V (G)$, then the requirement that all pairs of vertices from $X$ are $X$-positionable leads to the concept of general position sets. Alternatively, we can also require that other pairs of vertices of $G$ are $X$-positionable. Considering all natural cases leads us to the variety of general position sets as introduced in~\cite{TianKlavzar-2025}, which contains the following set types in addition to general position sets. Setting $\overline{X} = V (G)\setminus X$, the set $X$ is 
\begin{itemize}
\item a {\em total general position set}, if every pair $u, v \in V(G)$ is $X$-positionable;
\item an {\em outer general position set}, if every pair $u, v \in X$ is $X$-positionable, and every pair $u \in X$, $v \in \overline{X}$ is $X$-positionable; and
\item  a {\em dual general position set}, if every pair $u, v \in X$ is $X$-positionable, and every pair $u, v \in \overline{X}$ is $X$-positionable.
\end{itemize} 
The cardinality of a largest total general position set, a largest outer general position set, and a largest dual general position set are denoted by $\gpt(G)$, $\gpo(G)$, and $\gpd(G)$ respectively.

Total general position sets are fully understood: 

\begin{theorem}  {\rm \cite[Theorem~2.1]{TianKlavzar-2025}}
\label{thm:total-characterisation}
Let $G$ be a graph and $X\subseteq V(G)$. Then $X$ is a total general position set of $G$ if and only if $X\subseteq S(G)$. Moreover, $\gpt(G) = \s(G)$.
\end{theorem}

If $u$ and $v$ are vertices of a graph $G$ such that $d_G(u,v) = \diam(G)$, then $\{u,v\}$ is an outer general position set. Hence, if $G$ is a connected graph of order at least two, then $\gpo(G)\ge 2$. Recalling from Subsection~\ref{sec:strong-resolving} the definition of the strong resolving graph $G\sr$ of $G$, outer general position sets can be characterised as follows.

\begin{theorem} {\rm \cite[Theorem~2.3]{TianKlavzar-2025}}
\label{thm:characterise outer}
Let $G$ be a connected graph and $X\subseteq V(G)$. Then $X$ is an outer general position set of $G$ if and only if each pair of vertices from $X$ is mutually maximally distant. Furthermore,
$$\gpo(G) = \omega(G\sr).$$
\end{theorem}

As for the dual general position sets, they are exactly those general position sets which have convex complement: 

\begin{theorem} {\rm \cite[Theorem~3.1]{TianKlavzar-2025}}
\label{thm:dual-characterization}
Let $X$ be a general position set of a graph $G$. Then $X$ is a dual general position set if and only if $G-X$ is convex.
\end{theorem}

The paper~\cite[Theorem~3.1]{TianKlavzar-2025} pays particular attention to graphs $G$ with $\gpd(G) \in \{0,1\}$, and to the variety of general position sets in  Cartesian products. From the obtained results one can deduce, among other things, that
if $n\ge 3$, then 
\begin{align*}
\gp(K_n\cp K_{2n}) & = 3n - 2\,,\\
\gpd(K_n\cp K_{2n}) & = 2n\,, \\
\gpo(K_n\cp K_{2n}) & =  n\,, \\
\gpt(K_n\cp K_{2n}) & = 0\,.
\end{align*}
This demonstrates that the four general position invariants can be arbitrarily far apart. 

In~\cite{Dokyeesun-2026}, outer, dual, and total general position sets are investigated on strong products and on lexicographic products. For the strong product, sharp lower and upper bounds are proved for the outer and the dual version. For the lexicographic product, the outer general position number is determined in all cases, and the dual general position number in many cases.  Along the way some results on outer general position sets are also derived. 

As already mentioned, $\gp(S_p^2)$ and $\gp(S_3^n)$ are determined in\cite{Roy-2026}, where $S_p^n$  are the Sierpi\'nski graphs. In the same paper it is proved that if $p\geq 3$ and $\tau\in \{\gpd, \gpt, \gpo\}$, then 
$\tau(S_p^2) = p$ and $\# \tau(S_p^2) = 1$~\cite[Corollary 3.4, Theorem 3.5]{Roy-2026}. Furthermore, if $n\geq 1$, then $\gpt(S_3^n) = \gpd(S_3^n) = 3$, and $\#\gpt(S_3^n) = 1$ and $\#\gpd(S_3^n) = 1$~\cite[Corollary 4.3]{Roy-2026}. It remains open to determine the exact values of $\gpo(S_3^n)$ for $n\geq 3$.

The variety of general position sets was investigated also on glued binary trees, which are defined as follows. A {\em perfect binary tree} $T_{r}$ of depth $r\ge 1$ is a rooted binary tree in which all non-leaf vertices have two children, and all leaves have depth $r$. A {\em glued binary tree} $GT(r)$, $r\ge 1$, is obtained from two copies of $T_{r}$ by pairwise identifying their leaves. For these graphs we have: 

\begin{theorem} {\rm \cite[Theorem~2.7]{KlaLakRoy-2025}}
\label{thm:gp-gpo-glued binary tree}
If $r\geq 2$, then $\gpo (GT(r))  = \gp (GT(r)) = 2^r.$ Also, $\#\gpo(GT(r)) = 1$ and $\#\gp(GT(r)) = 2^{r-1} + 1$.
\end{theorem}

\begin{theorem} {\rm \cite[Theorem~2.10]{KlaLakRoy-2025}}
If $r\ge 2$, then $\gpt(GT(r)) = \gpd(GT(r)) = 0$.    
\end{theorem}

In Subsection~\ref{sec:vertex-edge-delete} we have described the general position number under vertex and edge removal. The paper~\cite{Tian-2026} does this for the entire variety of the general position problem. For vertex removal it is proved that (i) if $x$ is not a cut-vertex, then $\gpt(G) - 1 \le \gpt(G - x) \le \gpt(G) + \deg_G(x)$, (ii) the outer general position number of a graph when a vertex is removed cannot be bounded in terms of the outer general position number of the original graph, but if $x$ lies in some $\gpo$-set, then $\gpo(G) - 1 \le \gpo(G - x)$, and (iii) $\gpd(G - x)$ can also be arbitrarily larger/smaller than $\gpd(G)$, but if $x$ is not a cut-vertex and lies in some $\gpd$-set of $G$, then $\gpd(G) - 1 \le \gpd(G - x)$. For the edge removal it is proved that (i) $\gpt(G) -|S(G)_{e}| \le \gpt(G-e) \le \gpt(G) +2\,$, where $S(G)_{e}$ denotes the set of simplicial vertices in $G$ which are adjacent to both $u$ and $v$, (ii) $\gpo(G)/2 \le \gpo(G-e)\leq\ 2\gpo(G)$, and (iii) the difference $\gpd(G) - \gpd(G-e)$ can be arbitrarily large. All the bounds described in this paragraph are also demonstrated to be sharp in~\cite{Tian-2026}. 

\subsection{Edge general position}~\label{subsec:edgegp}

The edge version of the general position problem was introduced in~\cite{Manuel-2022}; a set $S$ of edges of a graph $G$ is an {\em edge general position set} if no shortest path of $G$ contains three edges of $S$.  An edge general position set of maximum cardinality is a {\em $\gpe$-set} of $G$, its cardinality is the \textit{edge general position number} (for short {\em $\gpe$-number}) of $G$ and is denoted by $\gpe(G)$.

The edge general position problem is intrinsically different from the standard general position problem; one outstanding example of this is the fact that one can determine the edge general position number of hypercubes exactly.  

\begin{theorem} {\rm \cite[Theorem 3.2]{Manuel-2022}}
\label{thm:edges-hypecubes}
If $r\ge 2$, then $\gpe(Q_r) = 2^r$.
\end{theorem}

The above theorem was derived in a broader context by considering edge general position sets in partial cubes (graphs which admit isometric embeddings into hypercubes). The lower bound $\gpe(Q_r) \ge 2^r$ follows from the fact~\cite[Lemma 3.1]{Manuel-2022} that in a partial cube, the union of two arbitrary so-called $\Theta$-classes forms an edge general position set. This approach was later used in~\cite{KlavzarTan-2023}, where edge general position sets were considered in Fibonacci cubes and Lucas cubes (for more information on these cubes see the book~\cite{egecyoglu-2023}), which form two important classes of partial cubes. In particular, the union of two largest $\Theta$-classes of a Fibonacci cube or a Lucas cube forms a maximal edge general position set~\cite[Theorem 4.3]{KlavzarTan-2023}. For grids we have: 

\begin{theorem} {\rm \cite[Theorems 4.1, 4.2]{Manuel-2022}}
\label{thm:edges-grids}
If $r\ge s\ge 2$, then 
$$\gpe(P_r \cp P_s) = 
\begin{cases}
r + \normalfont{2}; & s = 2\,,\\
2r; & s = 3\,,\\
2r + 2s - 8; & s \geq 4\,.
\end{cases}
$$
Moreover, if $r,s\ge 5$, then the $\gpe$-set of $P_r \cp P_s$ is unique.
\end{theorem}

The edge general position number of several additional graph products was investigated in~\cite{Hamed-2026b}. For hierarchical products with a complete first factor, a lower bound is proved in~\cite[Theorem 2.1]{Hamed-2026b}, while for the case when the second factor is a tree a formula is deduced in~\cite[Theorem 2.2]{Hamed-2026b}. Similarly, for corona products with a complete first factor a formula for the edge general position number is proved in~\cite[Theorem 3.1]{Hamed-2026b}, while in~\cite[Theorem 4.2]{Hamed-2026b} a related formula is deduced for arbitrary edge coronas. 

The paper~\cite{Hamed-2026b} also presents an integer linear programming model to address the edge
general position problem and applies it to certain specific midi fullerene graphs. 

It is clear that $\gpe(G) = m(G)$ if and only if $\diam(G) \le 2$. The graphs $G$ with $\gpe(G) = m(G) - 1$ can be classified into two families. The first family ${\cal G}_1$ consists of the graphs $G$ that can be obtained from an arbitrary graph $H$ with $\diam(H) = 2$ by attaching a leaf to a vertex $u$ of $H$ with $\deg_H(u) \le n(H) - 2$. The second family ${\cal G}_2$ consists of the graphs $G$ constructed as follows. Let $G_0$, $G_1$, and $G_2$ be arbitrary not necessarily connected graphs. Then $G$ is obtained from the disjoint union of $G_0$, $G_1$, $G_2$, and $K_2$, where $V(K_2) = \{x_1, x_2\}$, by adding all possible edges between $K_2$ and $G_0$, adding the edges between $x_1$ and all vertices from $G_1$, and adding the edges between $x_2$ and all vertices of $G_2$. 

\begin{theorem} {\rm \cite[Theorems 2.2]{Tian-2024b}}
\label{thm:edges-m(G)-1}
If $G$ is a connected graph with $n(G) \ge 4$, then $\gpe(G) = m(G) - 1$ if and only if $G\in {\cal G}_1 \cup {\cal G}_2$.
\end{theorem}

On the other hand, it is straightforward to see that $\gpe(G) = 2$ if and only if $G$ is a path. It is not much more difficult to prove that $\gpe(G) = 3$ if and only if $G$ is $K_3$ or a tree with three leaves~\cite[Proposition 2.4]{Tian-2024b}. The article~\cite{Tian-2024b} also reports that if $G$ is a graph with $\gpe(G) = 4$, then $\Delta(G) \le 4$ and, moreover, in the case when $\Delta(G) = 4$, the graph $G$ must be bipartite. This line of research was later taken up independently by two groups, Cao, Ji, and Wang~\cite{cao-2025}, and Li and Gong~\cite{li-2025}. In both papers a complete characterisation of graphs $G$ with $\gpe(G) = 4$ is presented. Here we follow the description from the latter paper. It is first proved that if $G$ is a graph with $\gpe(G) = 4$, then $G$ is a cactus graph. The list of all graphs $G$ with $\gpe(G) = 4$ then consists of the following families. 
\begin{itemize}
\item ~\cite[Theorem 3.1]{li-2025} If $T$ is a tree, then $\gpe(T) = 4$ if and only if $T$ contains exactly four leaves.
\item ~\cite[Theorem 3.4]{li-2025} If $U$ is a unicyclic graph, then $\gpe(U) = 4$ if and only if $U$ belongs to the family ${\cal U}$, which contains some of the unicyclic graphs with at most four vertices of degree at least three, each having only one or two leaves attached. For the exact definition of the family ${\cal U}$, see~\cite[p.~8]{li-2025}.
\item ~\cite[Theorem 3.5]{li-2025} If $B$ is a cactus graph with at least two cycles, then $\gpe(B) = 4$ if and only if $B$ belongs to the family ${\cal B}$, where this family contains even cycle chains in which cycles are attached at antipodal vertices, and some simple vertex- or edge-deleted subgraphs of even cycle chains that contain at least three cycles. For the exact definition of the family ${\cal B}$, see~\cite[p.~11]{li-2025}.
\end{itemize}
In addition to the characterisation of graphs $G$ with $\gpe(G) = 4$, the paper~\cite{cao-2025} of Cao, Ji and Wan presents a classification of graphs $G$ with $\diam(G) = 3$, $\gpe(G) = m(G)-2$, and having a $\gpe$-set $S$ such that if $e$ and $f$ are edges not in $S$, then $G - \{e, f\}$ is disconnected. Edge general position numbers of cactus graphs are investigated further in~\cite{Cao-2026}, which gives upper and lower bounds in terms of the number of cycles and leaves and characterises the equality cases. The paper~\cite{Cao-2026} gives an independent proof of the characterisation of cactus graphs with edge general position number two.

Recall that $\s(G)$ denotes the number of simplicial vertices of $G$. We further say that an edge of $G$ is \textit{simplicial} if it is incident with at least one simplicial vertex. Setting $\s'(G)$ to be the number of simplicial edges of $G$, we have: 

\begin{theorem} {\rm \cite[Theorem 3.2]{Tian-2024b}}
\label{thm:edges-block graphs}
If $G$ is a block graph, then $\s'(G)\leq \gpe(G)\leq \binom{\s(G)}{2}+1$. Both bounds are sharp.
\end{theorem}

In~\cite{Tian-2024b}, exact values of $\gpe$ are determined for several specific families of block graphs. We add that the edge $k$-general position problem was introduced and studied in~\cite{ManuelPrabhaKlavzar-2023+}. The problem is to find a largest set $S$ of edges of $G$ such that at most $k-1$ edges of $S$ lie on a common geodesic, and the paper explores this in torus graphs, hypercubes and Bene\v{s} networks.

\subsection{Steiner general position}\label{subsec:steiner}

The paper~\cite{KlaKuzPetYer-2021} introduces $k$-Steiner general position sets and $k$-Steiner general position numbers as follows. 

Let $G$ be a connected graph and $X\subseteq V(G)$. The {\em Steiner distance} $d_G(X)$ of $X$ is the minimum size of a connected subgraph of $G$ containing $X$~\cite{Chartrand-1989}. Such a subgraph is clearly a tree and is called a \emph{Steiner $X$-tree}. For a positive integer $k$, we say that $A\subseteq V(G)$ is a {\em $k$-Steiner general position set} if for every $B\subseteq A$, $|B| = k$, and for every Steiner $B$-tree $T_B$, it follows that $V(T_B)\cap A = B$. That is, $A$ is a $k$-Steiner general position set if no $k+1$ distinct vertices from $A$ lie on a common Steiner $B$-tree, where $B\subseteq A$ and $|B| = k$. Clearly, if $|A|\le k$, then $A$ is a $k$-Steiner general position set. The {\em $k$-Steiner general position number} $\sgp_k(G)$ of $G$ is the cardinality of a largest $k$-Steiner general position set in $G$. Note that $\sgp_2(G) = \gp(G)$.

In~\cite[Proposition 2.2]{KlaKuzPetYer-2021} it was observed that if $G$ is a graph and $k\in \{2, \ldots, n(G)-1\}$, then $\sgp_k(G)=n(G)$ if and only if $G$ is $(n(G)-k+1)$-connected. If $T$ is a tree, then $\sgp_k(T) = \ell(T)$ when $k\leq \ell(T)$, and $\sgp_k(T) = k$ when $k > \ell(T)$. The corresponding result for cycles is as follows. 

\begin{theorem}  {\rm \cite[Theorem 3.2]{KlaKuzPetYer-2021}}
If $n\ge 3$ and $k\in \{2,\ldots, n-1\}$, then
$$\sgp_k(C_n) = 
\left\{\begin{array}{ll}
k; & k\in \{\left\lfloor\frac{2n}{3}\right\rfloor, \ldots, n-2\}\,, \\
k+1; & \text{otherwise}\,.
\end{array}
\right.$$ 
\end{theorem}

In~\cite[Theorem 4.2]{KlaKuzPetYer-2021}, a formula for $\sgp_k$ of the join $G\vee H$ of graphs $G$ and $H$ is given in terms of two related invariants which we do not explain here. Instead, we just mention that the formula can be simplified in many cases. We extract the next two formulas from a list of five formulas given in~\cite[Corollary 4.4]{KlaKuzPetYer-2021}: 
\begin{itemize}
\item If $n\geq 6$ and $k\in \{2, \ldots, n-1\}$, then $$\sgp_k(W_n) = \sgp_k(K_1\vee C_{n-1}) = \max\left\{k+1,n-2-\left\lfloor \frac{n-2}{k+1}\right\rfloor\right\}\,.$$
\item If $r\leq s$ and $k\in \{2, \ldots, r+s-1\}$, then
$$\sgp_k(K_{r,s})=\sgp_k(\overline{K}_r\vee \overline{K}_s)=\left\{\begin{array}{ll}
\max\{s,\min\{k-1,r\}+k-1\}; & k\leq s\,, \\[0.1cm]
r+s; & k>s\,.    
\end{array}
\right.$$
\end{itemize}
The paper~\cite{KlaKuzPetYer-2021} further introduces bounds and exact results for the lexicographic product of graphs, for split graphs, and proves in~\cite[Theorem 6.2]{KlaKuzPetYer-2021} that $\sgp_k(P_{\infty}\cp P_{\infty})\ge 2k$.

Many open problems on the Steiner general position problem remain to be investigated. For instance, does the equality $\sgp_k(P_{\infty}\cp  P_{\infty})=2k$ hold for $k>2$?

\subsection{Detour position problem}\label{subsec:detour}
Motivated by the detour hull number (which we do not define here), rather than position problems, the paper~\cite{DelbinPrema-2018} considers the problem resulting from replacing `shortest path' in the general position problem by `longest path'. A longest path in a graph $G$ is called a \emph{detour}. For any two vertices $u,v \in V(G)$ the \emph{detour distance} $D(u,v)$ is the length of a longest path from $u$ to $v$; it is known that this distance function is in fact a metric on the graph, see~\cite{Chartrand-2004}. 

The authors of~\cite{DelbinPrema-2018} define a \emph{detour irredundant set} to be a subset $S \subseteq V(G)$ such that for any $u,v \in S$ no $u,v$-detour contains a third vertex of $S$. To fit the terminology of the rest of the literature, we will refer to these as \emph{detour position sets} and the largest number of vertices in a detour position set of $G$ as the \emph{detour position number} of $G$, which we will write as $\dpn (G)$.

The authors claim in~\cite[Theorem 3]{DelbinPrema-2018} that any longest path in a graph can contain at most two vertices of a detour position set, which would imply the upper bound $n(G)-D+1$, where $D$ is the length of a longest path in $G$ (known as the \emph{detour diameter}). However, the proof assumes that a subpath of a detour is also a detour, which is not true in general, so at present we regard the $n(G)-D+1$ bound as a conjecture. 

For any tree $T$, the longest paths and shortest paths coincide, so naturally $\dpn (T) = \gp (T) = \ell (T)$. The paper discusses the relation between the detour position number and the detour hull number and also shows that for any $k \geq 2$ there is a graph with $\dpn (G) = k$ and any feasible combination of detour radius and detour diameter. It is shown that a graph has detour position number $n(G)-1$ if and only if it is $K_3$ or a star, and detour position number $n(G)-2$ if and only if it is a star plus an edge or a double star (although both of these results depend on the $n(G)-D+1$ bound).

\section{Variations Not Based on Distance}\label{sec:variantsnotdistance}

In this section we examine variants of the general position problem that result from changing the path family to paths not related to distance. The main problem in this section is the \emph{monophonic position problem} (see Subsection~\ref{subsec:mp}), for which the path family is all induced paths. We also briefly discuss in Subsection~\ref{subsec:all paths} the effect of making the path family as large as possible.

\subsection{Monophonic position sets}\label{subsec:mp}

Graph convexity for monophonic (i.e.\ induced) paths has a very wide literature. This inspired the first investigation of the \emph{monophonic position problem} in~\cite{Thomas-2024b}: what is the largest possible number of vertices in a subset $S \subseteq V(G)$ such that no monophonic path of $G$ passes through three vertices of $S$? Such a set in a graph $G$ is a \emph{monophonic position set} of $G$ and the number of vertices in a largest such set is the \emph{monophonic position number} $\mono (G)$ of $G$. The number $\mono(G)$ is also referred to as the \emph{$\mono$-number} of $G$. This problem was later considered independently in~\cite{Araujo-2025}. A largest monophonic position set of the Petersen graph can be seen in the centre of Figure~\ref{fig:Petersen}.

For any distance-hereditary graph all induced paths are shortest paths, so that the monophonic and general position numbers coincide. It follows immediately that the monophonic position number is equal to the general position number for block graphs (in particular trees), complete multipartite graphs, cographs, etc. However, equality can also hold for non-distance-hereditary graphs. To illustrate this, consider the {\em house graph} $G$, which is obtained from the cycle $C_4$ by adding a new vertex and joining it to two adjacent vertices of $C_4$. Then it is easy to verify that $\gp(G)=\mono(G)=3$.  Straightforward upper bounds on $\mono (G)$ include $n(G)-\diam_{\rm m}(G)+1$ (where $\diam_{\rm m}(G)$ is the \emph{monophonic diameter}, that is, the length of a longest induced path in $G$), two times the induced path number (the smallest number of induced paths needed to cover $V(G)$) and $n(G)-\cut(G)$ (where $\cut(G)$ is the number of cut-vertices in $G$, since there always exists a largest monophonic position set that contains no cut-vertices). A lower bound for $\mono (G)$ is $\s(G)$. 

Trivially any monophonic position set is also in general position, so $\mono (G) \leq \gp (G)$ for any graph $G$. However, it is shown in~\cite{Thomas-2024b} that for any $2 \leq a \leq b$ there is a graph $G$ with $\mono (G) = a$ and $\gp (G) = b$ (the graph used is a gear graph, also known as a bipartite wheel graph). The article~\cite{TuiThoCha-2025} (and the extended abstract~\cite{TuiThoCha-2022}) seeks to optimise this result by finding the smallest such graph for given $a$ and $b$. If we denote the order of this smallest graph by $M(a,b)$, then we know that:

\begin{itemize}
			\item $M(2,3) = 5$ and for $b \geq 4$, we have $M (2,b) \leq \lceil \frac{3b}{2} \rceil +1  $, with equality for $4 \leq b \leq 8$.
			\item For $3 \leq a < b$ and $\frac{b}{2} \leq a$, we have $M(a,b) = b+2$ (proven extremal).
			\item For $3 \leq a < \frac{b}{2}$, we have $M(a,b) \leq b-a+2 + \lceil \frac{b}{2} \rceil $. 
		\end{itemize}
It is conjectured that the first and third estimates in the list above are exact. The extremal graphs for the second item are not unique. In a similar vein, estimates are given for the smallest size of a graph with given order, $\gp $-number and $\mono $-number (the answer is $n(G)+O(1)$, but is not known exactly). 

Whilst the structure of monophonic position sets is not characterised as in the case of general position sets, we have the following lemma.

\begin{lemma} {\rm \cite[Lemma 3.1]{Thomas-2024b}}\label{lem:mpstructure} 
	If $G$ is a connected graph and $M\subseteq V(G)$ a monophonic position set of $G$, then $G[M]$ is a disjoint union of $k \geq 1$ cliques $G[M]=\bigcup_{i=1}^{k} W_i$. If $k \geq 2$, then for every $i \in [k]$ any two vertices of $W_i$ have a common neighbour in $G -M$.
\end{lemma}

As a consequence of Lemma~\ref{lem:mpstructure}, the $\mono $-number of a triangle-free graph is bounded above by the independence number. It is shown in~\cite{Thomas-2024b} that the monophonic position number of a corona product is $\mono(G \odot H) = n(G)\mono(H)$ and the monophonic position number of the join of two graphs is given by $\mono(G \vee H) = \max\{\omega (G) + \omega (H), \mono(G), \mono(H) \} $. The article determines the $\mono $-numbers of unicyclic graphs, complements of bipartite graphs and split graphs, and characterises the split graphs that achieve equality $\mono (G) = \max \{ \omega (G),\alpha (G)\} $ using Hall's Theorem. Computational results also suggested the following two conjectures.

\begin{conjecture} {\rm \cite[Conjecture 2.9]{Thomas-2024b}}
The largest possible monophonic position number of a cubic graph with order $n$ is $\frac{n}{3}+O(n)$.
\end{conjecture}

\begin{conjecture} {\rm \cite[Conjecture 3.4]{Thomas-2024b}}
For sufficiently large $g$, the monophonic position number of a $(d,g)$-cage $G$ satisfies $\mono (G) < d$.
\end{conjecture}

The main topic of the paper~\cite{TuiThoCha-2025} (and extended abstract~\cite{TuiThoCha-2022}) is the extremal problem of the largest possible size of a graph with given position numbers. The function $\mex (n;a)$ is defined to be the largest size of a graph with order $n$ and monophonic position number $a$, and $\gex (n;a)$ is the largest size of a graph with order $n$ and general position number $a$. The asymptotic behaviour of these functions is completely different.

\begin{theorem} {\rm \cite[Theorems 3.5 and 3.7]{TuiThoCha-2025}}\label{thm:extremal}
If $n \geq a$, then 
\begin{enumerate}
    \item[(i)] $\mex (n;a) = (1-\frac{1}{a})\frac{n^2}{2}+O(n)$,
    \item[(ii)] $\gex(n;a) \leq \frac{R(a,a+1)-1}{2}n$,
\end{enumerate}
where $R(a,a+1)$ is a Ramsey number.
\end{theorem}
The reason the first result holds is that a graph $G$ with $\mono (G) = a$ cannot contain a $K_{a+1}$ and so has size restrained by Tur\'{a}n's Theorem, but one can remove the edges of a linear number of copies of $K_a$ from the Tur\'{a}n graph to obtain a graph with monophonic position number $a$. Thus in this sense the monophonic position number is tied to the clique number. For $\mono $-number two the answer is known exactly: $\mex (n;2) = \left \lceil \frac{(n-1)^2}{4} \right \rceil$ for $n \geq 6$. There is a unique extremal graph for each order, namely $K_{\left \lfloor \frac{n}{2} \right \rfloor ,\left \lceil \frac{n}{2} \right \rceil }$ with a matching of size $\left \lfloor \frac{n}{2} \right \rfloor $ deleted.

By contrast, the best upper bounds for $\gex (n;a)$ come from Ramsey's Theorem; if a vertex $u$ has large degree, then by Ramsey's Theorem $G[N_G(u)]$ will contain either a large clique or a large independent set, either of which will constitute a large general position set (and in the first case we can include the vertex $u$). Of course any independent union of cliques in $G[N_G(u)]$ will yield a general position set, so one can improve the multiplicative constant slightly. Nevertheless, the upper bound in Theorem~\ref{thm:extremal} seems much too large, and improving this estimate is an interesting problem. 

\begin{figure}[ht!]
		\centering
		\begin{tikzpicture}[x=0.4mm,y=-0.4mm,inner sep=0.2mm,scale=0.4,very thick,vertex/.style={circle,draw,minimum size=10}]
			\node at (70,0) [vertex] (x1) {$$};
			\node at (100,-100) [vertex] (x2) {$$};
			\node at (-100,-100) [vertex] (x3) {$$};
			\node at (-70,0) [vertex] (x4) {$$};
			\node at (-100,100) [vertex] (x5) {$$};
			\node at (100,100) [vertex] (x6) {$$};
			
			\node at (0,-140) [vertex] (y) {$$};
			\node at (0,0) [vertex] (z) {$$};
			
			\node at (0,-170) [vertex] (y') {$$};
			\node at (0,70) [vertex] (z') {$$};
			
			\path
			(x1) edge (x2)
			(x2) edge (x3)
			(x3) edge (x4)
			(x4) edge (x5)
			(x5) edge (x6)
			(x6) edge (x1)
			(x1) edge (x3)
			(x2) edge (x4)
			(x3) edge (x5)
			(x4) edge (x6)
			(x5) edge (x1)
			(x6) edge (x2)

			(z) edge (x2)
			(z) edge (x3)
			(z) edge (x5)
			(z) edge (x6)
			(x1) edge[bend right] (y)
			(x4) edge[bend left] (y)
			(y) edge (y')
			(z) edge (z')
			
			;
		\end{tikzpicture}
		\caption{A graph with order ten, gp-number three and largest size}
		\label{gex graph}
	\end{figure} 

The final section of~\cite{TuiThoCha-2025} classifies the possible diameters of a graph with $\mono $-number $a$. If $a \geq 3$, then all diameters in the range $2 \leq D \leq n-a+1$ are possible and if $a = 2$, then the possible diameters are $D = n-1$, the integers $D$ in the range $3 \leq D \leq \lfloor \frac{n}{2} \rfloor $ and, if $n \in \{ 3,4,5,8\} $ or $n \geq 11$, diameter $D = 2$. The case of graphs with $\mono $-number two and diameter two is rather interesting, and a new graph operation is defined that, given such a graph with order $n$, will produce a new graph with the same properties and order either $3n+2$, $3n+1$ or $3n$. An example of this operation applied to $C_4$ can be seen in Figure~\ref{fig:3r+2 construction}. This construction may be worthy of further study.

\begin{figure}[ht!]
		\centering
		\begin{tikzpicture}[x=0.4mm,y=-0.4mm,inner sep=0.2mm,scale=0.6,very thick,vertex/.style={circle,draw,minimum size=10}]
			\node at (-150,-50) [vertex] (x1) {$$};
			\node at (-50,-50) [vertex] (x2) {$$};
			\node at (50,-50) [vertex] (x3) {$$};
			\node at (150,-50) [vertex] (x4) {$$};

			\node at (-90,0) [vertex] (h1) {$$};
			\node at (-30,0) [vertex] (h2) {$$};
			\node at (30,0) [vertex] (h3) {$$};
			\node at (90,0) [vertex] (h4) {$$};

			\node at (-150,50) [vertex] (y1) {$$};
			\node at (-50,50) [vertex] (y2) {$$};
			\node at (50,50) [vertex] (y3) {$$};
			\node at (150,50) [vertex] (y4) {$$};
			
			\node at (0,-100) [vertex] (z1) {$$};
			\node at (0,100) [vertex] (z2) {$$};
			
			\path
			(x1) edge (y2)
			(x1) edge (y3)
			(x1) edge (y4)
			
			(x2) edge (y1)
			(x2) edge (y3)
			(x2) edge (y4)
			
			(x3) edge (y1)
			(x3) edge (y2)
			(x3) edge (y4)
			
			(x4) edge (y1)
			(x4) edge (y2)
			(x4) edge (y3)
			
			(x1) edge (h1)
			(y1) edge (h1)
			(x2) edge (h2)
			(y2) edge (h2)
			(x3) edge (h3)
			(y3) edge (h3)
			(x4) edge (h4)
			(y4) edge (h4)
			
			(h1) edge (h2)
			(h2) edge (h3)
			(h3) edge (h4)
			(h4) edge [bend left] (h1)
			
			(z1) edge (x1)
			(z1) edge (x2)
			(z1) edge (x3)
			(z1) edge (x4)
			
			(z2) edge (y1)
			(z2) edge (y2)
			(z2) edge (y3)
			(z2) edge (y4)
			
			(z1) edge (z2)
			
			;
		\end{tikzpicture}
		\caption{An example of the operation for making new graphs with diameter two and $\mono $-number two.}
		\label{fig:3r+2 construction}
	\end{figure} 

Computation on this problem suggested that it is always possible to find circulant graphs (i.e.\ a Cayley graph on a cyclic group) with this property.

\begin{conjecture} {\rm\cite[Conjecture 4.5]{TuiThoCha-2025}}
		For any $n \geq 11$, there is a circulant graph with order $n$, monophonic position number two and diameter two.
	\end{conjecture}

   The monophonic position numbers of the Cartesian and lexicographic product graphs are discussed in~\cite{ChandranKlavzar2025}. For the Cartesian product, as noted in~\cite[Observation 3.1, Corollary 3.6]{ChandranKlavzar2025}, we have that for any graphs \( G \) and \( H \), 
$$
\max \{\omega(G), \omega(H)\} \leq \mono(G \cp H) \leq \max \{\mono(G), \mono(H)\}.
$$ 
As these bounds coincide for products of paths and cycles, we have 
$\mono(P_m \cp P_n) = \mono(P_n \cp C_r) = \mono(C_r \cp C_s) = 2$, see~\cite[Corollary 3.7]{ChandranKlavzar2025}. As reported in~\cite[Corollary 3.8]{ChandranKlavzar2025}, if $H$ is a connected graph and $n \geq \mono(H)$, then $\mono(K_n \cp H) = n.$ 
This in particular yields $\mono (K_n \cp P_m) = n$ for $n \geq 2$, $\mono (K_n \cp C_m) = n$ for $n \geq 3$, and $\mono (K_n \cp K_m) = \max\{ n,m\} $. 
 
A vertex subset $S$ of a graph that is both an independent set and in monophonic position is called an \emph{independent monophonic position set}. The largest order of an independent monophonic position set is the \emph{independent monophonic position number} of $G$, denoted by $\mono_{\rm i}(G)$. With this notation in hand, the paper\cite{ChandranKlavzar2025} then reports several additional structural properties of monophonic position sets of Cartesian products, which lead to the following improved tight bounds. To this end, we fix $\sigma (G) = 1$ if $G$ contains a simplicial vertex and $\sigma (G) = 0$ otherwise.

\begin{theorem} {\rm\cite[Theorem 3.16]{ChandranKlavzar2025}}
\label{thm:general bounds}
If $G$ and $H$ are connected graphs, then $$\max\{\omega(G), \omega(H)\}\leq \mono (G \cp H) \leq \max \{ \omega (G), \omega (H), \sigma (G) \mono_i(H), \sigma (H) \mono_i(G)\}.$$ 
\end{theorem}
Note that Theorem~\ref{thm:general bounds} in particular implies that if neither $G$ nor $H$ has simplicial vertices, then $\mono (G \cp H) = \max \{ \omega (G), \omega (H)\}$.      

\cite[Theorem 3.18]{ChandranKlavzar2025} further provides improved bounds for the monophonic position number of Cartesian products of triangle-free graphs. 

\begin{theorem}{\rm \cite[Theorem 3.18]{ChandranKlavzar2025}}
If $G$ and $H$ are connected triangle-free graphs of order at least three, then 
$$\mono(G\cp H)\leq\max\{2, \sigma(G)\Delta(H), \sigma(H)\Delta(G)\}\,.$$ 
Moreover, the bound is tight when both $G$ and $H$ are star graphs.
\end{theorem}

For a monophonic position set $M$ of $G$ we denote the components of $G[M]$ by $A_1,A_2,\ldots,A_k,B_1,\dots ,B_r $, where $|A_i|\geq 2$ for each $i\in[k]$ and $|B_j|=1$ for each $j\in [r]$. Also we fix $n_M = \sum _{i=1}^k |A_{i}|$ and write $r_M = r$ to emphasise that $r$ is a function of $M$. Then for any monophonic position set $M$ of $G$ we have $|M| = n_M+r_M$. Using this notation,~\cite[Theorem 4.4]{ChandranKlavzar2025} provides a formula for the monophonic position number of the lexicographic product of arbitrary graphs. The result reads as follows.

\begin{theorem}{\rm \cite[Theorem 4.4]{ChandranKlavzar2025}}
\label{thm:6.5}
Let $G$ be a connected graph of order at least two and let $\mathcal{M}$ be the collection of all monophonic position sets of $G$. Then $$ \mono(G\circ H)=\max_{M\in \mathcal{M}}\{n_M\cdot\omega(H)+ r_M\cdot \mono(H)\}\,.$$
\end{theorem}

To discuss the complexity of the monophonic position problem, we shall define the decision version of the problem.
\begin{definition}\label{def:probmp}
        {\sc Monophonic position set} \\
        {\sc Instance}: A graph $G$, a positive integer $k\leq n(G)$. \\
        {\sc Question}: Is there a monophonic position set $S$ for $G$ such that $|S|\geq k$?
\end{definition}

\begin{theorem}\label{thm:NPM}{\rm\cite[Theorem 6.2]{Thomas-2024b}}
       The {\sc Monophonic position set} problem is NP-hard.  
\end{theorem}

To prove Theorem~\ref{thm:NPM} a reduction from the well-known NP-complete {\sc clique} problem to {\sc Monophonic Position Set} was used. However, it is not clear whether {\sc Monophonic Position Set} is NP-complete for general graphs, since proving this would also require demonstrating that a solution can be verified in polynomial time. As remarked in~\cite{Thomas-2024b}, for restricted classes of graphs, this verification can be performed in polynomial time, allowing us to show the NP-completeness of the problem for these specific cases. For example, if we restrict the problem to instances $(G,k)$ with $k>|V(G)|/2$ and $G=H \vee K$, where $H$ is a generic graph and $K$ is a clique graph having the same order as $H$, the problem is NP-complete. 
In such cases, a solution can be tested in polynomial time to verify if it forms a clique and if its order exceeds $k$.

Araujo et al.~\cite[Theorem 5.2]{Araujo-2025} have independently shown that {\sc Monophonic Position Set} is NP-hard even in graphs with diameter two.  Consequently, it is W[1]-hard (parametrised by the size of the solution) and $n^{1-\epsilon}$-inapproximable in polynomial time for any $\epsilon>0$ unless $P = NP$. (For the definition of the computational complexity class W[1]-hard, see~\cite{downey-2013}.)

\subsection{All paths position number}\label{subsec:all paths}

When varying the family of paths in Definition~\ref{def:gp problem}, it is natural to consider what happens when we place no restrictions on the paths at all, i.e.\ we take the family of paths to be the set of all paths. The `all-path position problem' was investigated in one article~\cite{Haponenko-2024} from 2024 by Haponenko and Kozerenko. Most of this paper concerns the all-path convexity problem (see~\cite{Pelayo-2013}), and hence this is one example of the interplay between convexity and position problems. 

It turns out that the all-path position problem is intimately connected to connectivity. For a 2-connected graph with order $n(G) \geq 2$ the all-path position number is two (this can be shown using the Fan Lemma), whereas for a graph $G$ with a cut-vertex the all-path position number is equal to the number of leaf blocks of $G$ (i.e.\ the number of blocks of $G$ containing a single cut-vertex of $G$). This result was also proved independently in unpublished work by Sumaiyah Boshar~\cite{Boshar}. We suggest that it would be of interest to consider the more general no-$k$-in-line problem for all paths.

\subsection{Position problem with respect to cycle convexity}\label{subsec:cycle convexity}

The variations of the general position problem discussed above suggest that the concept can be naturally extended by considering arbitrary interval functions. In this direction, the general position problem with respect to cycle convexity was investigated in~\cite{Araujo-2026}. Cycle convexity was originally introduced in~\cite{Araujo-2018}. Its development was primarily motivated by applications to the study of tunnel numbers of knots and links in Knot Theory.  For a graph $G$ and a set $S \subseteq V(G)$, the interval function associated with cycle convexity is defined by $I_{cc}[S]
= S \cup\left\{u \in V(G) : G[S \cup \{u\}] \text{ contains a cycle passing through } u
\right\}.$ A set $S \subseteq V(G)$ is said to be in \emph{general position with respect to cycle convexity} if $v \notin I_{cc}[S \setminus \{v\}]$
for every $v \in S$. The corresponding \emph{cycle convexity general position number}, denoted by $\gp_{cc}(G)$, is the maximum cardinality of such a set. In~\cite{Araujo-2026} it was shown that, under cycle convexity, the general position number of a graph $G$ coincides with its \emph{forest number} $\mathrm{MIF}(G)$, that is, the maximum order of an induced forest in $G$. Consequently, $
\gp_{cc}(G)=\mathrm{MIF}(G)$.
In addition, the authors obtained closed formulas for $\gp_{cc}(G)$ for several classes of graphs, including forests, cycles, complete graphs, and complete bipartite graphs. However, the primary focus of~\cite{Araujo-2026} was the investigation of the \emph{rank} of a graph, a parameter closely related to the general position problem.

\section{New Directions}\label{sec:otherdirections}

 In this section we give an overview into some new directions taken by recent research into position problems. 

\subsection{Mobile position sets}\label{subsec:mobilegp}

Recall that the general position problem was originally introduced in~\cite{manuel-2017arxiv} using the example of a collection of robots stationed at the vertices of a graph and communicating with each other by sending signals along shortest paths. Observing delivery robots that were being trialled in Milton Keynes, home to the Open University, inspired the authors of~\cite{KlaKriTuiYer-2023} to consider a dynamic version of the general position problem, in which the robots must visit the vertices of the graph whilst keeping communication channels free.

In the scenario considered in~\cite{KlaKriTuiYer-2023}, we begin with a swarm of robots stationed at the vertices of a general position set $S$, with no more than one robot on a vertex. At each stage, exactly one robot may move from a vertex $u \in S$ to an adjacent unoccupied vertex $v$ if $(S \setminus \{ u\} ) \cup \{ v\} $ is also in general position. We call this a \emph{legal move} and denote this step by $u \move v$. We require that from the starting configuration every vertex can be visited by some robot. Such a configuration is called a \emph{mobile general position set} of $G$ and the largest possible number of robots in a mobile general position set is the \emph{mobile general position number} of $G$, denoted $\mob (G)$.

As any pair of robots is trivially in general position, we have $\mob (G) \geq 2$ for any graph with $n(G) \geq 2$ (however, there is currently no characterisation of graphs with mobile general position number two). By way of illustration, consider $K_{n_1,\dots ,n_t}$, $n_1 \geq 2$, $t \geq 2$. If at any point there are two or more robots in a partite set $V_i$, then there can be no robot in a set $V_j$, where $i \neq j$. Hence in this case for there to be an available legal move there can be only two robots on the graph. Thus if there are at least three robots on the graph, then they must all lie in different partite sets, so that there are at most $t$ robots. If there are $t$ robots, then any move would result in two robots in the same partite set, whereas it is easily seen that $t-1$ robots can visit all vertices of the graph whilst remaining in different partite sets. Thus the mobile general position number of this graph is $\max \{ 2,t-1\} $ and it follows immediately that for any $2 \leq a \leq b$ there is a graph $G$ with $\mob (G) = a$ and $\gp (G) = b$.

For cycles with length $n \neq 4,6$ the mobile general position number is $\mob (C_n) = 3$, whereas $\mob (C_4) = \mob (C_6) = 2$ (in $C_6$ a set of three vertices is in general position only if they are at distance two from one another, and such a configuration cannot be maintained by any legal move).

The paper~\cite{KlaKriTuiYer-2023} begins by discussing separable graphs, i.e.\ graphs containing a cut-vertex. Let $G$ be a graph with cut-vertex $v$. Denote the components of $G-v$ by $H_1,\dots ,H_r$ and for $i \in [r]$ the subgraph induced by $V(H_i) \cup \{ v\} $ by $G_i$. At some point a robot must visit the cut-vertex $v$, so at this point all robots must be contained within some $G_i$. Hence at any point at most two components of $G-v$ can contain robots, and if there are robots in two such components, then one of the components contains just one robot; it follows that there is a component $H_j$ that contains all the robots with at most one exception, a robot $R$ that traverses the rest of the graph. This is the source of an interesting mistake in Lemma 2.1.ii) and Corollary 2.2 of~\cite{KlaKriTuiYer-2023}, which was pointed out by Ethan Shallcross. It was assumed that whilst $R$ is visiting other components, the robots remaining in $H_j$ must form a mobile general position set, which would imply that $\mob (G) = \mob (G_k)$ for some $k \in [r]$. In fact, this is not the case, as can be seen by considering the graph $G$ obtained from a path $P_6$ by adding an edge between its second and fifth vertices. It remains an open question how much larger $\mob (G)$ can be than $\max \{ \mob (G_i):\ i \in [r]\} $. 

However, these results are not used in the rest of the paper, and the conclusion in Theorem 2.3 that $\mob (G) = \omega (G)$ for a block graph remains intact. In particular, this implies that $\mob (T) = 2$ for any tree $T$. For rooted products, it is shown that 
$$\max\{ \mob(G), \mob(H)\} \le \mob( G\circ_x H) \le \max\{ \mob(H), n(G)\}\,.$$ 
These bounds are sharp, but there are also rooted products with mobile general position number strictly between the bounds. An especially interesting case is unicyclic graphs: if $G$ is a unicyclic graph with cycle length $\ell $ and $k$ vertices on the cycle having degree at least three in $G$, then if both $k$ and $\ell $ are even, we have $\mob (G) \leq \frac{k}{2}+2$, and this bound is achieved by a family of unicyclic graphs called \emph{jellyfish} in~\cite{KlaKriTuiYer-2023}. 

The paper~\cite{KlaKriTuiYer-2023} also solves the mobile general position problem for Kneser graphs $K(n,2)$ and line graphs of complete graphs $L(K_n)$.

\begin{theorem} {\rm \cite[Theorem~3.1]{KlaKriTuiYer-2023}}
\label{thm:mobK(n,2)}
For $n \geq 5$, \[ \mob (K(n,2)) = \max \left\{ 4,\left \lfloor \frac{n-3}{2} \right \rfloor \right\}.\]
\end{theorem}
For $n \geq 11$ the robots in $K(n,2)$ can begin by occupying almost all the vertices of a maximum clique. 
\begin{theorem} {\rm \cite[Theorem~3.2]{KlaKriTuiYer-2023}}
\label{thm:mobL(K_n)}
For $n \geq 4$, \[ \mob (L(K_n)) = n-2.\]
\end{theorem}
The initial configuration of the robots on $L(K_n)$ corresponds to all the edges incident to a fixed vertex in $K_n$, with one gap left for the robots to manoeuvre.

A problem at the end of~\cite{KlaKriTuiYer-2023} asked whether there is a lower bound for the mobile general position number in terms of the clique number. This was answered in the affirmative in~\cite{Shallcross-2025}. 

\begin{theorem} {\rm \cite[Theorem~3.1]{Shallcross-2025}}
    The mobile general position number of any non-complete graph $G$ with clique number $\omega(G) \geq 3$ is bounded below by 
\[ \mob (G) \geq 1+\left \lceil \frac{\omega (G)}{\diam (G)} \right \rceil -\theta (s),\] where $\omega (G) = q\diam (G)+s$, where $0 \leq s \leq \diam (G)-1$, and $\theta (s) = 1$ if $s = 1$ and zero otherwise. 
\end{theorem}

This bound is tight for diameters two and three, but it appears that it can be improved for fixed clique number. Theorem 3.3 of the same paper shows that any graph with clique number at least five has mobile general position number at least three, but it is not known if this is best possible.

It seems difficult to approach the complexity question for mobile general position sets, and this remains an open problem. Obviously it would be desirable to investigate mobile versions of other position problems, such as monophonic position. 

One of the problems at the the end of~\cite{KlaKriTuiYer-2023} asks what happens if it is required that every robot is able to visit every vertex (whereas~\cite{KlaKriTuiYer-2023} asks only that every vertex is visited by \emph{some} robot). This is called the \emph{completely mobile general position problem} and the corresponding position number is denoted by $\mob ^*(G)$. Small examples suggest that the mobile and completely mobile general position numbers are close to each other; however a construction based on half graphs by Shallcross et al.\ in the paper~\cite{Shallcross-2025} shows that for any $2 \leq a \leq b$ there is a graph $G$ with $\mob ^*(G) = a$ and $\mob (G) = b$. Another realistic scenario to consider is when several robots are allowed to move at the same time; such sets of robots have been called \emph{hypermobile general position sets} by Shallcross.

The paper~\cite{KlaKriKuzShaTuiYer-2026} focuses on  mobile general position sets in Cartesian products, coronas and joins. For the Cartesian product we have the following two lower bounds, where interestingly in the second bound the outer general position number appears. 

\begin{proposition} {\rm \cite[Proposition~2.1]{KlaKriKuzShaTuiYer-2026}}
\label{prop:bounds-for-cp}
For any connected graphs $G$ and $H$ of order at least two, the following hold. 
\begin{enumerate}
 \item[(i)] $\mob(G\cp H) \ge \max\{\mob(G), \mob(H)\}$.
 \item[(ii)] $\mob (G \cp H) \geq \max \{ \gpo(G), \gpo(H)\}$. 
\end{enumerate}
\end{proposition}

Both lower bounds in Proposition~\ref{prop:bounds-for-cp} are sharp, as demonstrated by the next example given in~\cite[Proposition~2.2]{KlaKriKuzShaTuiYer-2026}: if $r,s \geq 2$, then 
$$\mob(K_r \cp P_s) = r = \max \{ \mob(K_r), \mob(P_s) \} = \max \{ \gpo(K_r), \gpo(P_s) \}\,.$$
The paper~\cite{KlaKriKuzShaTuiYer-2026} gives exact values of $\mob(K_n\cp K_m)$ $(n, m\ge 1)$, $\mob(P_n\cp P_m)$ $(n,m\ge 3)$, $\mob(C_n\cp K_2) $ $(n\ge 3)$, $\mob (C_4 \cp P_s) = 3$ $(r \geq 3)$, $\mob(C_r \cp P_s)$, and $(r=9\ {\rm or}\ r \geq 11\ {\rm and}\ s \geq 5)$, and proves that  $\mob (T \cp K_2) = \ell (T)$ for any tree $T$ with order at least three~\cite[Theorem 3.1]{KlaKriKuzShaTuiYer-2026}. For the corona product, lower and upper bounds are given on $\mob(G \odot H)$ in~\cite[Theorem 4.1]{KlaKriKuzShaTuiYer-2026} and it is proved that $\mob(C_n \odot K_1) = \left\lceil \frac{n}{2}\right\rceil + 1$~\cite[Theorem 4.3]{KlaKriKuzShaTuiYer-2026}. Bounds on the mobile general position number of joins $G \vee H$ are also given.

\subsection{Lower position numbers}\label{subsec:lower}

In 1976 the famous recreational mathematician (and modern day Henry Dudeney) Martin Gardner posed the `worst-case' version of Dudeney's chessboard puzzle in his column in Scientific American~\cite{Gardner-1976}: what is the smallest number of pawns that can be placed on an $n \times n$ chessboard, such that adding any new pawn creates three-in-a-line? The results can be found as sequence A277433 in the Online Encyclopedia of Integer Sequences~\cite{OEIS}. Note that the pawns in this set do not necessarily themselves have the no-three-in-line property; when the no-three-in-line property is required, this is equivalent to a smallest maximal no-three-in-line set (see sequence A219760 in OEIS). This problem is also mentioned in~\cite{AdeHolKel}. The problem is taken up, under the name \emph{geometric dominating sets}, in the article~\cite{AicEppHai}, which presents a lower bound $\Omega (n^{2/3})$ and an upper bound $2 \left \lceil \frac{n}{2} \right \rceil $ (and the discrete torus is also investigated). When the pawns are required to be no-three-in-line (which the authors call \emph{independent geometric dominating sets}), they find the solutions for $n \leq 12$, but do not improve on the trivial $2n$ upper bound.

Whilst Dudeney's puzzle involves any straight line in the plane, Gardner points out that this `minimum' puzzle remains difficult even if we restrict attention to the set of lines corresponding to rows, columns and diagonals; this is often called the Queens version of the problem. This restricted version was treated in~\cite{CooPikSchWar}, in which the authors show using Combinatorial Nullstellensatz that the answer is at least $n$, except if $n \equiv 3 \pmod 4$, in which case the answer is at least $n-1$. The authors also include some details of Gardner's correspondence on this problem (which show that this result was proven earlier by John Harris, although his argument is not extant). It is shown in~\cite{Oh} that $n+1$ is a lower bound when $n \equiv 1 \pmod 4$.

This problem was introduced into graph theory in~\cite{DiKlKrTu-2024} by defining the \emph{lower general position number} $\gp ^-(G)$ of a graph $G$ to be the number of vertices in a smallest maximal general position set of $G$. This can be viewed as the worst-case output of a greedy algorithm for finding general position sets. There is an important connection here to \emph{universal lines} in graphs, as defined in~\cite{Rodriguez-2022}; these are related to the Chen-Chv\'{a}tal Conjecture for finite metric spaces~\cite{ChenChvatal}. It turns out that the existence of a universal line in a graph $G$ is equivalent to $\gp ^-(G) = 2$; hence the lower general position number can also be seen as quantifying how far away a graph is from containing a universal line. It is easy to see that a graph $G$ has a universal line if $G$ has a bridge, $\g (G) = 2$ or $G$ is bipartite. It is shown in~\cite{Rodriguez-2022} that a block graph has a universal line if and only if one of the blocks is $K_2$. The main theorem of~\cite{Rodriguez-2022} is a characterisation of Cartesian products with a universal line.

\begin{theorem} {\rm \cite[Theorem 4.1]{Rodriguez-2022}}
Let $G$ and $H$ be non-trivial, connected graphs. Then $\gp^-(G\cp H) = 2$ if and only if one of the following conditions holds.
\begin{enumerate}
    \item[(i)] $G$ or $H$ has a maximal general position set consisting of two adjacent vertices; 
    \item[(ii)] $\g(G) = 2$ and $\g(H) = 2$.
\end{enumerate}
\end{theorem}

It is not difficult to see that the lower general position number of odd cycles is three, and for a complete $r$-partite graph with smallest part of size $t$ the lower general position number is $\min \{ r,t\} $. The lower $\gp $-numbers of Kneser graphs, line graphs of complete graphs and rook graphs are given by~\cite{DiKlKrTu-2024} as follows.

\begin{theorem} {\rm \cite[Theorem 5.1]{DiKlKrTu-2024}}\label{thm:lowergp Kneser}
If $n\ge 3$, then
\[ \gp^- (K(n,2))=\begin{cases}
3; & n \in \{3,6,7\},  \\
4; & n \in \{5,8,9\}, \\
5; & n \in \{10,11\}, \\
6; & n = 4\ \mbox{or}\ n\geq 12.
\end{cases}\]
For $n \geq 12$, a lower $\gp $-set is the set of six $2$-subsets of four symbols from $[n]$.
\end{theorem}

 \begin{theorem} {\rm \cite[Theorem 5.2]{DiKlKrTu-2024}}
If $n\ge 2$, then
      \[ \gp^- (L(K_n))=\begin{cases}
 \frac{n}{2}; & n \text{ even}, \\
\frac{n+3}{2}; & n \text{ odd.}
 \end{cases}\]
 Lower $\gp $-sets correspond to a perfect matching in $K_n$ for even $n$ and a matching plus a triangle for odd $n$.
  \end{theorem}

  \begin{theorem} {\rm \cite[Theorem 5.3]{DiKlKrTu-2024}}
\label{thm:K-r-K-s}
If $r,s\ge 2$, then
      \[ \gp^- (K_r\cp K_s)=\min\{r,s\}.\]
  \end{theorem}

An interesting problem is to compare the lower general position number with the geodetic number (recall the definition from Subsection~\ref{subsec:terminology}). A geodetic set has the property that adding any vertex creates three-in-a-line, so intuition may suggest that if a graph has a small geodetic set, then it must also have a small maximal general position set. This turns out not to be the case, and~\cite{DiKlKrTu-2024} presents constructions that demonstrate that there is a graph $G$ with $\gp^-(G) = a$ and $\g (G) = b$ if and only if $2 \leq a \leq b$ or $4 \leq b \leq a$. 

The article~\cite{DiKlKrTu-2024} also considers the largest possible size of a graph with order $n$ and given lower $\gp $-number $k$; in contrast to the extremal question for the ordinary general position number, this problem is completely solved ~\cite[Theorem 10]{DiKlKrTu-2024}, and the unique extremal graph is formed by deleting the edges of a copy of the star $S_k$ from a clique $K_n$. Regarding complexity,~\cite{DiKlKrTu-2024} shows that the independent domination problem can be polynomially reduced to the lower general position problem, so that the latter is NP-complete.

\begin{theorem} {\rm \cite[Theorem 4.2]{DiKlKrTu-2024}}\label{thm:NP}
       The {\sc Lower General Position} problem is NP-complete.
\end{theorem}

The same paper also defines the lower version of the monophonic position number, i.e.\ smallest maximal monophonic position sets. For example, any pair of vertices at distance two in the Petersen graph constitutes a lower monophonic position set. Of course for the ordinary $\gp $- and $\mono $-numbers we have $\mono (G) \leq \gp (G)$ for any graph $G$, and since any monophonic position set is in general position, it may seem that this inequality should hold for the lower versions of these numbers. However, an interesting feature of the lower position numbers is that the inequality can be reversed, and indeed the lower $\mono $-number can be arbitrarily larger than the lower $\gp $-number.

\begin{theorem} {\rm \cite[Theorem 6.1]{DiKlKrTu-2024}}
For $a,b \geq 1$, there exists a graph $G$ with $\mono^-(G) = a$ and $\gp^-(G) = b$ if and only if $a = b$, $2 \leq a < b$ or $3 \leq b < a$.
\end{theorem}

The lower general position problem for Cartesian products is treated in~\cite{Kruft-2024}. By the result of~\cite{Rodriguez-2022}, if either of the factors has a lower general position set consisting of two adjacent vertices (in particular this will hold if a factor is bipartite), then the lower general position number of the product is trivially two. If $G$ is an odd cycle or a wheel, then~\cite{Kruft-2024} shows that $\gp ^-(G \cp H) = 3$, unless $H$ has the form just mentioned. More interesting results give the lower $\gp $-numbers of Cartesian products of complete graphs and complete multipartite graphs.

\begin{theorem} {\rm \cite[Theorem 2.6]{Kruft-2024}}\label{thm:Kruftcomplete}
If $k\ge 2$ and $n_1 \ge 2, \ldots, n_k\ge 2$, then 
\[ \gp ^-(K_{n_1} \cp K_{n_2} \cp \cdots \cp K_{n_k}) = \min \{ n_1,n_2,\dots ,n_k\} .\]
\end{theorem}

\begin{theorem} {\rm \cite[Theorem 4.6]{Kruft-2024}}
If $G = K_{m_1,\dots ,m_r}$ and $H = K_{n_1,\dots ,n_s}$, where $r,s \geq 2$, $m_1 \geq \dots \geq m_r$, $n_1 \geq \dots \geq n_s$, and $m_1,n_1 \geq 2$, then 
		\[ \min \{ r,s,m_r,n_s\} \leq \gp ^-(G \cp H) \leq  \min \{ r,s,\max \{ m_r,n_s\} \}.\]
Moreover, if either $m_r = n_s$ or $\min \{ m_r,n_s\} \geq 8$, then \[ \gp ^-(G \cp H) = \min \{ r,s,m_r,n_s\} .\]
	\end{theorem}
    
The proof of Theorem~\ref{thm:Kruftcomplete} suggests the notion of \emph{orthogonal} general position sets: two (not necessarily disjoint) general position sets $S_1$ and $S_2$ are \emph{orthogonal} if any shortest path starting in $S_1$ and ending in $S_2$ contains just two vertices of the multiset $S_1 \cup S_2$. In fact the vertices of a maximal general position set in the $G$-layers of $G \cp K_r$ correspond to orthogonal general position sets. The number of orthogonal maximal general position sets that a graph can contain appears worthy of study. 

Often the easiest way to find a small maximal general position set of a product $G \cp H$ is to look for such sets contained inside a single layer; a vertex subset of a layer $G^h$ is a maximal general position set of $G \cp H$ if and only if it corresponds to a \emph{terminal set} of $G$, where a set $S \subseteq V(G)$ is {\em terminal} if $S$ is a general position set of $G$ and adding any vertex $u \in V(G)-S$ to $S$ would create three-in-a-line with $u$ as an endpoint. It is not clear that terminal sets exist for every graph. However, the authors of~\cite{Kruft-2024} conjecture that such a set does always exist.

\begin{conjecture} {\rm \cite[Conjecture 3.3]{Kruft-2024}}
    Every graph has a terminal set.
\end{conjecture}
Theorems 3.1 and 3.2 of~\cite{Kruft-2024} provide an algorithmic proof that every graph with diameter at most three has a terminal set (although it is not clear how close this algorithm is to producing a smallest terminal set), and the conjecture is also proven for cographs and chordal graphs. We close with an unexpectedly deep conjectured lower bound for $\gp ^-(G \cp H)$ in terms of the lower $\gp $-numbers of the factors, which holds for all pairs of graphs with order at most six.

\begin{conjecture} {\rm \cite[Conjecture 2.10]{Kruft-2024}}
		For any graphs $G$ and $H$, \[ \gp ^-(G \cp H) \geq \min \{ \gp ^-(G),\gp ^-(H)\} .\]
	\end{conjecture}

\subsection{Colouring problems}\label{subsec:colouring}

Graph colourings in which each colour class is required to have some special property are very common in graph theory. To give a very few examples, in \emph{proper colourings} (the foundational colouring problem) each colour class must constitute an independent set, in \emph{cocolourings} each colour class must induce either a clique or an independent set, for \emph{$2$-distance colourings} any pair of vertices with the same colour must be at distance at least two, and the \emph{domatic number} of a graph is the largest possible number of colours when each colour class is a dominating set. 

\begin{center}
			\begin{figure}[ht!]
				\centering
				\includegraphics[scale=0.75]{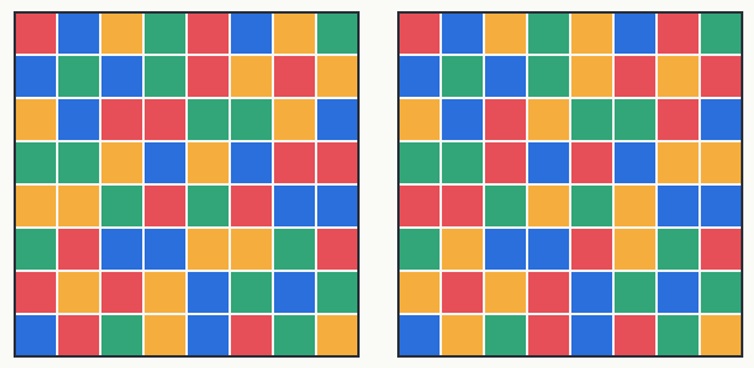}
				\caption{Two no-three-in-line colourings of the $8 \times 8$ grid with four colours found by T.~Prellberg.}
				\label{fig:8by8colouring}
			\end{figure}
		\end{center}

In discrete geometry colourings with a no-three-in-line property has long been an active area of research (see~\cite{AraujoPardo} for some references). A colouring version of Dudeney's geometrical puzzle was considered in the short note~\cite{Woodcolouring}. In this article Wood shows that if $p$ is the smallest prime $\geq n$, then the $n \times n$ chessboard can be covered by $n+p-1$ sets of pawns of different colours such that each colour class has the no-three-in-line property, from which it follows that for any $\epsilon > 0$ and sufficiently large $n$ at most $(2+\epsilon )n$ colours suffice. Two decades later, this was improved in~\cite{AraujoPardo} to $(1+\epsilon )n$ colours for any $\epsilon > 0$ and sufficiently large $n$ as a by-product of new results on no-three-in-line (or arc-proper) colourings in finite geometries. Two examples of optimal no-three-in-line colourings are shown in Figure~\ref{fig:8by8colouring}. A more general colouring conjecture for any collection of points in the plane is presented in~\cite{Payne} and some progress towards this is contained in the thesis~\cite{Paynethesis}.

It seems that the first person to suggest investigating graph colourings in which each colour class is a general position set was E.\ Sampathkumar, and this idea was taken up in the paper~\cite{ChaDiSreeThoTui2024+}. Such colourings are called $\gp$-colourings. The smallest number of colours needed for a $\gp$-colouring of $G$ is the $\emph{$\gp $-chromatic number}$ of $G$ and is written $\chi _{\gp }(G)$. For an example, in Figure~\ref{fig:Petersen} the red set on the left and the blue set on the right together constitute a $\gp $-colouring of the Petersen graph. Some simple bounds on the $\gp $-chromatic number are: 

\begin{lemma} {\rm \cite[Lemma 3.2]{ChaDiSreeThoTui2024+}}\label{lem:n/gp}
If $G$ is a graph, then 
\[ \left \lceil \frac{n(G)}{\gp (G)}\right \rceil \leq \chi _{\gp }(G) \leq \frac{n(G)-\gp (G)+2}{2}.\]
\end{lemma}

\begin{lemma} {\rm \cite[Lemma 3.6]{ChaDiSreeThoTui2024+}}\label{lem:chi_gpdiam}
If $G$ is a connected graph, then
\[ \left \lceil \frac{\diam (G)+1}{2}\right \rceil \leq \chi _{\gp }(G).\]
\end{lemma}

The lower bound in Lemma~\ref{lem:n/gp} follows from the fact that any colour class contains at most $\gp (G)$ vertices and the upper bound from a greedy colouring. The lower bound in Lemma~\ref{lem:chi_gpdiam} follows from the observation that any diametral path contains at most two vertices from each colour class. The lower bound in Lemma~\ref{lem:n/gp} is tight for paths, cycles and cliques, and Lemma~\ref{lem:chi_gpdiam} is tight for block graphs. The article also characterises graphs with very large or small values of the $\gp $-chromatic number.

As any clique is in general position, $\chi _{\gp }(G)$ is also bounded above by the \emph{clique cover number}. For similar reasons, the $\gp $-chromatic number is bounded above by the cochromatic number for graphs with diameter at most three. A more interesting bound is that for any $K_4^-$-free graph, the $\gp $-chromatic number is bounded above by the total domination number $\gamma_{\rm t}(G)$ (since any vertex neighbourhood is in general position and at most $\gamma_{\rm t}(G)$ such neighbourhoods are needed to cover $V(G)$). However, the following result shows that $\chi _{\gp }(G)$ is not tied to $\gamma_{\rm t}(G)$. 

\begin{theorem} {\rm \cite[Theorem 4.10]{ChaDiSreeThoTui2024+}}\label{thm:chi_gptotdom}
There exists a diamond-free graph $G$ with $\chi _{\gp}(G) = a$ and $\gamma_{\rm t}(G) = b$ if and only if $1 \leq a \leq b$.
\end{theorem}

The $\gp $-chromatic numbers of complete multipartite graphs, Kneser graphs and line graphs of complete graphs are determined. The most interesting of these is $L(K_n)$, for which the answer is connected to the existence of Steiner triple systems and nearly Kirkman triple systems.

\begin{theorem} {\rm \cite[Theorem 5.1]{ChaDiSreeThoTui2024+}}
If $r \geq 1$ and $n_1 \geq \dots \geq n_r$, then 
\[\chi _{\gp}(K_{n_1,\dots,n_r}) = \min \{ r,n_{r-i+1}+i-1:\ i \in [r]\} .\]
\end{theorem}

\begin{theorem} {\rm \cite[Theorem 5.3]{ChaDiSreeThoTui2024+}}
If $n \geq 5$, then $\chi _{\gp }(K(n,2)) = n-3$.
\end{theorem}

\begin{theorem} {\rm \cite[Theorem 4.10]{ChaDiSreeThoTui2024+}}
If $n \geq 3$, then
 \[ \chi _{\gp }(L(K_n)) =\begin{cases}
 \frac{n}{2}+1; & n \in \{ 6,12\},  \\
\frac{n+1}{2}; & n \equiv 1,5 \bmod 6,  \\
\frac{n}{2}; &  n \equiv 2,4 \bmod 6, \text{ or } n \equiv 0 \bmod 6 \text{ and } n \geq 18,  \\
\frac{n-1}{2}; & n \equiv 3 \bmod 6.  \\
\end{cases}\]
\end{theorem}

The article also treats the $\gp $-chromatic number of Cartesian products of paths and cycles. When both factors are large, maximum general position sets nearly tessellate the vertex set of the graph and the lower bound in Lemma~\ref{lem:n/gp} holds asymptotically. In the case of Theorem~\ref{thm:chigpgrid}, the upper bound follows from the total domination upper bound in Theorem~\ref{thm:chi_gptotdom}. A colouring of the type referred to in Theorem~\ref{thm:4byngrid} is given in Figure~\ref{fig:P_4 cp P_n}.

\begin{theorem} {\rm \cite[Theorem 7.4]{ChaDiSreeThoTui2024+}}\label{thm:4byngrid}
If $n \geq 4$, then $\chi _{\gp }(P_4 \cp P_n) = n+1$.
\end{theorem}

\begin{theorem} {\rm \cite[Theorem 7.5]{ChaDiSreeThoTui2024+}}\label{thm:chigpgrid}
If $n_1,n_2 \geq 16$, then 
\[ \frac{n_1n_2}{4} \leq \chi _{\gp }(P_{n_1} \cp P_{n_2} ) \leq \left \lfloor \frac{(n_1+2)(n_2+2)}{4} \right \rfloor -4. \]
\end{theorem}

\begin{theorem} {\rm \cite[Theorem 7.9]{ChaDiSreeThoTui2024+}}
    When $n_1 \geq 5$ and $n_2 = 7$ or $n_2 \geq 9$, \[ \chi _{\gp }(P_{n_1} \cp C_{n_2}) \leq \frac{n_1n_2}{5}+O(n_1+n_2).\]
\end{theorem}

\begin{theorem} {\rm \cite[Theorem 7.10]{ChaDiSreeThoTui2024+}}
    When $s \geq t \geq 7$, \[ \chi _{\gp }(C_{7s} \cp C_{7t}) = 7st.\] 
\end{theorem}
When the orders of the factors satisfy certain divisibility conditions, we can be more precise. For example, if $n_1 \geq 3$ is odd and $12|n_2$, then $\chi _{\gp }(P_{n_1} \cp P_{n_2}) = \frac{n_1n_2}{4}+\frac{n_2}{12}$ (Theorem 7.6 of~\cite{ChaDiSreeThoTui2024+}).

\begin{figure}[ht!]
\centering
		\begin{tikzpicture}[x=0.2mm,y=-0.2mm,inner sep=0.1mm,scale=1.30,
			thick,vertex/.style={circle,draw,minimum size=10,font=\tiny,fill=white},edge label/.style={fill=white}]
			\tiny
			
			\node at (-125,-75) [vertex,color=yellow] (11) {};
			\node at (-75,-75) [vertex,color=red] (12) {};
			\node at (-25,-75) [vertex,color=blue] (13) {};
			\node at (25,-75) [vertex,color=green] (14) {};
			\node at (75,-75) [vertex,color=pink] (15) {};
			\node at (125,-75) [vertex,color=yellow] (16) {};
			
			\node at (-125,-25) [vertex,color=red] (21) {};
			\node at (-75,-25) [vertex,color=blue] (22) {};
			\node at (-25,-25) [vertex,color=green] (23) {};
			\node at (25,-25) [vertex,color=pink] (24) {};
			\node at (75,-25) [vertex,color=yellow] (25) {};
			\node at (125,-25) [vertex,color=orange] (26) {};
			
			\node at (-125,25) [vertex,color=orange] (31) {};
			\node at (-75,25) [vertex,color=brown] (32) {};
			\node at (-25,25) [vertex,color=red] (33) {};
			\node at (25,25) [vertex,color=blue] (34) {};
			\node at (75,25) [vertex,color=green] (35) {};
			\node at (125,25) [vertex,color=pink] (36) {};
			
			\node at (-125,75) [vertex,color=brown] (41) {};
			\node at (-75,75) [vertex,color=red] (42) {};
			\node at (-25,75) [vertex,color=blue] (43) {};
			\node at (25,75) [vertex,color=green] (44) {};
			\node at (75,75) [vertex,color=pink] (45) {};
			\node at (125,75) [vertex,color=brown] (46) {};
			
			\path
			
			(11) edge (12)
			(12) edge (13)
			(13) edge (14)
			(14) edge (15)
			(15) edge (16)
		
			(21) edge (22)
			(22) edge (23)
			(23) edge (24)
			(24) edge (25)
			(25) edge (26)
			
			(31) edge (32)
			(32) edge (33)
			(33) edge (34)
			(34) edge (35)
			(35) edge (36)
			
			(41) edge (42)
			(42) edge (43)
			(43) edge (44)
			(44) edge (45)
			(45) edge (46)
			
			(11) edge (21)
			(21) edge (31)
			(31) edge (41)
			
			(12) edge (22)
			(22) edge (32)
			(32) edge (42)
			
			(13) edge (23)
			(23) edge (33)
			(33) edge (43)
			
			(14) edge (24)
			(24) edge (34)
			(34) edge (44)
			
			(15) edge (25)
			(25) edge (35)
			(35) edge (45)
			
			(16) edge (26)
			(26) edge (36)
			(36) edge (46)

			;

		\end{tikzpicture}
		\caption{A $\gp $-colouring of a $4 \times 6$ Cartesian grid}
		\label{fig:P_4 cp P_n}
	\end{figure}

However, when one of the factors is small such a `nearly perfect' tessellation is not always possible.

\begin{theorem} {\rm \cite[Theorem 7.2]{ChaDiSreeThoTui2024+}}
If $n \geq 3$, then 
\[ \chi _{\gp }(P_2 \cp P_n) =\begin{cases}
2r; &  n = 3r,  \\
2r+1; &  n = 3r+1, \\
2r+2; & n = 3r+2. \\
\end{cases}\]
\end{theorem}

\begin{theorem} {\rm \cite[Theorem 7.3]{ChaDiSreeThoTui2024+}}
If $n\ge 1$, then 
\[ \chi _{\gp }(P_3 \cp P_n) = \frac{5n}{6}+O(1).\] 
Moreover, if $12|n$, then $\chi _{\gp }(P_3 \cp P_n) = \frac{5n}{6}$.
\end{theorem} 

Finding the $\gp $-chromatic number is, predictably, a computationally hard problem.

\begin{theorem} {\rm \cite[Theorem 8.2]{ChaDiSreeThoTui2024+}}\label{theo:gpc-NP}
 Deciding if a graph $G$ has a $\gp $-colouring with $\leq k$ colours is NP-complete even for instances $(G,k)$ with $\diam(G) = 2$ and $k=3$.  
\end{theorem} 

The article~\cite{ChaDiSreeThoTui2024+} also discusses the analogous chromatic numbers corresponding to colour classes that are in monophonic position, independent general position and independent monophonic position (denoted by $\chi_{\mono}(G)$, $\chi_{\gp_{\rm i}}(G)$ and $\chi_{\mono_{\rm i}}(G)$ in the paper respectively). One interesting question is finding the largest size of graphs with given order and position chromatic numbers. For independent general or monophonic position colourings the extremal graphs are given by Tur\'{a}n's Theorem. By contrast, graphs with fixed $\gp $- or $\mono $-chromatic number can have almost all edges present.

\begin{theorem} {\rm \cite[Theorem 4.12]{ChaDiSreeThoTui2024+}}
For $a \geq 2$, the largest size of a graph $G$ of order $n$ with $\chi _{\gp }(G) = a$ is at least ${n \choose 2}-\frac{(a-1)a(a+1)}{6}$ for $n \geq \frac{a(a+1)}{2}$ and the largest size of a graph with $\chi _{\mono }(G) = a$ is at least ${n \choose 2}-(a-2)(2a-1)$ for $n \geq 2a$.
\end{theorem}
The exact values of the extremal sizes remains an open question. Another open problem is given in this conjectured Nordhaus-Gaddum relation.

\begin{conjecture} {\rm \cite[Conjecture 9.3]{ChaDiSreeThoTui2024+}}\label{conj:NGrelation}
 If $G$ is a graph and $\pi \in \{\gp, \mono\}$, then $\chi _{\pi_{\rm i}}(G) +\chi_{\pi_{\rm i}}(\overline{G}) \leq n(G)+1$. 
\end{conjecture}

In~\cite{Cody-2026} general position colourings were generalised to colourings in which each each colour class is a $k$-general $d$-position set (see Subsection~\ref{subsec:d-position-sets} for discussion of the $k$-general $d$-position problem, although we note that in~\cite{Cody-2026} these subsets are called $k,d$-independent sets and use a different notation). The smallest number of colours needed to colour a graph $G$ is the \emph{$k,d$-chromatic number} $\chi ^k_d(G)$ of $G$. Two greedy bounds are derived in~\cite{Cody-2026}. The first is in terms of graph powers: $\chi ^k_D(G) \leq \Delta (G^{d-k+2})+1$ for $ 1 \leq k-1 \leq d$. The second is $\chi ^k_d(G) \leq \Delta (\mathbb{H}_{k,d}(G))+1$, where $\mathbb{H}_{k,d}(G)$ is the $k$-uniform hypergraph derived from $G$ with vertex set $V(G)$ and hyperedges consisting of all $k$-subsets that are contained on a common shortest path of $G$ of length $\leq d$.

An elementary bound for a graph $G$ with order $n$ is \[ \max \left \{ \left \lceil \frac{n}{\gp ^k_d(G)} \right \rceil , \left \lceil \frac{\omega ^k_d(G)}{k-1} \right \rceil \right \} \leq \chi ^k_d(G) \leq \left \lceil \frac{n-\gp ^k_d(G)}{k-1} \right \rceil +1 , \] where $\omega ^k_d(G)$ is the largest number of vertices in a subset $S \subseteq V(G)$ such that every $k$-element subset of $S$ is contained in a shortest path of length at most $d$. Inspired by perfect graphs, the authors define a graph $G$ to be \emph{$k,d$-perfect} if $\chi ^k_d(H) = \left \lceil \frac{\omega ^k_d(H)}{k-1} \right \rceil $ for every induced subgraph $H$ of $G$. For example, it is shown that chordal graphs are $2,d$-perfect for odd $d$, and strongly chordal graphs (i.e.\ a chordal graph with the property that every even cycle of length $\geq 6$ has a chord) are $2,d$-perfect for all $d \geq 1$. The latter class includes block graphs, and in particular trees, and it is also proven that trees are $3,2$-perfect.

The article~\cite{Cody-2026} finds the exact value of the $k,d$-chromatic number for paths, cycles and powers of paths. All finite paths and the two-way infinite path are $k,d$-perfect for $d+1 \geq k \geq 2$, and a more complicated classification of $k,d$-perfect cycles is provided. The $k,d$-chromatic number for powers of cycles and the question of $k,d$-perfection for powers of paths and cycles are left open. The article suggests investigating these problems in other families, such as hypercubes and Cartesian products of paths and cycles. An interesting problem raised is whether there is an analogue of the Strong Perfect Graph Theorem for $k,d$-perfect graphs.

\subsection{Fractional position problems}\label{subsec:fractional}

The fractional version of the general position problem turns out to be of great interest.

\begin{definition}
  A \emph{general position labelling} of a connected graph $G$ with $n(G) \geq 2$ is a non-negative real-valued function $f$ on $V(G)$ such that for any shortest path $P$ of $G$, the sum $\sum _{u \in V(P)}f(u)$ is at most two. The \emph{fractional general position number} $\gp _f(G)$ of $G$ is the largest possible sum of labels of a general position labelling of $G$.  
\end{definition}

A fractional $\gp $-labelling is obtained by assigning label 1 to every vertex of a $\gp $-set and 0 to every other vertex. Furthermore, as the sum of the labels along any shortest path is at most two, it follows that $2\rho (G)$ is an upper bound for $\gp _f(G)$. This yields the following pleasing bound.

\begin{lemma}
For any graph $G$, the fractional $\gp $-number is bounded by \[ \gp (G) \leq \gp _f(G) \leq 2\rho (G).\]
\end{lemma}

An equivalent problem was considered in~\cite{Fitzpatrick-2006} as the dual problem of the `fractional isometric path number', and without any connection to the general position problem. For technical reasons our definition of fractional $\gp $-number is always double the dual of the fractional isometric path number, so we take the liberty of multiplying all results of~\cite{Fitzpatrick-2006} by two. Note that the fractional general position problem can be construed as a linear program, and hence the fractional general position number is guaranteed to be a rational number. Several of the results below were proven using properties of the dual program.

Firstly, since any shortest path has at most $\diam (G)+1$ vertices, one can assign every vertex of a graph the label $\frac{2}{\diam (G)+1}$ to give the following lower bound.

\begin{lemma} {\rm \cite[Theorem 1.1]{Fitzpatrick-2006}}
    If $G$ is a graph, then $\gp _f(G) \geq \frac{2n(G)}{\diam (G)+1}$.
\end{lemma}

\begin{theorem} {\rm \cite[Theorem 1.2]{Fitzpatrick-2006}}
 If $n\ge 2$, then $\gp _f(K_n) = n$, and if $n \geq 3$, then 
 \[ \gp _f(C_n) = 
\begin{cases}
 \frac{4n}{n+2}; & \text{$n$ is even}\,, \\
 \frac{4n}{n+1}; & \text{$n$ is odd}\,.
 \end{cases}\]
\end{theorem}

\begin{theorem} {\rm \cite[Theorem 1.3]{Fitzpatrick-2006}}
If $T$ is a tree, then $\gp _f(T) = \gp (T) = \ell (T)$.
\end{theorem}

\begin{theorem} {\rm \cite[Theorem 1.4]{Fitzpatrick-2006}}
If $n\ge 1$, then $\gp _f(Q_n) = \frac{2^{n+1}}{n+1}$.
\end{theorem}

The majority of the paper focusses on the square grid $P_n \cp P_n$. Finding the fractional general position number of this graph turns out to be a deep problem. The following lower bound, which is the best known construction for even $n$, is obtained by assigning label 0 to every vertex at distance $< t$ from a corner and label $\frac{2}{2n-2t-1}$ to the remaining vertices, for some integer in the range $1 \leq t \leq \frac{n}{2}$. 

\begin{theorem} {\rm \cite[Lemma 2.1]{Fitzpatrick-2006}}\label{thm:fitzpatrick1}
 For $n \geq 3$ and $0 \leq t \leq \frac{n}{2}$,
 \[ \gp _f(P_n \cp P_n) \geq \frac{n^2-2t(t+1)}{n-t-\frac{1}{2}}.\]
\end{theorem}
This construction can be improved on for odd $n$ by labellings in which vertices within a given distance from a corner are assigned label 0, and for the remaining vertices the zero and non-zero labels create a chequerboard pattern.

\begin{theorem} {\rm \cite[Lemma 2.2]{Fitzpatrick-2006}}\label{thm:fitzpatrick2}
   For odd $n \geq 3$, 
   \begin{itemize}
       \item $\gp _f(P_n \cp P_n) \geq \frac{n^2+1-8k^2}{n-2k}$ for $0 \leq k \leq m$, where $m = \frac{n-1}{4}$ when $n \equiv 1 \pmod 4$ and $m = \frac{n+1}{4}$ when $n \equiv 3 \pmod 4$,
       \item $\gp _f(P_n \cp P_n) \geq \frac{n^2-1-8(k^2+k)}{n-2k-1}$ for $0 \leq k \leq m$, where $m = \frac{n-1}{4}$ when $n \equiv 1 \pmod 4$ and $m = \frac{n-3}{4}$ when $n \equiv 3 \pmod 4$.
   \end{itemize}
\end{theorem}

The upper bound given in~\cite{Fitzpatrick-2006} is quite complex, but it is asymptotically equal to $2n$. Hence the ratio of the upper and lower bounds tends to approximately 1.0082 as $n \rightarrow \infty $. The paper~\cite{Fitzpatrick-2006} gives exact values only for $n \leq 7$; by making use of the symmetry of the grid, it is shown in~\cite{fractional} that the lower bound is best possible for $n \leq 10$. Some values are shown in Table~\ref{Tab:gp_fG_n}. We therefore make: 

\begin{conjecture}
    The lower bound from Theorem~\ref{thm:fitzpatrick1} is optimal for even $n$ and the bound from Theorem~\ref{thm:fitzpatrick2} is optimal for odd $n$. 
\end{conjecture}

 \begin{table}[ht!]
\centering
\renewcommand{\arraystretch}{2}
\begin{tabular}{|c| c| c|c|c|c|c|c|c|c|c|c|c|c|c|c|}
\hline
$n$ & 3 & 4 & 5 & 6 & 7 & 8 & 9 & 10 & 11 & 12 & 13 & 14 & 15 & 16 & 17 \\ 
\hline 
Lower & 4 & $\frac{24}{5}$ & 6 & $\frac{64}{9}$ & $\frac{42}{5}$ & $\frac{104}{11}$ & $\frac{32}{3}$ & $\frac{176}{15}$ & $13$ & $\frac{240}{17}$ & $\frac{46}{3}$ & $\frac{312}{19}$ & $\frac{194}{11}$ & $\frac{432}{23}$ & 20 \\ 
\hline
Upper & 4 & $\frac{24}{5}$ & 6 & $\frac{64}{9}$ & $\frac{42}{5}$ & $\frac{104}{11}$* & $\frac{32}{3}$* & $\frac{176}{15}$* & 14 & $\frac{128}{9}$ & $\frac{1602}{103}$ & $\frac{19942}{1198}$ & 18 & $\frac{96}{5}$ & $\frac{101318}{5001}$\\
\hline
\end{tabular}
\renewcommand{\arraystretch}{1}
\caption{Best known bounds for $\gp _f(P_n \cp P_n)$ from~\cite{Fitzpatrick-2006}. Entries with a * are from~\cite{fractional}.}
\label{Tab:gp_fG_n}
\end{table}

We note that~\cite{AraujoPardo} introduces a fractional version of arc-proper colourings in finite geometry. This suggests investigating fractional versions of the problems discussed in the preceding subsection.

\subsection{Games on graphs}\label{subsec:games}

When subsets of a graph exhibit a special property, a common question is how to identify optimal subsets with this property, particularly when they emerge through adversarial play. The articles~\cite{AicEppHai} and~\cite{Gardner-1976} suggested looking at such games in the geometric No-Three-In-Line Problem. Klav\v{z}ar, Neethu and Chandran~\cite{KlaNeeCha-2022} introduced two games on graphs in 2021: the general position achievement game and the general position avoidance game ({\em gp-achievement} and {\em gp-avoidance} for short).

In both games, two players (A and B) take turns selecting free vertices of a graph $G$ such that at any time the set of selected vertices is in general position; the first player unable to move is the loser in the achievement game, but is the winner in the avoidance game. In game terminology, the $\gp$-achievement game corresponds to the normal game (where the last player to move wins), while the $\gp$-avoidance game corresponds to the \textit{mis\`ere} game (where the last player to move loses). By the classical Zermelo-von Neumann Theorem, in both games one of the two players must have a winning strategy, as these are finite, perfect-information games without the possibility of a draw.

The key question, therefore, is: given a graph in one of these general position games, which player has a winning strategy? It is important to note that the two games are independent, meaning that winning the $\gp$-achievement game on a given graph does not imply that the same player will lose the $\gp$-avoidance game on that graph. Additionally, both games result in a maximal general position set. This suggests that positional games provide an effective means of constructing maximal general position sets. The paper~\cite{KlaSam} studies this problem further, especially with regard to the complexity of the problem, see Subsubsection~\ref{subsubsec:complexity of games}. 

The sequence of vertices played in the position games on a graph $G$ is denoted by $a_1, b_1, a_2, b_2, \ldots$. The vertices played by player A are $a_1, a_2, \ldots$, and the vertices played by player B are $b_1, b_2, \ldots$. For example, we can say that player A starts the game by playing $a_1 = x$, where $x \in V(G)$. Suppose that $x_1, \ldots, x_j$ are the vertices played so far on $G$. We say that $y \in V(G)$ is a {\em playable vertex} if $y \notin \{x_1, \ldots, x_j\}$ and the set $\{x_1, \ldots, x_j\} \cup \{y\}$ is a general position set of $G$. Let $ \Pl_G(\ldots x_j)$ represent the set of all playable vertices after the vertices $x_1, \ldots, x_j$ have already been played. The following implicit strategies for the playable vertices, as proven in~\cite{KlaSam, KlaNeeCha-2022}, are useful in many contexts.

\begin{theorem}{\rm \cite[Theorem 2.2]{KlaNeeCha-2022}}
\label{thm:up-to-gp-set1} 
Let $G$ be a graph. Then the following holds. 

(i) If A has a strategy such that after the vertex $a_k$, $k\ge 1$, is played, the set $\Pl_G(\ldots a_k)\cup \{a_1, b_1, \ldots, a_k\}$ is a general position set and $|\Pl_G(\ldots a_k)|$ is even, then A wins the $\gp$-achievement game. 

(ii) If B has a strategy such that after the vertex $b_k$, $k\ge 1$, is played, the set $\Pl_G(\ldots b_k)\cup \{a_1, b_1, \ldots, b_k\}$ is a general position set and $|\Pl_G(\ldots b_k)|$ is even, then B wins the $\gp$-achievement game. 
\end{theorem}

\begin{theorem}{\rm \cite[Theorem 2.4]{{KlaSam}}}\label{thm:up-to-gp-set}
Let $G$ be a graph. Then the following holds. 

(i) If A has a strategy such that after the vertex $a_k$, $k\ge 1$, is played, the set $\Pl_G(\ldots a_k)\cup \{a_1, b_1, \ldots, a_k\} $ is a general position set and $|\Pl_G(\ldots a_k)|$ is odd, then A wins the $\gp $-avoidance game. 

(ii) If B has a strategy such that after the vertex $b_k$, $k\ge 1$, is played, the set $\Pl_G(\ldots b_k)\cup \{a_1, b_1, \ldots, b_k\}$ is a general position set and $|\Pl_G(\ldots b_k)|$ is odd, then B wins the $\gp$-avoidance game. 
\end{theorem}

As a consequence of Theorem~\ref{thm:up-to-gp-set1}, the $\gp$-achievement game is solved for complete multipartite graphs and bipartite graphs. 

\begin{proposition}{\rm \cite[Theorem 2.3]{KlaNeeCha-2022} }
\label{prp:complete-multi}
If $k\ge 2$ and $n_i\ge 2$, $i\in [k]$, then A wins the $\gp $-achievement game on $K_{n_1,\ldots, n_k}$ if and only if $k$ is odd and at least one $n_i$ is odd.  
\end{proposition}

\begin{theorem}{\rm \cite[Theorem 2.5]{KlaNeeCha-2022} }
\label{thm:bipartite} 
Let $G$ be a bipartite graph. Then A wins the $\gp $-achievement game on $G$ if and only if the number of isolated vertices in $G$ is odd.
\end{theorem}

As an application of Theorem~\ref{thm:up-to-gp-set}, the $\gp $-avoidance game is solved in~\cite[Theorem 2.5]{KlaSam} for complete multipartite graphs. 

\begin{proposition}\label{prp:complete-bipartite}
If $k\ge 2$ and $n_i\ge 2$, $i\in [k]$, then A wins the $\gp$-avoidance game on $K_{n_1,\ldots, n_k}$ if and only if $k$ is even and at least one $n_i$ is even.  
\end{proposition}

For the generalised wheel $W_{n,m} = \overline{K_n} \vee C_m$, we have the following result for the $\gp $-avoidance game. This contrasts with the observation that solving the $\gp $-achievement game on $W_{ 
n,m}$ appears to be challenging.

\begin{theorem}{\rm \cite[Theorem 2.6]{KlaSam}}
\label{thm:4.2}
If $n\ge 1$ and $m\ge 3$, then B wins the $\gp$-avoidance game on $W_{n,m}$ if and only if $m\ge 4$. 
\end{theorem}

In~\cite{KlaNeeCha-2022}, the $\gp$-achievement game is further investigated on Cartesian and lexicographic products. This game is resolved for Cartesian products when one factor is a bipartite graph and for rook graphs.

\begin{theorem}{\rm \cite[Theorem 3.6] {KlaNeeCha-2022}}
\label{thm:p-bipartite} 
Let $G$ be a connected graph and let $H$ be a connected bipartite graph with at least one edge. Then B wins the $\gp $-achievement game on $G\cp H$. 
\end{theorem}

\begin{theorem} {\rm \cite[Theorem 3.7] {KlaNeeCha-2022}}
\label{thm:3.11}
If $n, m\ge 2$, then A wins the $\gp $-achievement game on $K_n\cp K_m$ if and only if both $n$ and $m$ are odd.
\end{theorem}

As remarked in~\cite[Corollary 3.5]{KlaNeeCha-2022}, if $n$ is even and $G$ is a connected graph, then B wins the $\gp $-achievement game on $K_n\cp G$. On the other hand, if $n$ is odd and $G$ is a connected graph, the outcome of the $\gp $-achievement game on $K_n\cp G$ appears to be difficult to determine. This statement is justified by the proof of the following result. 

\begin{theorem}{\rm \cite[Theorem 3.8]{KlaNeeCha-2022}}
If $m\ge 3$, then A wins the $\gp $-achievement game on $K_3 \cp C_m$ if and only if $m\in \{3,5\}$.
\end{theorem}

For the lexicographic product, we have: 
 
\begin{theorem} {\rm \cite[Theorem 4.3]{KlaNeeCha-2022}}
If $G$ is a connected graph, then B wins the $\gp$-achievement game on $G\circ K_n$ if and only if B wins on $G$ or $n$ is even.
\end{theorem}

The paper~\cite{KlaSam} presents several additional results for the $\gp $-avoidance game for Cartesian and lexicographic products.
We recall that player A wins the $\gp$-achievement game on the rook graph $K_n\cp K_m$ if and only if both $n$ and $m$ are odd.  For the $\gp$-avoidance game, the following holds. 

 \begin{theorem}{\rm \cite[Theorem 4.4]{KlaSam}}\label{thm:p-complete} If $m \geq n \geq 2$, then B wins the $\gp$-avoidance game on $K_n\cp K_m$ if and only if either $n=2$ and $m$ is odd, or $n=3$ and $m$ is even.
\end{theorem}

Recall that player B always wins the $\gp$-achievement game for any connected bipartite graph. However, solving the $\gp$-avoidance game for arbitrary connected bipartite graphs presents significant challenges. This difficulty is partly supported by results obtained for grids and cylinders.

\begin{theorem}{\rm \cite[Theorem 4.5]{KlaSam}}
\label{theorem:4.4} If $n\geq 3$ and $m\geq 2$, then B wins the $\gp$-avoidance game on $P_n\cp P_m$.
\end{theorem}

\begin{theorem}{\rm \cite[Theorem 4.7]{KlaSam}}\label{theorem:3.13b}
If $n\geq 3$ and $m\geq 2$, then B wins the $\gp$-avoidance game on $C_n\cp P_m$ if and only if $n$ is odd.
\end{theorem}

For the lexicographic product, the next result holds for the $\gp$-avoidance game. The game is solved for a few more classes in~\cite{KlaSam}.

\begin{theorem} {\rm \cite[Theorem 5.4]{KlaSam}}
If $G$ is a connected graph and $n\ge 1$, then B wins the $\gp$-avoidance game on $G\circ K_n$ if and only if B wins the $\gp$-avoidance game on $G$ and $n$ is odd.
\end{theorem}

A different approach is taken in~\cite{KlavzarTianTuite-2026}, inspired by the domination game. In this \emph{Builder/Blocker} game, two players, called Builder and Blocker, take it in turns to add a new vertex to a general position set and play finishes when no further vertices can be added without creating three-in-a-line. However the goals of the players are diametrically opposed: Builder wishes the resulting general position set to be as large as possible, whereas Blocker wants to keep the set as small as possible. If Builder starts the game, this is called the B-game, but if Blocker starts it is the B'-game. The number of vertices in the set built by optimal play in the B-game on a graph $G$ is the \emph{Builder-game general position number} $\gp_{\rm g}(G)$ and in the B'-game it is the \emph{Blocker-game general position number} $\gp_{\rm g}'(G)$. 

Note that, whilst the achievement and avoidance $\gp $-games can result in sets of different sizes, the number of vertices in the output of the Builder-Blocker games is always the same. As the game finishes with a maximal general position set, both $\gp_{\rm g}(G)$ and $\gp_{\rm g}'(G)$ will lie between $\gp ^-(G)$ and $\gp (G)$. In fact Corollary 2.8 of~\cite{KlavzarTianTuite-2026} shows that for any triple $2 \leq c \leq b \leq a$ there exists a graph $G$ with $\gp (G) = a$, $\gp_{\rm g}(G) = b$ and $\gp ^-(G) = c$. 

The article determines these game numbers for different families of graphs. For example, if $t\ge 2$ and $n_1 \geq \cdots \geq n_t\ge 2$, then $\gp_{\rm g}(K_{r_1,\dots,r_t}) = \min \{ r_1,t\} $ and $\gp_{\rm g}'(K_{r_1,\dots,r_t}) = \max \{ r_t,t\} $. For Kneser graphs we have \[\gp_{\rm g}(K(n,2)) = \gp_{\rm g}'(K(n,2)) = 6,\] so that for $n \geq 12$ the game general position numbers coincide with the lower general position number. It is shown that for trees $\gp_{\rm g}'(T)\le \ell(T) - \Delta(T) + 2$ and the equality case is characterised. However, the result for $L(K_n)$ and $K(n,k)$ for $k > 2$ is at present unknown.

It is easily seen that a graph satisfies $\gp_{\rm g}(G) = 2$ if and only if every vertex is contained in a maximal general position set of order two. This is true, for example, of any bipartite graph or any vertex-transitive graph with $\gp ^-(G) = 2$. It would be desirable to have a more elementary characterisation of these graphs. A graph with $\gp_{\rm g}'(G) = 2$ must have a special vertex $u$ that Blocker can take advantage of such that $\{ u,v\} $ yields a universal line for any $v \in V(G \setminus \{ u\}$; these graphs are characterised in~\cite{KlavzarTianTuite-2026} and an example is shown in Figure~\ref{fig:gpg=2}. Every member of this family also satisfies $\gp_{\rm g}(G) = 2$, so $\gp_{\rm g}'(G) = 2$ implies $\gp_{\rm g}(G) = 2$. 

\begin{figure}[ht!]
		\centering
		\begin{tikzpicture}[x=0.2mm,y=-0.2mm,inner sep=0.2mm,scale=0.7,thick,vertex/.style={circle,draw,minimum size=10}]
			\node at (0,0) [vertex] (v1) {$u$};
                \node at (100,-100) [vertex] (v2) {};
                \node at (100,-50) [vertex] (v3) {};
                \node at (100,0) [vertex] (v4) {};

                \node at (100,50) [vertex] (v5) {};
                \node at (100,100) [vertex] (v6) {};

               \node at (200,-50) [vertex] (v7) {};
               \node at (200,75) [vertex] (v8) {};

               \node at (300,-50) [vertex] (v9) {};
               \node at (400,-50) [vertex] (v10) {};
               \node at (500,-50) [vertex] (v11) {};

               \node at (-100,-75) [vertex] (v12) {};
               \node at (-100,-25) [vertex] (v13) {};
               \node at (-100,25) [vertex] (v14) {};
               \node at (-100,75) [vertex] (v15) {};

                \node at (-200,-75) [vertex] (v16) {};
                \node at (-300,-75) [vertex] (v17) {};
                \node at (-200,-25) [vertex] (v18) {};

			\path

			(v1) edge (v2)
                (v1) edge (v3)
                (v1) edge (v4)
                (v1) edge (v5)
                (v1) edge (v6)
                
                (v7) edge (v2)
                (v7) edge (v3)
                (v7) edge (v4)

                (v8) edge (v5)
                (v8) edge (v6)
			
			(v7) edge (v9)
                (v9) edge (v10)
                (v10) edge (v11)

                (v1) edge (v12)
                (v1) edge (v13)
                (v1) edge (v14)
                (v1) edge (v15)

                (v12) edge (v16)
                (v16) edge (v17)
                (v13) edge (v18)
			
			;
		\end{tikzpicture}
		\caption{A graph with Blocker game general position number two.}
		\label{fig:gpg=2}
\end{figure}

An unexpectedly interesting question is: does it matter who goes first? Based on the domination game and its `Continuation Principle', see~\cite{bresar-2010, book-2021}, intuition may suggest that the order of the players does not make a big difference. This turns out to be completely false. We call a pair $(a,b)$ \emph{realisable} if there is a graph $G$ with $\gp_{\rm g}(G) = a$ and $\gp_{\rm g}'(G) = b$. Proving that pairs with $2 \leq a \leq b$ are realisable is quite simple, but for $a < b$ the problem is harder. A construction involving identifying 4-cycles and an odd number of 3-cycles along an edge shows that $(a,b)$ is realisable if $a > b$ and $b \geq 3$ is odd. Some such pairs with even $b$ were also found. The article therefore raised the question: is every pair $(a,b)$ with $a > b > 2$ and even $b$ realisable? This question was completely settled in the affirmative in a forthcoming article~\cite{Evans}. Finding the smallest graphs that realise a given pair also seems a worthwhile problem.

\subsubsection{Complexity of games}\label{subsubsec:complexity of games}

The algorithmic complexity of deciding the winning player of the general position achievement/avoidance games of a given graph was studied in~\cite{KlaSam}. It was proved that the complexity of deciding the winning player of the general position achievement/avoidance games is in the class of PSPACE-complete problems even in graphs with diameter at most $4$. Formally, they consider the games as decision problems: \emph{ Given a graph $G$, does player A have a winning strategy?}

The paper~\cite{KlaSam} presents a reduction to the general position achievement game from the clique-forming game, which is PSPACE-complete (see~\cite{schaefer-1978} for the complexity result on this problem). In this game, two players alternately select vertices from a graph $G$, ensuring the selected vertices form a clique. The player who plays the last vertex of a maximal clique wins.
This game is similar to Node Kayles, which is PSPACE-complete (see~\cite{schaefer-1978}). In Node Kayles, the goal is to form an independent set. The clique-forming game is Node Kayles on the complement of the graph, and vice versa.

\begin{theorem}{\rm \cite[Theorem 3.2]{KlaSam}}\label{teo-pspace1}
The game $\gp$-achievement is PSPACE-complete even on graphs with diameter at most four.
\end{theorem}

An analogous result, that is, PSPACE-completeness, is true for the $\gp$-avoidance game even on graphs of diameter at most four. By using a similar setting as for the $\gp$-achievement game, a reduction to the $\gp$-avoidance games  from the \textit{mis\`ere} clique-forming game is presented in~\cite{KlaSam}. For this purpose, the PSPACE-hardness of the \textit{mis\`ere} clique-forming game and the \textit{mis\`ere} Node Kayles game was used. The \textit{mis\`ere} clique-forming game is the same as the clique-forming game, but the first player unable to play wins the game. The \textit{mis\`ere} clique-forming game is strongly related to the \textit{mis\`ere} Node Kayles game, whose objective is to obtain an independent set and the first player unable to play wins the game.

\begin{theorem}{\rm\cite[Theorem 3.3]{KlaSam}}\label{teo-pspace2}
The game $\gp$-avoidance is PSPACE-complete even in graphs with diameter at most four.
\end{theorem}

The complexity of the decision version of the Builder/Blocker game is completely unknown, but~\cite{KlavzarTianTuite-2026} gives this conjecture: 

\begin{conjecture} {\rm \cite[Conjecture 5.5]{KlavzarTianTuite-2026}}
The decision versions of the B- and B'-games are PSPACE-complete.
\end{conjecture}

\subsection{General position in directed graphs}\label{subsec:digraphs}

The paper~\cite{digraphs} considers the general position problem in the setting of directed graphs (for notation, we refer to the latter paper, or to the encyclopedic text~\cite{BangJensenGutin}). By a natural extension of the definition for undirected graphs, a vertex subset $S$ of a digraph $D$ is in general position if for any pair $u,v$ in $S$ no vertex of $S \setminus \{ u,v\} $ lies on a shortest $u,v$-path. We note that in a digraph the shortest $u,v$- and $v,u$-paths must be considered separately, and these may have different length. Also, by the definition of~\cite{digraphs}, it is allowed that a general position set contains vertices $u,v$ such that $D$ has no directed path from $u$ to $v$. An example from~\cite{digraphs} is shown in Figure~\ref{fig:d2k3n20}.

\begin{figure}
	\centering
	\begin{tikzpicture}[x=0.4mm,y=-0.4mm,inner sep=0.2mm,scale=0.25,very thick,vertex/.style={circle,draw,minimum size=10,fill=lightgray}]
	\node at (181.4,173.2) [vertex] (v1) {};
	\node at (222.6,300.4) [vertex,color=red] (v2) {};
	\node at (317.2,211.5) [vertex,color=red] (v3) {};
	\node at (95.9,379.5) [vertex] (v4) {};
	\node at (679.4,379.5) [vertex] (v5) {};
	\node at (387.6,87.8) [vertex] (v6) {};
	\node at (387.6,671.2) [vertex] (v7) {};
	\node at (118.1,267.9) [vertex] (v8) {};
	\node at (118.1,491.1) [vertex,color=red] (v9) {};
	\node at (657.2,267.9) [vertex,color=red] (v10) {};
	\node at (657.2,491.1) [vertex] (v11) {};
	\node at (499.3,110) [vertex,color=red] (v12) {};
	\node at (276,110) [vertex] (v13) {};
	\node at (499.3,649) [vertex] (v14) {};
	\node at (276,649) [vertex,color=red] (v15) {};
	\node at (181.4,585.8) [vertex] (v16) {};
	\node at (593.9,173.2) [vertex] (v17) {};
	\node at (593.9,585.8) [vertex] (v18) {};
	\node at (547.6,297.7) [vertex] (v19) {};
	\node at (456.8,212.6) [vertex] (v20) {};
	\path
	(v1) edge[->] (v2)
	(v1) edge[->] (v3)
	(v8) edge[->] (v1)
	(v13) edge[->] (v1)
	(v2) edge[->] (v4)
	(v2) edge[->] (v5)
	(v18) edge[->] (v2)
	(v3) edge[->] (v6)
	(v3) edge[->] (v7)
	(v18) edge[->] (v3)
	(v4) edge[->] (v8)
	(v4) edge[->] (v9)
	(v19) edge[->] (v4)
	(v5) edge[->] (v10)
	(v5) edge[->] (v11)
	(v19) edge[->] (v5)
	(v6) edge[->] (v12)
	(v6) edge[->] (v13)
	(v20) edge[->] (v6)
	(v7) edge[->] (v14)
	(v7) edge[->] (v15)
	(v20) edge[->] (v7)
	(v8) edge[->] (v12)
	(v15) edge[->] (v8)
	(v13) edge[->] (v9)
	(v9) edge[->] (v14)
	(v9) edge[->] (v16)
	(v10) edge[->] (v13)
	(v14) edge[->] (v10)
	(v10) edge[->] (v17)
	(v12) edge[->] (v11)
	(v11) edge[->] (v15)
	(v11) edge[->] (v18)
	(v12) edge[->] (v17)
	(v14) edge[->] (v18)
	(v15) edge[->] (v16)
	(v16) edge[->] (v19)
	(v16) edge[->] (v20)
	(v17) edge[->] (v19)
	(v17) edge[->] (v20)
	;
	\end{tikzpicture}
	\caption{A $\gp $-set in a digraph (vertices in red)}
	\label{fig:d2k3n20}
\end{figure}

Some elementary bounds and structural information is given in~\cite{digraphs}. For example, the general position number of a digraph is equal to that of its converse. For any digraph, the set of sources and sinks is always a general position set. More generally,~\cite{digraphs} defines a \emph{near-sink} to be a vertex that is not a sink, but all of its out-neighbours are sinks, and likewise a \emph{near-source} is a vertex that is not a source, but all of its in-neighbours are sources.

\begin{theorem}{\rm \cite[Lemma~2.8]{digraphs}}\label{thm:sourcesandsinks}
Let $\sink (D)$, $\nsink (D)$, $\source (D)$ and $\nsource (D)$ be the sets of sinks, near-sinks, sources, and near-sources of a digraph $D$, respectively. Then the sets $\sink (D) \cup \source (D)$, $\nsink (D) \cup \sink (D)$ and $\nsource (D) \cup \source (D)$ are in general position.
\end{theorem}

Of note is the following extension of the characterisation of general position sets in Theorem~\ref{thm:gpsets-characterised} to directed graphs. This theorem implies that a general position set in an oriented graph induces an acyclic subdigraph.

\begin{theorem} {\rm \cite[Theorem~2.16]{digraphs}}\label{thm:characterisation_digraph}
   Let $S$ be a vertex subset of the digraph $D$ and denote the strong components of the subdigraph induced by $S$ by $S_1,\ldots,S_r$. Then $S$ is a general position set of $D$ if and only if 
    \begin{itemize}
        \item each strong component $S_i$ is a complete digraph, 
        \item $\{ S_1,\dots ,S_r\} $ is \emph{distance-constant}, i.e. for any $1 \leq i < j \leq r$, we have $d(u,v) = d(u',v')$ for any $u,u' \in S_i$ and $v,v' \in S_j$ (where we allow this distance to be infinite),
        \item for no $1 \leq i,j,k \leq r$ it is the case that \[ d(S_i,S_k) = d(S_i,S_j)+d(S_j,S_k) < \infty .\]
    \end{itemize} 
\end{theorem}

It is trivial that for any digraph with order $n \geq 2$ we have $2 \leq \gp (G) \leq n$. One feature of general position in digraphs that contrasts markedly with the undirected case is that there are many graphs meeting the trivial upper and lower bounds. It follows from Theorem~\ref{thm:characterisation_digraph} that a digraph $D$ with order $n$ satisfies $\gp (D) = n$ if and only if $D$ is transitive. A sufficient condition for $\gp (D) = 2$ is that $D$ has a Hamilton path, although a family of non-traceable digraphs with general position number two is constructed.

By replacing each edge by a pair of oppositely directed arcs, the general position problem for undirected graphs can be construed as a special case of the general position problem for digraphs. As such, Theorem~\ref{thm:fundamental_complexity_result} implies immediately that the decision version of the general position problem for digraphs is NP-complete. More interestingly, the complexity of the general position problem is analysed in~\cite{digraphs} within the class of oriented graphs via a reduction from the independent set problem. 

\begin{theorem}\label{thm:NP}
    The decision version of the general position problem is NP-complete even for oriented graphs $D$ with diameter three.
\end{theorem}

The general position problem has been investigated in four common families of digraphs, namely circulant digraphs, oriented trees, permutation digraphs and Kautz digraphs. The circulant digraph $\Circ(n, [d])$ is the digraph with vertex set $\mathbb{Z}_n$ and arcs $x \rightarrow x+i$ for all $x \in \mathbb{Z}_n$ and $i \in [d]$. The following bound is conjectured to be best possible.

\begin{theorem} {\rm \cite[Theorems~3.2 \& 3.3]{digraphs}} \newline
	Set $n = qd+r$, where $q \geq 0$ and $2 \leq r \leq d+1$. The general position number of $\Circ(n, [d])$ satisfies
	\[r\leq \gp(\Circ(n,[d]))\leq d+1.\] When $n\equiv 1\pmod{d}$, we have $\gp(\Circ(n,[d])) = d+1$.

    If $n \equiv a \pmod{d}$ where $a\geq 2$, then $\gp(\Circ(n,[d])) \geq \left \lceil \frac{d+1}{a} \right \rceil $.
\end{theorem}

Finding the general position numbers of oriented trees is a surprisingly difficult problem. It is shown in~\cite{digraphs} that the set of sinks and near-sinks from Theorem~\ref{thm:sourcesandsinks} is optimal for out-arborescences (and sources and near-sources for in-arborescences). However, Theorem~\ref{thm:sourcesandsinks} is in general not best possible, and the problem remains open.

The permutation digraph $\Pe (d,k)$ is defined on the alphabet $[0,\ldots ,d+k-1]$. The vertices of $\Pe (d,k)$ are words $x_0x_1\ldots x_{k-1}$ of length $k$ of distinct letters from $[0,\ldots ,d+k-1]$ and the arcs have the form $x_0x_1\ldots x_{k-1} \rightarrow x_1\ldots x_{k-1}x_k$, where $x_k$ is distinct from the letters $x_0,x_1,\ldots ,x_{k-1}$. The permutation digraphs have some applications in extremal digraph theory. The general position numbers are known exactly for $k = 2$ and a lower bound is known for larger $k$ that is conjectured to be the correct answer.

\begin{theorem} {\rm \cite[Theorems~3.5 \& 3.6]{digraphs}}
For $d \geq 2$, $\gp (\Pe(d,2)) = \left \lceil \frac{(d+2)^2}{4} \right \rceil $. For $k \geq 3$ and $d \geq 2k$, $\gp (\Pe (d,k)) \geq 2(d+k-2)_{k-1}$.
\end{theorem}

The vertices of the Kautz digraph $\Ka (m,k)$ are the words $x_0x_1\ldots x_{k-1}$ of length $k$ from an alphabet of size $m$ such that any pair of consecutive letters are distinct. The arcs are of the form $x_0x_1\ldots x_{k-1} \rightarrow x_1\ldots x_{k-1}x_k$, where $x_{k-1} \neq x_k$. Again, the exact value is known for $k = 2$. However, the value for larger $k$ is not so clear.

\begin{theorem} {\rm \cite[Theorem~3.4]{digraphs}}
    For $m \geq 3$, $\gp (\Ka(m,2)) = \left \lfloor \frac{m^2}{3} \right \rfloor $.
\end{theorem}

Several interesting problems are associated with the orientations of an undirected graph. The \emph{general position spectrum} $S_{\gp}(G)$ of an undirected graph $G$ is the set of general position numbers of all orientations of $G$. The largest entry of $S_{\gp}(G)$ is the \emph{upper oriented general position number} $\ugp(G)$ and the smallest entry is the \emph{lower oriented general position number} $\lgp(G)$. A graph $G$ with order $n$ is \emph{$\gp $-full} if $S_{\gp }(G) = [2,n]$. Numerical evidence suggests two natural conjectures.

\begin{conjecture} {\rm \cite[Conjecture~4.3]{digraphs}}\label{conj:spectrum1}
    The general position spectrum of any graph $G$ with order $n \geq 3$ contains at least two distinct entries, that is, $\lgp (G) < \ugp (G)$.
\end{conjecture}

\begin{conjecture} {\rm \cite[Conjecture~4.30]{digraphs}}
    The general position spectrum of any graph is an interval.
\end{conjecture}

Conjecture~\ref{conj:spectrum1} is proven for comparability graphs (recall that a comparability graph is a graph with a transitive orientation), triangle-free sub-cubic graphs, regular graphs with degree at most five, and graphs with large size. The $\gp $-spectrum of cycles and some families of trees are determined. For complete bipartite graphs, the balanced complete bipartite graph $K_{n,n}$ is $\gp $-full, whereas for $r \geq 4$ the $\gp $-spectrum is $K_{r,2}$ is $[\left \lceil \frac{r}{2} \right \rceil ,r+2]$. It would be interesting to determine the $\gp $-spectrum of the remaining complete bipartite graphs.

By Theorem~\ref{thm:characterisation_digraph}, the upper general position number satisfies $\ugp (G) = n$ if and only if $G$ is a comparability graph. At the other extreme, it was conjectured that it is NP-hard to determine if an undirected graph has an orientation with general position number two. 

\section{Open problems}\label{sec:openproblems}
We conclude this survey with a list of what the authors consider to be some of the most important open problems in the field of position problems. When applicable, we also refer to the paper in which the problem was first posed.

\begin{itemize}
\item Find a value of $n$ for which $2n$ pawns cannot be placed in general position in the No-Three-In-Line Problem, or show that $2n$ pawns is always possible (perhaps for large $n$). This remains a challenging problem.
\item \cite{AraujoPardo,Woodcolouring} What is the smallest value of $c$ such that $(c+\epsilon )n$ colours suffice for any $\epsilon > 0$ and sufficiently large $n$ to cover an $n \times n$ chessboard with coloured pawns, such that each colour class has the no-three-in-line property? (It is now known that $\frac{1}{2} \leq c \leq 1$, see~\cite{AraujoPardo}).
\item \cite{AicEppHai} Are the orders of both types of `geometric dominating sets' from Subsection~\ref{subsec:lower} monotone in $n$ for the $n \times n$ grid? And are smallest geometric dominating sets with the no-three-in-line property always larger than the smallest geometric dominating sets without this property?
\item Investigate position problems in signed graphs and hypergraphs.
\item \cite[Problem~4.8]{Klavzar-2019} For which graphs $G$ and $H$ does $\gp(G\boxtimes H) = \gp(G)\gp(H)$ hold? (It was shown in~\cite{Vesel} that the conjecture that this holds for all pairs of graphs is false.) 
\item \cite{ThaChaTuiThoSteErs-2024} Is the ratio $\frac{\vp (G)}{\vp^- (G)} $ bounded for connected graphs?
\item What is the largest possible size of a graph with given order and equidistant number?
\item \cite[Conjecture 9.3]{ChaDiSreeThoTui2024+} Is it true that if $G$ is a graph and $\pi \in \{\gp, \mono\}$, then $\chi _{\pi_{\rm i}}(G) +\chi _{\pi_{\rm i}}(\overline{G}) \leq n(G)+1$? 
\item Prove that the detour position number of a graph is bounded above by $n-D+1$, where $D$ is the detour diameter.
\item \cite[Conjecture 3.3]{Kruft-2024} Does every graph contain a terminal set?
\item \cite[Conjecture 2.10]{Kruft-2024} Is it true that $\gp ^-(G \cp H) \geq \min \{ \gp ^-(G),\gp ^-(H)\} $ for any graphs $G,H$?
\item \cite{TuiThoCha-2025} Find an exact formula for $\mex (n;k)$ for $k \geq 3$, and improve the Ramsey-theoretic upper bound for $\gex (n;k)$.
\item \cite[Conjecture~4.5]{TuiThoCha-2025} Prove that for any $n \geq 11$ there is a circulant graph with order $n$, monophonic position number two and diameter two.
\item Identify families of graphs with unimodal general position polynomial.
\item \cite[Problem~9.2]{ChaDiSreeThoTui2024+} (Packing problem) For a connected graph $G$, what is the largest number of disjoint maximal general position sets contained in $G$?
\item \cite{TianKlavzar-2025} Characterise the graphs $G$ with $\gp(G) = \gpo(G) = \gpd(G) = \gpt(G)$. 
\item Which graphs contain a dominating set that is in general position? Or a metric resolving set that is in general position?
\item \cite[Conjecture 11]{KlavzarTianTuite-2026} What is the complexity of the decision version of the Builder/Blocker general position game?
\item \cite{KlaKriTuiYer-2023} What is the complexity of finding $\mob (G)$?
\item \cite{Fitzpatrick-2006} Find the fractional general position number of Cartesian grids.
\item Is there a graph with general position number three and mutual visibility number greater than seven? (It is not hard to show that there is a graph $G$ with $\gp (G) = a$ and $\mu (G) = b$ if $4 \leq a \leq b$, but what happens for $a = 3$, $b> 7$ is unknown).
\item \cite{digraphs} Find a formula for the general position number of oriented trees.
\item \cite{digraphs} Determine the general position number of the permutation digraphs $\Pe (d,k)$ and Kautz digraphs $\Ka (d,k)$ for $k \geq 3$.
\item \cite{digraphs} Find the general position spectrum of all complete bipartite (or even multipartite) graphs.
\item \cite{digraphs} Prove that the general position spectrum of an undirected graph with order $n \geq 3$ is an interval and contains at least two distinct entries. 
\item \cite{digraphs} Investigate the complexity of determining the lower oriented general position number. In particular, is it NP-hard to decide if an undirected graph has an orientation with general position number two?
\end{itemize}

\section*{Acknowledgements}

The three authors thank the anonymous referees for their valuable feedback. Sandi Klav\v{z}ar was supported by the Slovenian Research Agency (ARIS) under the grants P1-0297, N1-0285, N1-0355, N1-0431, J1-70045.
The research of James Tuite was partly funded by LMS Research in Pairs grant 42235. He also gratefully acknowledges the hospitality of the Institute of Mathematics, Physics and Mechanics, Ljubljana, the University of Ljubljana, and travel funding from the OU Crowther Fund. The authors thank Grahame Erskine for helpful discussions, in particular for finding the $\gp $-sets of the Clebsch and Frucht graphs, and Thomas Prellberg for introducing them to new developments in the geometric No-Three-In-Line Problem.



\begin{thebibliography}{999}
\bibliographystyle{plain}

\bibitem{Abiad2024} A.~Abiad, A.J.~Ameli, L.~Reijnders, 
The clique number of the exact distance $t$-power graph: Complexity and eigenvalue bounds, 
Discrete Appl.\ Math.\ 363 (2025) 55--70.
		
\bibitem{AdeHolKel}
M.A.~Adena, D.A.~Holton, P.A.~Kelly,
Some thoughts on the no-three-in-line problem,
Lecture Notes in Math.\ 403 (1974) 6--17.

\bibitem{Adhikary2019} 
R.~Adhikary, K.~Bose, M.K.~Kundu, B.~Sau, 
Mutual visibility by asynchronous robots on infinite grid, 
in Algorithms for Sensor Systems: 14th International Symposium on Algorithms and Experiments for Wireless Sensor Networks, Springer International Publishing (2018) 83--101.
		
\bibitem{AicEppHai}
O.~Aichholzer, D.~Eppstein, E.M.~Hainzl,
Geometric dominating sets-a minimum version of the No-Three-In-Line Problem,
Comput.\ Geom.\ 108 (2023) Paper 101913.

\bibitem{connectfour} 
L.V.~Allis, 
A knowledge-based approach of Connect-Four, 
J.\ Int.\ Comput.\ Games Assoc.\ 11 (1988) 165--165.

\bibitem{Go} 
L.V.~Allis, 
Searching for Solutions in Games and Artificial Intelligence. 
Ph.D.\ thesis, University of Limburg, The Netherlands (1994) 121--154.
  
\bibitem{AnaChaChaKlaTho}
B.S.~Anand, U.~Chandran S.V., M.~Changat, S.~Klav\v{z}ar, E.J.~Thomas,
Characterization of general position sets and its applications to cographs and bipartite graphs,
Appl.\ Math.\ Comput.\ 359 (2019) 84--89.

\bibitem{Araujo-2018} 
J.~Araujo, G.~Ducoffe, N.~Nisse, and K.~Suchan, On interval number in cycle convexity, Discrete Math.\ Theor.\ Comput.\ Sci. 20(1) (2018) Paper 13. 

\bibitem{Araujo-2025}
J.~Araujo, M.C.~Dourado, F.~Protti, R.~Sampaio, The iteration time and the general position number in graph convexities, 
Appl.\ Math.\ Comput.\ 487 (2025) Paper 129084. 

\bibitem{Araujo-2026} 
J.~Araujo, S.N.~Araujo, P.P.~Medeiros, N.~Nisse, C.~Silva,
On the rank and the general position number in cycle convexity,
Lecture Notes in Comput.\ Sci.\ 16587 (2026) 73--87. 

\bibitem{AraujoPardo}
G.~Araujo-Pardo, Gabriela, L.~Mart\'{i}nez-Sandoval,
The arc chromatic number for Galois projective planes, affine planes and Euclidean grids,
\url{arXiv:2601.19043} (2026).

\bibitem{BaloghLuo} 
J.~Balogh, H.~Luo, 
Maximum number of points in general position in a random subset of finite 3-dimensional spaces, 
Comput.\ Geom.\ (2026) Paper 102244.

\bibitem{Balogh} 
J.~Balogh, J.~Solymosi, 
On the number of points in general position in the plane, 
Discrete Anal.\ 16 (2018) 1--20.

\bibitem{BangJensenGutin}
J.~Bang-Jensen, G.~Gutin,
Digraphs: Theory, Algorithms and Applications, 
Springer, London, 2001. 

\bibitem{Bhagat2020} 
S.~Bhagat, 
Optimum algorithm for the mutual visibility problem, 
Lecture Notes in Comput.\ Sci.\ 12049 (2020) 31--42. 

\bibitem{Bhagat2019} 
S.~Bhagat, K.~Mukhopadhyaya, 
Mutual visibility by robots with persistent memory, 
Lecture Notes in Comput.\ Sci.\ 11458 (2019) 144--155. 

\bibitem{Bhagat2017} 
S.~Bhagat, K.~Mukhopadhyay, 
Optimum algorithm for mutual visibility among asynchronous robots with lights, 
Lecture Notes in Comput.\ Sci.\ 10616 (2017) 341--355.

\bibitem{Boshar} 
S.~Boshar, 
personal communication (2022).

\bibitem{Brassetal} 
P.~Brass, E.~Cenek, C.A.~Duncan, A.~Efrat, C.~Erten, D.P.~Ismailescu, S.G.~Kobourov, A.~Lubiw, J.S.B.~Mitchell,  
On simultaneous planar graph embeddings, 
Comput.\ Geom.\ 36 (2007) 117--130.

\bibitem{Brass} 
P.~Brass, W.~Moser, J.~Pach, 
Packing lattice points in subspaces. 
Section 10.1.\ in Research Problems in Discrete Geometry. Springer, New York (2005) 417--421.

\bibitem{bresar-2010}
B.~Bre{\v{s}}ar, S.~Klav\v{z}ar, D.F.~Rall,
Domination game and an imagination strategy,
SIAM J.\ Discrete Math.\ 24 (2010) 979--991.
	
\bibitem{book-2021}
B.~Bre\v{s}ar, M.A.~Henning, S.~Klav\v{z}ar, D.F.~Rall, 
Domination Games Played on Graphs, 
Springer Briefs in Mathematics, 2021.  

\bibitem{Cao-2026}
Y.~Cao, S.~Ji,
A note on the edge general position number of cactus graphs,
Open Math.\ 23 (2025) 20250207.

\bibitem{cao-2025}
Y.~Cao, S.~Ji, L.~Wang,
The edge general position number of some graphs,
Bull.\ Malays.\ Math.\ Sci.\ Soc.\ 48 (2025) Paper 95.

\bibitem{digraphs}
U.~Chandran S.V., G.~Di Stefano, G.~Erskine, H.~Sreelatha, E.J.~Thomas, J.~Tuite, 
The general position number of digraphs,
\url{arXiv:2604.15909} (2026).

\bibitem{ChaDiSreeThoTui2024+} 
U.~Chandran S.V., G.~Di Stefano, H.~Sreelatha, E.J.~Thomas, J.~Tuite, 
Colouring a graph with position sets, 
Ars.\ Math.\ Contemp.\ (2025) \url{doi.org/10.26493/1855-3974.3454.a3c}. 

\bibitem{KlaSam}
U.~Chandran S.V., S.~Klav\v{z}ar, P.K.~Neethu, R.~Sampaio,
The general position avoidance game and hardness of general position games,
Theoret.\ Comput.\ Sci.\ 988 (2024) Paper 114370.

\bibitem{ChandranKlavzar2025}
U.~Chandran S.V., S.~Klav\v{z}ar, P.K.~Neethu, J.~Tuite,
Monophonic position sets of Cartesian and lexicographic products of graphs, 
\url{arXiv:2412.09837v2} (2025).

\bibitem{ullas-2016}
U.~Chandran S.V., G.J.~Parthasarathy,
The geodesic irredundant sets in graphs,
Int.\ J.\ Math.\ Combin.\  4 (2016) 135--143.

\bibitem{Graphsdigraphs} 
G.~Chartrand, L.~Lesniak, P.~Zhang, 
Graphs \& Digraphs, Chapman \& Hall, London (1996).

\bibitem{Chartrand-1989}
G.~Chartrand, O.R.~Oellermann, S.~Tian, H.B.~Zou,
Steiner distance in graphs,
\v{C}asopis P\v{e}st.\ Mat.\ 114 (1989) 399--410.

\bibitem{Chartrand-2004} 
G.~Chartrand, P.~Zhang, T.W.~Haynes, 
Distance in graphs---taking the long view, 
AKCE Int.\ J.\ Graphs Comb.\ 1 (2004) 1--13.

\bibitem{ChenChvatal}
X.~Chen, V. Chv\'{a}tal, 
Problems related to a De Bruijn-Erd\H{o}s theorem, 
Discrete Appl.\ Math.\ 156 (2008) 2101--2108.

\bibitem{Chen} 
Y.~Chen, X.~Liu, J.~Nie, J.~Zeng, 
Random Tur\'an and counting results for general position sets over finite fields, 
Sci.\ China Math.\ 26 (2025) 3043--3062. 

\bibitem{Chen2}
Y.~Chen, J.~Nie, J.~Yu, W.~Zhang, 
Maximum in-general-position set in a random subset of $\mathbb {F}^ d_q$, 
J.\ Comb.\ Theory Ser.\ B 179 (2026) 335--366. 

\bibitem{chen-2026}
H.~Chen, Y.~Yuan, 
Coefficients of the matching polynomial of a self-complementary graph,
Graphs Combin.\ 42 (2026) Paper 22. 

\bibitem{CiDiDiNav2023} 
S.~Cicerone, A.~Di Fonso, G.~Di Stefano, A.~Navarra, 
The geodesic mutual visibility problem: Oblivious robots on grids and trees, 
Pervasive Mobile Comput.\ 95 (2023) Paper 101842.

\bibitem{Cody-2026}
B.~Cody, R.~Detore, 
Metric general position extensions of classical graph invariants,
\url{arXiv:2601.04351} (2026).

\bibitem{Cody-2025}
B.~Cody, G.~Moore,
The $k$-general $d$-position problem for graphs,
Discrete Appl.\ Math.\ 366 (2025) 135--151.

\bibitem{CooPikSchWar}
A.S.~Cooper, O.~Pikhurko, J.R.~Schmitt, G.S.~Warrington,
Martin Gardner's minimum no-3-in-a-line problem,
Amer.\ Math.\ Monthly 121 (2014) 213--221.

\bibitem{Cooper} 
J.N.~Cooper, J.~Solymosi, 
Collinear points in permutations, 
Ann.\ Comb.\ 9 (2005) 169--175.

\bibitem{DelbinPrema-2018} 
S.~Delbin Prema, C.~Jayasekaran, 
The detour irredundant number of a graph, 
Int.\ J.\ Pure Appl.\ Math.\ 119 (2018) 53--62.

\bibitem{Wolfram} \url{https://demonstrations.wolfram.com/NoThreeInLineProblem/} (accessed 16/09/2024).

\bibitem{DiLuna-2017}
G.A.~Di Luna, P.~Flocchini, S.G.~Chaudhuri, F.~Poloni, N.~Santoro, G.~Viglietta,
Mutual visibility by luminous robots without collisions,
Inf.\ Comput.\ 254 (2017) 392--418.

\bibitem{DiLuna2014} 
G.A.~Di Luna, P.~Flocchini, F.~Poloni, N.~Santoro, G.~Viglietta, 
The mutual visibility problem for oblivious robots, 
in: CCCG 2014, Halifax, Nova Scotia, August 11–13 (2014) Paper 51.  

\bibitem{DiStefano-2022}
G.~{Di Stefano},
Mutual visibility in graphs,
Appl.\ Math. Comput.\ 419 (2022) Paper 126850.

\bibitem{DiKlKrTu-2024}
G.~{Di Stefano}, S.~Klav\v{z}ar, A.~Krishnakumar, J.~Tuite, I.G.~Yero,
Lower general position sets in graphs,
Discuss.\ Math.\ Graph Theory 45 (2025) 509--531.

\bibitem{Dokyeesun-2025}
P.~Dokyeesun, S.~Klav\v zar, J.~Tian,
The general position number under vertex and edge removal,
Quaest.\ Math.\ (2025) 48 (2025) 1277--1290.

\bibitem{Dokyeesun-2026}
P.~Dokyeesun, S.~Klav\v zar, D.~Kuziak, J.~Tian,
General position problems in strong and lexicographic products of graphs,
Comput.\ Appl.\ Math.\ 45 (2026) Paper 97. 

\bibitem{DongXu} 
Z.~Dong, Z.~Xu, 
Large grid subsets without many cospherical points,
\url{arXiv:2506.18113} (2025).

\bibitem{downey-2013}
R.G.~Downey, M.R.~Fellows, 
Fundamentals of Parameterized Complexity,
Springer, London, 2013. 
       
\bibitem{dudeney-1917}
H.E.~Dudeney,
Amusements in Mathematics,
Nelson, Edinburgh (1917).

\bibitem{egecyoglu-2023}
\"O.~E\u gecio\u glu, S.~Klav\v{z}ar, M.~Mollard,
Fibonacci Cubes---With Applications and Variations,
World Scientific, Singapore, 2023. 

\bibitem{Erdos} 
P.~Erd\H{o}s, 
On some metric and combinatorial geometric problems, 
Discrete Math.\ 60 (1986) 147--153.

\bibitem{Evans} 
A.~Evans, E.~Shallcross, J.~Tuite, 
Position games in graphs, 
manuscript (2025).

\bibitem{fialho-2026}
P.M.S.~Fialho, E.~Juliano, A.~Procacci, 
On the zero-free region for the chromatic polynomial of graphs with maximum degree $\Delta$ and girth $g$, 
Discrete Math.\ 349 (2026) Paper 114825.

\bibitem{Fitzpatrick-2006} 
S.L.~Fitzpatrick, J.~Janssen, R.J.~Nowakowski,
Fractional isometric path number,
Australas.\ J.\ Comb.\ 34 (2006) 281--298.

\bibitem{Flammenkamp1} 
A.~Flammenkamp,
Progress in the no-three-in-line-problem,
J.\ Comb.\ Theory Ser.\ A 60 (1992) 305--311.

\bibitem{Flammenkamp2} 
A.~Flammenkamp,
Progress in the no-three-in-line problem, II,
J.\ Comb.\ Theory Ser.\ A 81 (1998) 108--113.

\bibitem{Flammenkampwebsite}
A.~Flammenkamp,
\url{https://wwwhomes.uni-bielefeld.de/~achim/no3in/readme.html}
(accessed 22/07/2026).

\bibitem{foucaud_cliques_2021} 
F.~Foucaud, S.~Mishra, N.~Narayanan, R.~Naserasr, P.~Valicov, 
Cliques in exact distance powers of graphs of given maximum degree, 
Procedia Comput.\ Sci.\ 195 (2021) 427--436. 

\bibitem{Froese} 
V.~Froese, I.~Kanj, A.~Nichterlein, R.~Niedermeier, 
Finding points in general position, 
Int. J. Comput. Geom. Appl. 27 (4) (2017) 277--296.

\bibitem{Gabriel}
T.~G\'{a}briel, M.~J\'{a}nosik, D.~Melj\'{a}n, B.~N\'{a}dor,
The extensible no-$(k(n)+1)$-in-line problem,
\url{arXiv:2606.02843} (2026).

\bibitem{Gardner-1976}
M.~Gardner,
Mathematical games: Combinatorial problems, some old, some new and all newly attacked by computer,
Sci.\ Amer.\ 235 (1976) 131--137.

\bibitem{Gasser} 
R.~Gasser, 
Solving nine men's morris, 
Comput.\ Intell.\ 12 (1996) 24--41.

\bibitem{Ghorbani-2021}
M.~Ghorbani, H.R.~Maimani, M.~Momeni, F.R.~Mahid, S.~Klav\v{z}ar, G.~Rus,
The general position problem on Kneser graphs and on some graph operations,
Discuss.\ Math.\ Graph Theory 41 (2021) 1199--1213.

\bibitem{Ghosal}
A.~Ghosal, 
A note on the extensible no-three-in-line problem,
\url{arXiv:2605.07000} (2026).

\bibitem{Ghosal2026}
A.~Ghosal, R.~Goenka, A.~Grebennikov, P.~Keevash, M.~Kwan, H.T.~Pham,  
No $(k+1)$-in-line problem for $k \geq 3$,
\url{arxiv:2607.05255} (2026).

\bibitem{GhosalGoenkaKeevash}
A.~Ghosal, R.~Goenka, P.~Keevash, 
On subsets of lattice cubes avoiding affine and spherical degeneracies,
Discrete Comput.\ Geom.\ (2026) \url{https://doi.org/10.1007/s00454-026-00853-7}.

\bibitem{Kaplan} 
\url{https://11011110.github.io/blog/2018/11/10/random-no-three.html} (accessed 16/09/2024).

\bibitem{GrebennikovKwan}
A.~Grebennikov, M.~Kwan,
No-$(k+1)$-in-line problem for large constant $k$,
\url{arXiv:2510.17743} (2025).

\bibitem{Green}
B.~Green,
100 open problems,
\url{https://people.maths.ox.ac.uk/greenbj/papers/open-problems.pdf} (2024).

\bibitem{Guy} 
R.K.~Guy, 
The no-three-in-a-line problem, 
in Unsolved Problems in Number Theory, 2nd ed. 
Springer, New York (1994) 240--244.
		
\bibitem{Guy-1968} 
R.K.~Guy, P.A.~Kelly,
The no-three-in-line problem,
Can.\ Math.\ Bull.\ 11 (1968) 527--531.
		
\bibitem{Hall-1975} 
R.R.~Hall, T.H.~Jackson, A.~Sudbery, K.~Wild,
Some advances in the no-three-in-line problem,
J.\ Comb.\ Theory Ser.\ A. 18 (1975) 336--341.

\bibitem{Hamed-2026b}
Z.~Hamed-Labbafian, M.A.~Henning, M.~Tavakoli,
Edge general position in graphs: Graph products, integer linear programming and some applications,
Discrete Appl.\ Math.\ 386 (2026) 1--8. 

\bibitem{Hamed-2026}
Z.~Hamed-Labbafian, N.~Sabeghi, M.~Tavakoli, S.~Klav\v{z}ar,
Three algorithmic approaches to the general position problem,
Bull.\ Aust.\ Math.\ Soc.\ 113 (2026) 1--9.  

\bibitem{Haponenko-2024}
V.~Haponenko, S.~Kozerenko, 
All-path convexity: Two characterizations, general position number, and one algorithm,
Discrete Math.\ Lett.\ 13 (2024) 58--65. 

\bibitem{haritha-2026}
Haritha S., U.~Chandran S.V., 
General position and mutual-visibility in shadow graphs, 
\url{arXiv:2601.19769} (2026).

\bibitem{haynes-2007}
T.W.~Haynes, M.A.~Henning, P.J.~Slater, L.C.~van der Merwe, 
The complementary product of two graphs,
Bull.\ Inst.\ Combin.\ Appl.\ 51 (2007) 21--30.

\bibitem{positionalgames} 
D.~Hefetz, M.~Krivelevich, M.~Stojakovi\'{c}, T.~Szab\'{o}, 
Positional Games, Birkhäuser, Basel, 2014.

\bibitem{hinz-2017}
A.M.~Hinz, S.~Klav\v{z}ar, S.S.~Zemlji\v{c}, 
A survey and classification of {S}ierpi\'nski-type graphs,
Discrete Appl.\ Math.\ 217 (2017) 565--600. 

\bibitem{irsic-2024}
V.~Ir\v si\v c, S.~Klav\v zar, G.~Rus, J.~Tuite, 
General position polynomials,
Results Math.\ 79 (2024) Paper 110.

\bibitem{fractional}
V.~Ir\v si\v c Chenoweth, G.~Erskine, S.~Klav\v zar, G.~Rus, J.~Tian, J.~Tuite,
Fractional position problems on graphs,
in preparation (2026).

\bibitem{Janosik}
M.~J\'{a}nosik, A.~N\'{a}dor, N.Z.~L\'{o}r\'{a}nt, L.B.~Simon,
Avoiding configurations of small size in the square grid,
\url{arXiv:2601.14465} (2026).

\bibitem{kadrawi-2025}
O.~Kadrawi, V.~Levit, 
The independence polynomial of trees is not always log-concave starting from order 26,
Ars Math.\ Contemp.\ 25(4) (2025) Paper 3. 

\bibitem{KlaKriKuzShaTuiYer-2026} 
S.~Klav\v zar, A.~Krishnakumar, D.~Kuziak, E.~Shallcross, J.~Tuite, I.G.~Yero, 
Moving through Cartesian products, coronas and joins in general position,
Discrete Appl.\ Math.\ 379 (2026) 768--780.

\bibitem{KlaKriTuiYer-2023} 
S.~Klav\v zar, A.~Krishnakumar, J.~Tuite, I.G.~Yero, 
Traversing a graph in general position, 
Bull.\ Aust.\ Math.\ Soc.\ (2023) 1--13.        

\bibitem{KlaKuzPetYer-2021}
S.~Klav\v zar, D.~Kuziak, I.~Peterin, I.G.~Yero,
A Steiner general position problem in graph theory,
Comput.\ Appl.\ Math.\ 40 (2021) 1--15.

\bibitem{KlaLakRoy-2025}
S.~Klav\v{z}ar, A.S.~Lakshmanan, D.~Roy, 
Counting largest mutual-visibility and general position sets of glued $t$-ary trees,
Results Math. 80 (2025) Paper 207.

\bibitem{KlaNeeCha-2022}
S.~Klav\v{z}ar, P.K.~Neethu, U.~Chandran S.V.,
The general position achievement game played on graphs,
Discrete Appl.\ Math.\ 317 (2022) 109--116.

\bibitem{KlaPatRusYero-2021}
S.~Klav\v{z}ar, B.~Patk\'{o}s, G.~Rus, I.G.~Yero, 
On general position sets in {C}artesian products,
Results Math.\ 76 (2021) Paper 123.

\bibitem{KlaRalYer-2021}
S.~Klav\v zar, D.F.~Rall, I.G.~Yero,
General $d$-position sets.
Ars Math.\ Contemp.\ 21 (2021) Paper P1.03.

\bibitem{KlavzarRus-2021}
S.~Klav\v zar, G.~Rus,
The general position number of integer lattices, 
Appl.\ Math.\ Comput.\ 390 (2021) Paper 125664.

\bibitem{KlavzarTan-2023}
S.~Klav\v zar, E.~Tan,
Edge general position sets in Fibonacci and Lucas cubes,
Bull.\ Malays.\ Math.\ Sci.\ Soc.\ 46 (2023) Paper 120.

\bibitem{KlavzarTianTuite-2026}
S.~Klav\v zar, J.~Tian, J.~Tuite,
Builder-Blocker general position games,
Commun.\ Comb.\ Optim.\ 11 (2026)  11 465--486. 

\bibitem{Klavzar-2019}
S.~Klav\v{z}ar, I.G.~Yero,
The general position problem and strong resolving graphs,
Open Math.\ 17 (2019) 1126--1135.

\bibitem{Korner-1995}
J.~K\"orner,
On the extremal combinatorics of the Hamming space,
J.\ Comb.\ Theory Ser.\ A 71 (1995) 112--126.

\bibitem{KorzeVesel-2023}
D.~Kor\v{z}e, A.~Vesel, 
General position sets in two families of Cartesian product graphs,
Mediterr.\ J.\ Math.\ 20 (2023) Paper 203.

\bibitem{KorzeVesel-2024}
D.~Kor\v{z}e, A.~Vesel, 
Mutual-visibility and general position sets in Sierpi\'nski triangle graphs,
Bull.\ Malays.\ Math.\ Sci.\ Soc.\ 48 (2025) Paper 106. 

\bibitem{Kovacs-2025 algebraic} 
~Kov\'{a}cs, Z.L.~Nagy, D.R.~Szab\'{o},
Randomised algebraic constructions for the no-$(k+1)$-in line problem,
\url{arxiv:2508.07632} (2025).

\bibitem{Kovacs-2025} 
B.~Kov\'{a}cs, Z.L.~Nagy, D.R.~Szab\'{o},
Settling the no-$(k+1)$-in-line problem when $k$ is not small, 
\url{arXiv:2502.00176} (2025).

\bibitem{kovic-2023}
J.~Kovi\v{c}, T.~Pisanski, S.S.~Zemlji\v{c}, A.~\v{Z}itnik, 
The {S}ierpi\'nski product of graphs,
Ars Math.\ Contemp.\ 23 (2023) Paper 1.

\bibitem{Kruft-2024}
E.~Kruft Welton, S.~Khudairi, J.~Tuite,
Lower general position in Cartesian products,
Commun.\ Comb.\ Optim.\  10 (2025) 110--125.

\bibitem{Ku} 
C.Y.~Ku, K.B.~Wong, 
On no-three-in-line problem on $m$-dimensional torus, 
Graphs Combin.\ 34  (2018) 355--364.

\bibitem{RamseyTheoryontheIntegers}
B.M.~Landman, A.~Robertson,
Ramsey Theory on the Integers,
Vol. 73. American Mathematical Soc.\ (2014).

\bibitem{li-2025} 
Z.-L.~Li, S.-C.~Gong,
Graphs whose edge general position number is $4$,
Bull.\ Malays.\ Math.\ Sci.\ Soc.\ 48 (2025) Paper 87.

\bibitem{ManuelKlavzar-2018a}
P.~Manuel, S.~Klav\v{z}ar,
A general position problem in graph theory,
Bull.\ Aust.\ Math.\ Soc.\ 98 (2018) 177--187.

\bibitem{manuel-2017arxiv} 
P.~Manuel, S.~Klav\v{z}ar, 
Graph theory general position problem,
\url{arXiv:1708.09130} (2017).

\bibitem{ManuelKlavzar-2018b}
P.~Manuel, S.~Klav{\v z}ar,
The graph theory general position problem on some interconnection networks,
Fund.\ Inform.\ 163 (2018) 339--350.

\bibitem{Manuel-2022}
P.~Manuel, R.~Prabha, S.~Klav\v{z}ar,
The edge general position problem,
Bull.\ Malays.\ Math.\ Sci.\ Soc.\ 45 (2022) 2997--3009.

\bibitem{ManuelPrabhaKlavzar-2023+}
P.~Manuel, R.~Prabha, S.~Klav\v{z}ar,
Generalization of edge general position problem,
Art Discrete Appl.\ Math.\ 8 (2025) \#P1.02.

\bibitem{mertens-2024}
S.~Mertens, 
Domination polynomial of the rook graph,
J.\ Integer Seq.\ 27 (2024) Paper 24.3.7. 

\bibitem{Misiak} 
A.~Misiak, Z.~St\k{e}pie\'{n}, A.~Szymaszkiewicz, L.~Szymaszkiewicz, M.~Zwierzchowski, 
A note on the no-three-in-line problem on a torus, 
Discrete Math.\ 339 (2016) 217--221.

\bibitem{Nagy-2023} 
D.T.~Nagy, Z.L.~Nagy, R.~Woodroofe, 
The extensible No-Three-In-Line problem, 
European J.\ Combin.\ 114 (2023) Paper 103796.

\bibitem{Neethu-2021}
P.K.~Neethu, U.~Chandran S.V., M.~Changat, S.~Klav{\v z}ar,
On the general position number of complementary prisms,
Fund.\ Inform.\ 178 (2021) 267--281.

\bibitem{Oellermann2007}
O.R.~Oellermann, J.~Peters-Fransen,
The strong metric dimension of graphs and digraphs,
Discrete Appl.\ Math.\ 155 (2017) 356--364.

\bibitem{Oh}
S.~Oh, J.R.~Schmitt, X.~Wang, 
Repeatedly applying the Combinatorial Nullstellensatz for Zero-sum Grids to Martin Gardner’s minimum no-3-in-a-line problem, 
European J.\ Combin.\ 125 (2025) Paper 104095.

\bibitem{OEIS} 
https://oeis.org/search?q=A219760\&language=english\&go=Search.

\bibitem{Patkos-2019}
B.~Patk\'{o}s, 
On the general position problem on Kneser graphs,
Ars Math.\ Contemp.\ 18 (2020) 273--280.

\bibitem{Paynethesis}
M.S.~Payne,
Combinatorial Geometry of Point Sets with Collinearities,
PhD Thesis, University of Melbourne (2014).

\bibitem{Payne} 
M.S.~Payne, D.R.~Wood, 
On the general position subset selection problem, 
SIAM J.\ Discrete Math. 27 (2013) 1727--1733.

\bibitem{Pelayo-2013} 
I.M.~Pelayo, Geodesic Convexity in Graphs, 
Springer, New York, 2013.

\bibitem{Por} 
A.~P\'{o}r, D.R.~Wood, 
No-three-in-line-in-3D, 
Algorithmica 47 (2007) 481--488.

\bibitem{Prabha-2020}
R.~Prabha, S.~Renukaa Devi, 
General position problem of hyper tree and shuffle hyper tree networks, 
AIP Conf.\ Proc.\ 2277 (2020) Paper 100015. 

\bibitem{Prabha-2021}
R.~Prabha,  S.~Renukaa Devi,
General position problem of hexagonal derived networks,
Adv.\ Appl.\ Math.\ Sci.\ 21 (2021) 145--154.

\bibitem{Prabha-2023} 
R.~Prabha, S.~Renukaa Devi, P.~Manuel, 
General position problem of butterfly networks, 
\url{arXiv:2302.06154} (2023).

\bibitem{Prellberg2026}
T.~Prellberg, 
Constraint satisfaction programming for the no-three-in-line problem,
J.\ Combin.\ Theory 225 (2027) 106244.

\bibitem{Prellbergcheckerboard}
T.~Prellberg, 
No-three-in-line sets on the checkerboard grid,
\url{arXiv:2605.09215} (2026).

\bibitem{Prellberg} 
T.~Prellberg, 
\url{https://oeis.org/A272651/a272651_13.txt} (accessed 21/01/26).

\bibitem{Ramanathan2025}
P.~Ramanathan, T.~Prellberg, M.~Lewis, P.D.~Joshi, R.A.~Dandekar, R.~Dandekar, S.~Panat,
Three methods, one problem: Classical and AI approaches to no-three-in-line,
\url{arXiv:2512.11469} (2025).

\bibitem{randriambololona-2013}
H.~Randriambololona,
$(2,1)$-separating systems beyond the probabilistic bound,
Israel J.\ Math.\ 195 (2013) 171--186.

\bibitem{rather-2026+}
B.A.~Rather,
Explicit formulas and unimodality phenomena for general position polynomials,
\url{arXiv:2603.06930} (2026).

\bibitem{Reijnders} 
L.~Reijnders, 
Eigenvalue Bounds for Distance-type Graph Parameters, 
MSc Thesis Eindhoven University of Technology (2023).

\bibitem{Rodriguez-2022}
J.A.~Rodr\'{\i}guez-Vel\'{a}zquez,
Universal lines in graphs,
Quaest.\ Math.\ 45 (2022) 1485--1500.

\bibitem{Roth} 
K.F.~Roth, On a problem of Heilbronn, 
J.\ Lond.\ Math.\ Soc.\ 26 (1951) 198--204.

\bibitem{Roy-2025}
D.~Roy, S.~Klav\v{z}ar, A.S.~Lakshmanan,
Mutual-visibility and general position in double graphs and in Mycielskians,
Appl.\ Math.\ Comput.\ 488 (2025) Paper 129131.

\bibitem{Roy-2026}
D.~Roy, S.~Klav\v{z}ar, A.S.~Lakshmanan, J.~Tian,
Varieties of mutual-visibility and general position on Sierpi\'nski graphs,
Discuss.\ Math.\ Graph Theory (2026) \url{doi.org/10.7151/dmgt.2625}.

\bibitem{schaefer-1978} 
T.~Schaefer, 
On the complexity of some two-person perfect-information games,
J.\ Comp.\ Sys.\ Sci.\ 16 (1978) 185–-225.

\bibitem{Shallcross} 
E.~Shallcross, 
personal communcation (2024).

\bibitem{Shallcross-2025} 
E.~Shallcross, J.~Tuite, A.~Evans, A.~Krishnakumar, S.~Boshar,
Solution to some conjectures on mobile position problems,
\url{arXiv:2507.16622} (2025).

\bibitem{Suk-spheres}
A.~Suk, E.P.~White, A note on the no-$(d + 2)$-on-a-sphere problem, Vol. 332 Leibniz International Proceedings in Informatics (2025) 76:1--76:8.
 
\bibitem{Suk} 
A.~Suk, J.~Zeng, 
On higher dimensional point sets in general position, 
Comb.\ Prob.\ Comput.\ 35 (2026) 134--148.

\bibitem{Szabo}
D.R.~Szab\'{o}, 
Rational normal curves as no-$(d+2)$-on-$Q$-quadric sets, \url{arXiv:2511.03526} (2025).

\bibitem{ThaChaTuiThoSteErs-2024} 
M.~Thankachy, U.~Chandran S.V., J.~Tuite, E.~Thomas, G.~Di Stefano, G.~Erskine, 
On the vertex position number of graphs, 
Discuss.\ Math.\ Graph Theory 44 (2024) 1169--1188.  

\bibitem{Thiele}
T.~Thiele,
Geometric Selection Problems and Hypergraphs,
PhD Thesis,  Freie Universit\"{a}t Berlin (1995).

\bibitem{Thiele2}
T.~Thiele, 
The no-four-on-circle problem,
J.\ Comb.\ Theory Ser.\ A 71 (1995) 332--334.

\bibitem{Thomas-2020}
E.J.~Thomas, U.~Chandran S.V.,
Characterization of classes of graphs with large general position number,
AKCE Int.\ J.\ Graphs Comb.\ (2020) 1--5.

\bibitem{Thomas-2021}
E.J.~Thomas, U.~Chandran S.V.,
On independent position sets in graphs,
Proyecciones 40 (2021) 385--398. 

\bibitem{Thomas-2024a}
E.J.~Thomas, U.~Chandran S.V., J.~Tuite, G.~Di Stefano,
On the general position number of Mycielskian graphs,
Discrete Appl.\ Math.\ 353 (2024) 29--43. 

\bibitem{Thomas-2024b}
E.J.~Thomas, U.~Chandran S.V., J.~Tuite, G.~Di Stefano,
On monophonic position sets in graphs,
Discrete Appl.\ Math.\ 354 (2024) 72--82. 

\bibitem{Tian-2026}
J.~Tian, P. Dokyeesun, S.~Klav\v{z}ar,
On the variety of general position problems under vertex and edge removal, 
 Discrete Appl.\ Math.\ 388 (2026) 56--64. 

\bibitem{TianKlavzar-2025}
J.~Tian, S.~Klav\v{z}ar,
Variety of general position problems in graphs,
Bull.\ Malays.\ Math.\ Sci.\ Soc.\ 48 (2025) Paper 5. 

\bibitem{TianKlavzar-2024++}
J.~Tian, S.~Klav\v{z}ar,
General position sets, colinear sets, and Sierpi\'nski product graphs,
Ann.\ Comb.\ (2025) 29 (2025) 837--852.

\bibitem{Tian-2024b}
J.~Tian, S.~Klav\v{z}ar, E.~Tan,
Extremal edge general position sets in some graphs,
Graphs Combin.\ 40 (2024) Paper 40.

\bibitem{Tian-2021a}
J.~Tian, K.~Xu,
The general position number of Cartesian products involving a factor with small diameter,
Appl.\ Math.\ Comput.\ 403 (2021) Paper 126206.

\bibitem{powers} 
J.~Tian, K.~Xu, 
On the general position number of the $k$-th power graphs, 
Quaest.\ Math.\ 47 (2024) 2215--2230.

\bibitem{TianXuChao-2023}
J.~Tian, K.~Xu, D.~Chao,
On the general position numbers of maximal outerplane graphs,
Bull.\ Malays.\ Math.\ Sci.\ Soc.\ 46 (2023) Paper 198.

\bibitem{Tian-2021b}
J.~Tian, K.~Xu, S.~Klav\v{z}ar,
The general position number of the Cartesian product of two trees,
Bull.\ Aust.\ Math.\ Soc.\ 104 (2021) 1--10.

\bibitem{Tsuchiya95}
Y.~Tsuchiya, Y.~Takefuji, 
A neural network algorithm for the no-three-in-line problem, 
Neurocomputing 325 (1995) 85--98.

\bibitem{TuiThoCha-2022} 
J.~Tuite, E.~Thomas, U.~Chandran S.V., 
Some position problems for graphs,
Lecture Notes in Comput.\ Sci.\ 13179 (2022) 36--47. 

\bibitem{TuiThoCha-2025} 
J.~Tuite, E.~Thomas, U.~Chandran S.V., 
On some extremal position problems for graphs, 
Ars Math.\ Contemp.\ 25 (2025) \#P1.09. 

\bibitem{Vesel}
A.~Vesel,
General position sets in strong products with paths and cycles,
\url{arXiv:2607.26844} (2026).

\bibitem{Voutier}
P.M.~Voutier,
On the Guy-Kelly Conjecture for the no-three-in-line problem,
\url{arxiv:2603.00215} (2026).

\bibitem{Wang2026}
X.~Wang,
A complete multipartite counterexample to corona-unimodality of general position polynomials,
\url{https://zenodo.org/records/20411641} (2026).

\bibitem{Wang2026-2}
X.~Wang,
Complete bipartite counterexamples and the rational-slope phase diagram for corona-unimodality of general position polynomials,
\url{https://zenodo.org/records/20453868} (2026).

\bibitem{west-2021}
D.B.~West,
Combinatorial Mathematics,
Cambridge University Press, Cambridge, 2021.

\bibitem{Woodcolouring} 
D.R.~Wood, 
A note on colouring the plane grid, 
Geombinatorics 13 (2004) 193--196.

\bibitem{wooddrawing} 
D.R.~Wood, 
Grid drawings of $k$-colourable graphs, 
Comp.\ Geom.\ Theory Appl.\ 30 (2005) 25--28.

\bibitem{yao-2022}
Y.~Yao, M.~He, S.~Ji,
On the general position number of two classes of graphs,
Open Math.\ 20 (2022) 1021--1029.

\bibitem{Zarankiewicz-1951}
K.~Zarankiewicz,
Problem P 101,
Colloq.\ Math.\ 2 (1951) 301.

\bibitem{Zaslavsky}
C.~Zaslavsky, 
Tic Tac Toe: And Other Three-In-A Row Games from Ancient Egypt to the Modern Computer, 
Crowell, 1982.

\end{thebibliography}
\end{document}